\documentclass[preprint]{imsart}

\RequirePackage[OT1]{fontenc}
\RequirePackage{amsthm,amsmath}
\RequirePackage{natbib}
\RequirePackage[colorlinks,citecolor=blue,urlcolor=blue]{hyperref}
\startlocaldefs
    \numberwithin{equation}{section}
    \theoremstyle{plain}
    
    \usepackage{mathtools}
    \usepackage{amsmath, amsthm, amssymb, dsfont, hyperref}
    \usepackage{natbib}
    \usepackage{notation}
    \usepackage{subcaption}
\endlocaldefs

\begin{document}

\begin{frontmatter}
    \title{Minimax semi-supervised confidence sets for multi-class classification}
    %Optimal rates for semi-supervised confidence set for multi-class classification}
    \runtitle{Confidence Sets}
% \thankstext{T1}{Footnote to the title with the ``thankstext'' command.}

    \begin{aug}
        \author{\fnms{Evgenii} \snm{Chzhen}\thanksref{t1}\ead[label=e1]{evgenii.chzhen@u-pem.fr}},
        \author{\fnms{Christophe} \snm{Denis}\ead[label=e2]{christophe.denis@u-pem.fr}}
        \and
        \author{\fnms{Mohamed} \snm{Hebiri}
        \ead[label=e3]{mohamed.hebiri@u-pem.fr}}

        \thankstext{t1}{This work was partially supported by ``Labex B\'ezout'' of Universit\'e Paris-Est}
        % \thankstext{t2}{First supporter of the project}
        % \thankstext{t3}{Second supporter of the project}
        \runauthor{Chzhen, Denis and Hebiri}

        \affiliation{LAMA, Universit\'e Paris-Est -- Marne-la-Vall\'ee}

        \address{Universit\'e Paris-Est -- Marne-la-Vall\'ee\\
        Cit\'e Descartes,
        Bâtiment Copernic\\
        5 boulevard Descartes\\
        77454 Marne-la-Vall\'ee cedex 2\\
        \printead{e1}\\
        \phantom{E-mail:\ }\printead*{e2}\\
        \phantom{E-mail:\ }\printead*{e3}}

        % \address{Address of the Third author\\
        % Usually a few lines long\\
        % Usually a few lines long\\
        % \printead{e3}\\
        % \printead{u1}}
    \end{aug}

    \begin{abstract}
    In this work we study the semi-supervised framework of confidence set classification with controlled expected size in minimax settings.
    We obtain semi-supervised minimax rates of convergence under the margin assumption and a H\"older condition on the regression function.
    Besides, we show that if no further assumptions are made, there is no supervised method that outperforms the semi-supervised estimator proposed in this work.
    We establish that the best achievable rate for any supervised method is $n^{-1/2}$, even if the margin assumption is extremely favorable.
    On the contrary, semi-supervised estimators can achieve faster rates of convergence provided that sufficiently many unlabeled samples are available.
    We additionally perform numerical evaluation of the proposed algorithms empirically confirming our theoretical findings.
    \end{abstract}



\end{frontmatter}

%%%%%%%%%%%%%%%%%%%%%%%%%%%%%%%%%%%%%%%%%%%%%%%%%%%%%%%%%%%%%%%%%%%%%%%%%%%%%%%
    %%%%%%%%%%%%%%%%%%%%%%%%%%%%%%%%%%%%%%%%%%%%%%%%%%%%%%%%%%%%%%%%%%%%%%%%%%%%%%%
%%%%%%%%%%%%%%%%%%%%%%%%%%%%%%%%%%%%%%%%%%%%%%%%%%%%%%%%%%%%%%%%%%%%%%%%%%%%%%%
\section{Introduction}
\label{sec:problem}
%%%%%%%%%%%%%%%%%%%%%%%%%%%%%%%%%%%%%%%%%%%%%%%%%%%%%%%%%%%%%%%%%%%%%%%%%%%%%%%
%%%%%%%%%%%%%%%%%%%%%%%%%%%%%%%%%%%%%%%%%%%%%%%%%%%%%%%%%%%%%%%%%%%%%%%%%%%%%%%

Let $K \geq 2$ and $(X ,Y) \in \bbR^d \times [K] \eqdef \ens{1, \ldots, K}$ be a random couple distributed according to a distribution $\Prob$ on $\bbR^d \times [K]$, where $X \in \bbR^d$ is seen as the feature vector and $Y \in [K]$ as the class.
This problem falls within the scope of the multi-class setting where the goal is to predict the label $Y$ for a given feature.
Commonly, prediction is performed by a classifier that outputs a single label.
However, in the confidence set framework, the objective differs: we aim at predicting a \emph{set} of labels instead of a \emph{single} one.
This problem has been studied in a few works, and we consider in this contribution the setup put forward by~\cite{Denis_Hebiri17}.
The essential feature of their perspective is the control of the size of confidence sets in expectation.
While they provided a procedure to build confidence sets based on Empirical Risk Minimization (ERM) and established upper bounds, the present work aims at giving a general analysis of the confidence problem in the minimax sense.

%%%%%%%%%%%%%%%%%%%%%%%%%%%%%%%%%%%%%%%%%%%%%%%%%%%%%%%%%%%%%%%%%%%%%%%%%%%%%%%
%%%%%%%%%%%%%%%%%%%%%%%%%%%%%%%%%%%%%%%%%%%%%%%%%%%%%%%%%%%%%%%%%%%%%%%%%%%%%%%
\subsection{Problem statement}
\label{subsec: Theproblem}
%%%%%%%%%%%%%%%%%%%%%%%%%%%%%%%%%%%%%%%%%%%%%%%%%%%%%%%%%%%%%%%%%%%%%%%%%%%%%%%
%%%%%%%%%%%%%%%%%%%%%%%%%%%%%%%%%%%%%%%%%%%%%%%%%%%%%%%%%%%%%%%%%%%%%%%%%%%%%%%

All along the paper, we denote by $\Prob_X$ the marginal distribution of $X \in \bbR^d$ and by $p(\cdot) \eqdef (p_1(\cdot), \ldots, p_K(\cdot))^\top$ the regression function defined for all $k \in [K]$ and all $x \in \bbR^d$ as $p_k(x) \eqdef \Prob(Y = 1 | X = x)$.
For any sets $A, A' \subset [K]$ we denote by $A \triangle A'$ their symmetric difference.
We assume that two data samples $\data_n, \data_N$ are available.
The first sample $\data_n = \{(X_i, Y_i)\}_{i = 1}^n$ consists of $n \in \bbN$ \iid copies of $(X, Y) \in \bbR^d \times [K]$ and the second sample $\data_N = \{X_i\}_{i = n + 1}^{n + N}$ consist of $N \in \bbN$ \iid copies of $X \in \bbR^d$.

A confidence set classifier $\Gamma$ is a measurable function from $\bbR^d$ to $2^{[K]} \eqdef \enscond{A}{A \subset [K]}$, that is, $\Gamma : \bbR^d \rightarrow 2^{[K]}$ and we denote by $\Upsilon$ the set of all such functions.
For any confidence set $\Gamma : \bbR^d \rightarrow 2^{[K]}$ we define its error and its information as
\begin{align*}
    \underbrace{\perf\pare{\Gamma} = \Prob\pare{Y \notin \Gamma(X)}}_{\text{error}}\,,\qquad\qquad
    \underbrace{\info\pare{\Gamma} = \Exp_{\Prob_X}\abs{\Gamma(X)}}_{\text{information}}\enspace,
\end{align*}
respectively, where $\Exp_{\Prob_X}$ stands for the expectation \wrt the marginal distribution of $X \in \bbR^d$ and $\abs{\Gamma(x)}$ is the cardinal of $\Gamma$ at $x\in \bbR^d$.

For a fixed integer $\beta \in [K]$ a $\beta$-Oracle confidence set $\Gamma^*_\beta$ is defined as
\begin{align*}
    \Gamma^*_\beta \in \argmin\enscond{\perf\pare{\Gamma}}{\Gamma \in \Upsilon \text { \st } \info(\Gamma) = \beta}\enspace.
\end{align*}
The set $\enscond{\Gamma \in \Upsilon}{\info(\Gamma) = \beta}$ is always non-empty, as it always contain those confidence sets whose cardinal is equals to $\beta$ for every $x \in \bbR^d$.

The description of $\beta$-Oracle confidence set in general situation might be complicated.
Hence, we introduce the following mild assumption, which allows to obtain an explicit expression.
\begin{ass}[Continuity of CDF]
    \label{ass:continuity_cdf}
    For all $k \in [K]$ the cumulative distribution function (CDF) $F_{p_k}(\cdot) \eqdef \Prob_X(p_k(X) \leq \cdot)$ of $p_k(X)$ is continuous on $(0, 1)$.
\end{ass}
\begin{prop}[$\beta$-Oracle confidence set]
    \label{prop:oracle_set}
    Fix $\beta \in [K - 1]$, and let the function $G: [0, 1] \rightarrow [0, K]$ be defined for all $t \in [0, 1]$ as
    \begin{align*}
        G(t) \eqdef \sum_{k = 1}^K \pare{1 - {F}_{p_k}(t)} = \sum_{k = 1}^K \Prob_X(p_k(X) > t)\enspace,
    \end{align*}
    then under Assumption~\ref{ass:continuity_cdf} a $\beta$-Oracle confidence set $\Gamma^*_\beta$ can be obtained as
    \begin{align}
        \label{eq:bayesConfSet}
        \Gamma_\beta^* (x)  =  \enscond{k \in [K]}{p_{k}(x)\geq G^{-1}(\beta)}\enspace,
    \end{align}
    where we denote by $G^{-1}$ the generalized inverse of $G$ defined for all $\beta \in [0, K]$ as
    \begin{align*}
        G^{-1}(\beta) \eqdef \inf\enscond{t\in [0,1]}{G(t) \leq \beta} \enspace.
    \end{align*}
\end{prop}
\begin{prop}
    \label{prop:minimizer_of_penalized_risk}
    Assume that Assumption~\ref{ass:continuity_cdf} is fulfilled, then the $\beta$-Oracle defined in Eq.~\eqref{eq:bayesConfSet} is a minimizer of the following risk
    \begin{align}
        \label{eq:penalized_risk}
        \risk_\beta(\Gamma) = \perf(\Gamma) + G^{-1}(\beta)\info(\Gamma)\enspace.
    \end{align}
\end{prop}
These propositions have been proven in~\citep[Proposition~4 and Proposition~7]{Denis_Hebiri17}.
Consequently, the accuracy of a confidence set $\Gamma$ can be for instance quantified according its excess risk
\begin{equation*}
\risk_\beta(\Gamma) - \risk_\beta(\Gamma^*_\beta) = \sum_{k =1}^K \Exp_{\Prob_X}\left[|p_k(X)-G^{-1}(\beta)|\1_{\{k \in \Gamma(X) \triangle \Gamma^{*}_{\beta}(X)\}}\right].
\end{equation*}
The statistical learning problem is then to estimate $\Gamma_{\beta}^*$ given the data sample $\data_n$ and $\data_N$.
The formulation in Eq.~\eqref{eq:bayesConfSet} of the $\beta$-Oracle appears to be closely related to the level set estimation problem~\citep{Hartigan87, Polonik95, Tsybakov97,Rigollet_Vert09}.
Hence at first sight, the introduction of an unlabeled sample may be surprising.
However, in our setup the estimation of the $\beta$-Oracle does not only rely on the regression function but also on the threshold $G^{-1}(\beta)$ which is \emph{unknown} beforehand and can be estimated in a semi-supervised way~\citep{Denis_Hebiri17}.
To fix these ideas, we give some examples of possible estimation procedures of $\Gamma_{\beta}^{*}$.
%%%%%%%%%%%%%%%%%%%%%%%%%%%%%%%%%%%%%%%%%%%%%%%%%%%%%%%%%%%%%%%%%%%%%%%%%%%%%%%
%%%%%%%%%%%%%%%%%%%%%%%%%%%%%%%%%%%%%%%%%%%%%%%%%%%%%%%%%%%%%%%%%%%%%%%%%%%%%%%
\subsection{Confidence set estimators}
\label{subsec: confidenceEStim}
%%%%%%%%%%%%%%%%%%%%%%%%%%%%%%%%%%%%%%%%%%%%%%%%%%%%%%%%%%%%%%%%%%%%%%%%%%%%%%%
%%%%%%%%%%%%%%%%%%%%%%%%%%%%%%%%%%%%%%%%%%%%%%%%%%%%%%%%%%%%%%%%%%%%%%%%%%%%%%%
An estimator $\hGamma$ is a measurable function that maps any given data samples into a confidence set classifier.
We shall distinguish two types of estimators: supervised and semi-supervised whose formal definition is provided below.
\begin{definition}[Supervised and semi-supervised estimators]
    \label{def:supervised_semi_supervised}
    A measurable mapping
    \begin{align*}
    \hGamma : \bigcup_{n, N \in \bbN}\pare{\bbR^d \times [K]}^n \times \pare{\bbR^d}^N \rightarrow \Upsilon\enspace,
    \end{align*}
    is called a \emph{supervised} estimator if for any $n, N \in \bbN$ and any data samples $\data_n = \{(X_i, Y_i)\}_{i = 1}^n$, $\data_N = \{X_i\}_{i = n + 1}^{n + N}$, and $\data'_N = \{X'_i\}_{i = n + 1}^{n + N}$ it holds that
    \begin{align*}
        \hGamma(x; \data_n, \data_N) = \hGamma(x; \data_n, \data'_N), \quad \text{ a.e. } x \in \bbR^d \text{ w.r.t. the Lebesgue measure}  \enspace.
    \end{align*}
    Otherwise the estimator is called \emph{semi-supervised}.
    In the sequel, for the simplicity of notation we write $\hGamma(x)$ instead of $\hGamma(x; \data_n, \data_N)$ where no ambiguity is present.
\end{definition}
Intuitively, the supervised estimators do not take into account the information that is provided by the unlabeled sample.
Besides, if we denote by $\hat{\Upsilon}$ the set of all estimators, Definition~\ref{def:supervised_semi_supervised} generates a natural partition of $\hat\Upsilon$ into two disjoint sets: the \emph{supervised} estimators $\hat\Upsilon_{\SE}$ and the \emph{semi-supervised} estimators $\hat\Upsilon_{\SSE}$.

Hereafter, we provide three different examples of estimation procedures which are the core of our study. All these methods rely on
{\it plug-in} principle.
\begin{itemize}
    \item {\it Top-$\beta$ procedure.}
    This method is the most intuitive estimator in the considered context.
    It is a supervised procedure, that is, based only on $\data_n$.
    Let consider an estimator $\hat{p}$ of the regression function $p$.
    Let $\left(\hat{p}_{\sigma_k(X)}\right)_{k \in [K]}$ be the order statistic associated to $\hat{p}(X)$, such that for all $x\in\bbR^d$ we have $\hat{p}_{\sigma_1(x)}(x) \geq \ldots \geq \hat{p}_{\sigma_K(x)}(x)$.
    A top-$\beta$ confidence set is then defined as
    \begin{equation}
        \label{eq:IntroTopBetaEstimator}
        \hGamma_{{\rm top}}(x) = \left\{\sigma_1 (x), \ldots, \sigma_{\beta}(x)  \right\},\qquad \forall x\in \bbR^d \enspace.
    \end{equation}
    \item {\it Supervised procedure.} Formally, in this type of methods, we only care about $\data_n$ (we forget about $\data_N$).
    We split $\data_n$ into two independent samples such that  $\data_n = \data_{\floor{n / 2}} \bigcup \data_{\ceil{n / 2}}$.
    Based on the first sample $\data_{\floor{n / 2}}$, we consider an estimator $\hat{p}$ of the regression function $p$.
    Furthermore, we define
    \begin{equation*}
        \hat{G}(\cdot) = \frac{1}{\ceil{n / 2}} \sum_{i \in \data_{\ceil{n / 2}}} \sum_{k = 1}^K \1_{\{\hat{p}_k(X_i) \geq \cdot\}}\enspace,
    \end{equation*}
    and one type of supervised estimator is then defined as follows
    \begin{equation}
        \label{eq:IntroSupEstimator}
        \hGamma_{\SE}(x) = \left\{ k \in [K] : \ \hat{p}_{k}(x) \geq \hat{G}^{-1}(\beta)  \right\},\qquad \forall x\in \bbR^d \enspace .
    \end{equation}
    Interestingly, conditional on the data sample $\data_{\floor{n / 2}}$, the definition of the estimator $\hat{G}$
    does not involves the labels associated to $\data_{\ceil{n / 2}}$. As a consequence, we can naturally consider a semi-supervised version
    of this estimator.
    \item {\it Semi-supervised procedure.}
    Based on $\data_n$, we consider an estimator $\hat{p}$  of the regression function $p$.
    Furthermore, we define
    \begin{equation*}
    \hat{G}(\cdot) = \frac{1}{N} \sum_{i \in  \data_N}\sum_{k = 1}^K \1_{\{\hat{p}_k(X_i) \geq \cdot\}}\enspace,
    \end{equation*}
    and one type of semi-supervised estimator is then defined as follows
    \begin{equation}
        \label{eq:IntroSemiSupEstimator}
        \hGamma_{\SSE}(x) = \left\{ k \in [K] : \hat{p}_{k}(x) \geq \hat{G}^{-1}(\beta)  \right\},\qquad \forall x\in \bbR^d\enspace.
    \end{equation}
\end{itemize}
One can note that these procedures are based on a preliminary estimator of
$p$ built from $\data_n$, that is, all of them are plug-in type procedures.
However, these procedures differ by the construction of the output set.
The top-$\beta$ procedure and the supervised procedure rely only on the labeled data while the semi-supervised estimator takes advantage of the information provided by the unlabeled data.
The top-$\beta$ procedure is the simplest among them, it naturally satisfies $|\hGamma(x)| = \beta$ for all $x \in \bbR^d$.
At the same time, the others are more involved and can have different cardinals for different values of $x \in \bbR^d$.
Nevertheless, for the other two procedures one can guarantee $\info(\hGamma) \approx \beta$.

These examples give a rise to natural questions which form the core our theoretical study and which are summarized below.
% Bellow we summarize these questions our which form the core of our study  we summarize below and which are in the core of our theoretical study.
\begin{enumerate}
    \item The first question is the statistical performance of these plug-in procedures which is assessed through rates of convergence and their optimality in the minimax sense.
    \item The second question focuses on the benefit of the semi-supervised approach.
    Roughly speaking, are there situations where the semi-supervised approach outperforms the supervised one and how can it be quantified?
    \item The third question concentrates on the reason why it is more relevant for this problem to consider more involved estimators than the simple top-$\beta$ method.

    % Finally as suggested in~\citep{Denis_Hebiri17}, why it is more relevant for this problem to consider estimators more involved than the simple top-$\beta$ method?
\end{enumerate}

%%%%%%%%%%%%%%%%%%%%%%%%%%%%%%%%%%%%%%%%%%%%%%%%%%%%%%%%%%%%%%%%%%%%%%%%%%%%%%%
%%%%%%%%%%%%%%%%%%%%%%%%%%%%%%%%%%%%%%%%%%%%%%%%%%%%%%%%%%%%%%%%%%%%%%%%%%%%%%%
\subsection{Minimax estimation}
\label{subsec:minimax_estimation}
%%%%%%%%%%%%%%%%%%%%%%%%%%%%%%%%%%%%%%%%%%%%%%%%%%%%%%%%%%%%%%%%%%%%%%%%%%%%%%%
%%%%%%%%%%%%%%%%%%%%%%%%%%%%%%%%%%%%%%%%%%%%%%%%%%%%%%%%%%%%%%%%%%%%%%%%%%%%%%%
For a given family $\class{P}$ of joint distributions on $\bbR^d \times [K]$, a given estimator $\hGamma \in \hat\Upsilon$, and fixed integers $K \geq 2$,  $\beta \in [K]$, $n, N \in \bbN$ we are interested in the following maximal risks
\begin{align*}
    \excess^{\Ham}_{n, N}(\hGamma;\class{P})
    &\eqdef \sup_{\Prob \in \class{P}} \Exp_{(\data_n, \data_N)}\Exp_{\Prob_X}\abs{\hGamma(X) \triangle \Gamma^*_\beta(X)} & \text{(Hamming risk)}\enspace,\\
    \excess^{\Exc}_{n, N}(\hGamma; \class{P})
    &\eqdef
    \sup_{\Prob \in \class{P}} \Exp_{(\data_n, \data_N)}\risk_\beta(\hGamma) - \risk_\beta(\Gamma^*_\beta) & \text{(Excess risk)}\enspace,\\
    \excess_{n, N}^{\Disc}(\hGamma; \class{P})
    &\eqdef
    \sup_{\Prob \in \class{P}} \Exp_{(\data_n, \data_N)}\left[\abs{\perf(\hGamma) - \perf(\Gamma^*_\beta)} +\abs{\beta - \info(\hGamma)}\right]& \text{(Discrepancy)}\enspace,
\end{align*}
where $\Exp_{(\data_n, \data_N)}$ denotes the expectation \wrt $\Prob^{\otimes n} \otimes \Prob_X^{\otimes N}$.
These maximal risks are arising in a natural way in the context of the confidence set estimation with controlled expected size.
The risk $\excess^{\Ham}_{n, N}(\hGamma; \class{P})$ corresponds to the estimation of the $\beta$-Oracle through the Hamming distance.
The second risks is directly connected with Proposition~\ref{prop:oracle_set}, which gives a description of the $\beta$-Oracle as a minimizer of $\risk_\beta(\cdot)$.
As the goal in this problem is to construct a procedure $\hGamma$ that exhibits a low error $\perf(\hGamma)$ and low cardinal discrepancy $\absin{\beta - \info(\hGamma)}$, it is natural to consider $\excess^{\Disc}_{n, N}(\hGamma; \class{P})$ which is composed of both.

Finally, we are in position to define the notion of the minimax rate.
% As one of our goals is to quantify the benefits of the semi-supervised estimators over the supervised estimators, the notion of minimax rate slightly differs from the conventional one.
% In fact,
The minimax rate in this context is not only determined by the family of distributions $\class{P}$ but also by the family of estimators $\classGamma \subset \hat\Upsilon$ that we consider.
\begin{definition}[Minimax rate of convergence]
    \label{def:minimax_optimal_rate}
    For a given family $\class{P}$ of joint distributions on $\bbR^d \times [K]$ and a given family of estimators $\classGamma \subset \hat\Upsilon$ the minimax rates are defined as
    \begin{align*}
        \excess^{\square}_{n, N}(\classGamma;\class{P})
        \eqdef
        \inf_{\hGamma \in \classGamma} \excess^{\square}_{n, N}(\hGamma;\class{P})\enspace,
    \end{align*}
    where $\square$ is $\Ham$, $\Exc$ or $\Disc$.
\end{definition}
The main families of estimators that we study are the \emph{supervised} $\hat\Upsilon_{\SE}$ and the \emph{semi-supervised} $\hat\Upsilon_{\SSE}$ estimators.
Obviously, since $\hat\Upsilon = \hat\Upsilon_{\SE} \bigcup \hat\Upsilon_{\SSE}$ and $\hat\Upsilon_{\SE} \bigcap \hat\Upsilon_{\SSE} = \emptyset$, we have the following relation
\begin{align*}
   \excess^{\square}_{n, N}(\hat\Upsilon;\class{P})
   =
   \excess^{\square}_{n, N}(\hat\Upsilon_{\SE};\class{P}) \bigwedge \excess^{\square}_{n, N}(\hat\Upsilon_{\SSE};\class{P})\enspace.
\end{align*}
As a consequence, a lower and an upper bounds on $\excess^{\square}_{n, N}(\hat\Upsilon_{\SE};\class{P})$, $\excess^{\square}_{n, N}(\hat\Upsilon_{\SSE};\class{P})$ yield the bounds on the minimax rate over all estimators.

%%%%%%%%%%%%%%%%%%%%%%%%%%%%%%%%%%%%%%%%%%%%%%%%%%%%%%%%%%%%%%%%%%%%%%%%%%%%%%%
%%%%%%%%%%%%%%%%%%%%%%%%%%%%%%%%%%%%%%%%%%%%%%%%%%%%%%%%%%%%%%%%%%%%%%%%%%%%%%%
\subsection{Related works}
\label{subsec:related_works}
%%%%%%%%%%%%%%%%%%%%%%%%%%%%%%%%%%%%%%%%%%%%%%%%%%%%%%%%%%%%%%%%%%%%%%%%%%%%%%%
%%%%%%%%%%%%%%%%%%%%%%%%%%%%%%%%%%%%%%%%%%%%%%%%%%%%%%%%%%%%%%%%%%%%%%%%%%%%%%%
Confidence set approach for classification was pioneered by \cite{Vovk02_transductive,Vovk02IndepError,Vovk_Gammerman_Shafer05} by the means of conformal prediction theory.
They rely on non-conformity measures which are based on some pattern recognition methods, and develop an asymptotic theory.
In this work, we consider a statistical perspective of confidence set classification and put our focus on non-asymptotic minimax theory.

The problem of confidence set multi-class classification has strong ties with the binary classification with reject option, also known as binary classification with abstention in machine learning literature.
In the binary classification with rejection, a classifier is allowed to output some special symbol, which indicates the rejection.
Such type of classifiers can be seen as confidence sets, which are allowed to output $\emptyset$ or $\{0, 1\}$ and are interpreted as reject.
This line of research was initiated by~\cite{Chow57,Chow70} in the context of information retrieval, where a predefined cost of rejection was considered.
An extensive statistical study of this framework was carried in~\citep{Herbei_Wegkamp06,Bartlett_Wegkamp08,Wegkamp_Yuan11}.

Instead of considering a fixed cost for rejection, which might be too restrictive, one may define two entities: probability of rejection and the probability of missclassification.
In the spirit of conformal prediction,~\cite{Lei14} aims at minimizing the probability rejection provided a fixed upper bound on the probability of missclassification.
In contrast,~\cite{Denis_Hebiri15} consider a reversed problem of minimizing the probability of missclassification given a fixed upper bound on the probability of rejection.

Once the multi-class classification is considered, there are several possible ways to extend the binary case: the confidence set approach and the rejection approach.
The reject counterpart is a more studied and known version, though it lacks statistical analysis.
To the best of our knowledge the only work which provides statistical guarantees is~\citep{Ramaswamy_Tewari_Agarwal18}.

As for the confidence set approach there are again two possibilities, similar to the binary case.
The one that is considered in this work was proposed by~\cite{Denis_Hebiri17}, where the authors analyse an ERM algorithm and derive oracle inequalities under the margin assumption~\citep{Tsybakov04}.
More specifically, they consider a convex surrogate of the error $\perf(\cdot)$ which relies on a convex real valued loss function $\phi$.
For a suitable choice of the convex function $\phi$ they show that, under Assumption~\ref{ass:continuity_cdf}, their $\beta$-Oracle satisfies
\begin{align*}
    \Gamma_{\beta}^*(\cdot) = \enscond{k \in [K]}{f^{*}_{k}(\cdot)\geq G_{f^{*}}^{-1}(\beta)}\enspace,
\end{align*}
where the function $f^{*}$ depends on $\phi$ and the value of $G_{f^{*}}^{-1}(\beta)$ is defined similarly to the present manuscript.
They propose a two step estimation procedure of the $\beta$-Oracle set.
Based on the ERM algorithm, they first estimate $f^{*}$ and in the second step, they estimate the threshold $G_{f^{*}}^{-1}(\beta)$ with an unlabeled sample.
This procedure is in the same spirit as the semi-supervised procedure~\eqref{eq:IntroSemiSupEstimator}.
Under mild assumptions, they provide an upper bound on the excess risk and obtain a rate of convergence of order $({n}/{\log n})^{-\alpha/(\alpha+s)} + {N}^{-1/2}$, with $s$ being a parameter that depends on the function $\phi$ and $\alpha$ being the margin parameter.
Note that this rate is slower than the rate obtained in the standard classification framework.

The conformal prediction theory \citep{Vovk_Gammerman_Shafer05} suggests to minimize the information level with a fixed budget on the error level.
Statistical properties of this framework were considered in the work of~\cite{Mauricio_Jing_Wasserman18}.
Their objective is formulated for some $a \in (0, 1)$ as
\begin{align*}
    \Gamma^*_a \in \argmin\enscond{\info(\Gamma)}{\Gamma \in \Upsilon \text{ \st } \perf(\Gamma) \leq a}\enspace,
\end{align*}
and such a confidence set is called a least ambiguous confidence set with bounded error rate.
The authors show that under Assumption~\ref{ass:continuity_cdf} this oracle set can be described as a thresholding of the regression function
\begin{align*}
    \Gamma^*_a(\cdot) = \enscond{k \in [K]}{p_k(\cdot) \geq t_a}\enspace,
\end{align*}
where the threshold $t_a$ is defined as
\begin{align*}
    t_a
    &=
    \sup\enscond{t \in [0, 1]}{\sum_{k = 1}^L \Prob(p_k(X) \geq t\,|\, Y = k)\Prob(Y = k) \geq 1 - a}\enspace.
\end{align*}
Notice that this framework is very similar to~\citep{Denis_Hebiri17} in the treatment of the Bayes optimal confidence set, as in both cases they are obtained via thresholding of the posterior distribution of the labels.
\cite{Mauricio_Jing_Wasserman18} also proceed in two steps as here, that is, they first estimate the posterior distribution $p_k(\cdot)$ for all $k \in [K]$ and estimate the threshold $t_a$ after.
However, they require the second \emph{labeled} dataset for the estimator of $t_a$, due to the presence of $\Prob(Y = k)$, the marginal distribution of the labels.
Besides, their theoretical analysis is carried out under a different set of assumptions on the joint distribution $\Prob$.
Apart from the standard margin assumption, they require a so-called detectability, that is, they require that the upper bound in the margin assumption is tight.
Under these assumptions they provide an upper bound on the Hamming excess risk and obtain a rate of convergence of order $\bigO((n / \log n)^{-1/2})$.

Interestingly, both approaches can be encompassed into the constrained estimation framework~\citep{Anbar77,Lepskii90,Brown_Low96}, where one would like to construct an estimator with some prescribed properties.
These properties are typically reflected by the form of the risk which in our case is the discrepancy measure, that is, the sum of error and information discrepancies.
Thus, both frameworks of~\cite{Mauricio_Jing_Wasserman18, Denis_Hebiri17} can be seen as an extension of the constrained estimation to the classification problems.
From the modeling point of view, we believe that the two frameworks can co-exist nicely and a particular choice depends on the considered application.
The major difference between the present work and those by~\cite{Denis_Hebiri17} and~\cite{Mauricio_Jing_Wasserman18} is the minimax analysis which we provide here and our treatment of semi-supervised techniques.

As already pointed out, the confidence set estimation problem is closely related to the level set estimation setup~\citep{Hartigan87, Polonik95, Tsybakov97, Rigollet_Vert09}.
This problem focuses on the estimation of a level set defined as
\begin{align*}
    \Gamma_p(\lambda) = \{x \in \mathbb{R}^{d} \; : \; p(x) \geq \lambda\},
\end{align*}
where $p$ is the density of the observations and $\lambda > 0$ is some fixed value.
Given a sample $X_1, \ldots, X_n$ distributed according the density $p$ the goal is to estimate $\Gamma_p(\lambda)$.
In~\citep{Rigollet_Vert09}, the authors study plug-in density level set estimators through the measure of symmetric differences and {\it the excess mass}.
In confidence set estimation the measure of symmetric differences is the Hamming risk whereas {\it the excess mass} is the excess risk.
They show that kernel based estimators are optimal in the minimax sense over a H\"older class of densities and under a margin type assumption~\citep{Polonik95, Tsybakov04}.
In particular, they derive fast rates of convergence, that is faster than $n^{-1/2}$,  for {\it the excess mass}.
In the level set estimation problem, the threshold $\lambda$ is chosen beforehand; whereas in our work, the threshold $G^{-1}(\beta)$ depends on the distribution of the data which makes the statistical analysis more difficult.

On the other part, the confidence set estimation problem is directly related to the standard classification settings.
This problem has been widely studied from a theoretical point of view in the binary classification framework.
\cite{Audibert_Tsybakov07} study the statistical performance of plug-in classification rules under assumptions which involve the smoothness of the regression function and the margin condition.
In particular, they derive fast rates of convergence for plug-in classifiers based on local polynomial estimators~\citep{Stone77,Tsybakov86,Audibert_Tsybakov07} and show their optimality in the minimax sense.
One of the aim of present work is to extend these results to the confidence set classification framework.

Another part of our work is to provide a comparison between supervised and semi-supervised procedures.
Semi-supervised methods are studied in several papers~\citep{Vapnik98, Rigollet07, Singh_Nowak_Zhu09, Bellec_Dalalyan_Grappin_Paris18} and references therein.
A simple intuition can be provided on whether one should or not expect a superior performance of the semi-supervised approach.
Imagine a situation when the unlabeled sample $\data_N$ is so large that one can approximate $\Prob_X$ up to any desired precision, then, if the optimal decision is independent of $\Prob_X$, the semi-supervised estimators are not to be considered superior over the supervised estimation.
This is the case in a lot of classical problems of statistics, where the inference is solely governed by the behavior of the conditional distribution $\Prob_{Y | X}$ (for instance regression or binary classification).
The situation might be different once the optimal decision relies on the marginal distribution $\Prob_X$.
In this case, as suggested by our findings, the semi-supervised approach might or not outperform the supervised one even in the context of the same problem.
Similar conclusions were stated by~\cite{Singh_Nowak_Zhu09} in the context of learning under the cluster assumption~\citep{Rigollet07}.

\subsection{Main contributions}
\label{subsec:contribution}

Bellow we summarize our contributions.
\begin{itemize}
    \item  Our results focus on the case where the regression $p$ belongs to a H\"older class and satisfy the margin condition.
    Under these assumptions, we establish lower bounds on the minimax rates, defined in Section~\ref{subsec:minimax_estimation} in the confidence set framework.
    \item As important consequences of our results, we first show that top-$\beta$ type procedures are in general inconsistent.
    Furthermore, by providing a rigorous definition of the semi-supervised and supervised estimators, we describe the situations when the semi-supervised estimation should be considered superior to its supervised counterpart.
    Interestingly, our analysis suggests that these regimes are governed by the interplay of the family of distributions and by the considered measure of performance.
    Besides, we show that in our settings supervised procedures cannot achieve fast rates, that is, its rate cannot be faster than $n^{-1/2}$.
    In contrast, some other classical settings~\citep{Audibert_Tsybakov07,Rigollet_Vert09,Herbei_Wegkamp06} allow to achieve faster rates for supervised methods.
    \item  We provide supervised and semi-supervised estimation procedures, which are optimal or optimal up to an extra logarithmic factor.
    Importantly, our results show that semi-supervised plug-in procedure based on local polynomial estimators can achieve fast rates, provided that the size of the unlabeled samples is large enough.
    %This rate is better than the rate obtained by~\citep{Denis_Hebiri17}.
    \item Finally, we perform a numerical evaluation of the proposed plug-in algorithms against the top-$\beta$ counterparts.
    This part supports our theoretical results and empirically demonstrates the reason to consider more involved procedures.
\end{itemize}
%%%%%%%%%%%%%%%%%%%%%%%%%%%%%%%%%%%%%%%%%%%%%%%%%%%%%%%%%%%%%%%%%%%%%%%%%%%%%%%
%%%%%%%%%%%%%%%%%%%%%%%%%%%%%%%%%%%%%%%%%%%%%%%%%%%%%%%%%%%%%%%%%%%%%%%%%%%%%%%
\subsection{Organization of the paper}
\label{subsec:organization_of_the_paper}
%%%%%%%%%%%%%%%%%%%%%%%%%%%%%%%%%%%%%%%%%%%%%%%%%%%%%%%%%%%%%%%%%%%%%%%%%%%%%%%
%%%%%%%%%%%%%%%%%%%%%%%%%%%%%%%%%%%%%%%%%%%%%%%%%%%%%%%%%%%%%%%%%%%%%%%%%%%%%%%

The paper is organized as follow.
In Section~\ref{sec:classConfSet}, we put some additional notation and introduce the family of distributions $\class{P}$ that we consider.
Section~\ref{sec:lower_bounds} is devoted to the lower bounds on the minimax rates and their implications.
In Section~\ref{sec:upper_bounds} we introduce the proposed algorithm, establish upper bounds for it, and evaluate its numerical performance.
We conclude this paper by Sections~\ref{sec:discussion} and~\ref{sec:conclusion} where we discuss and sum-up our results.

%%%%%%%%%%%%%%%%%%%%%%%%%%%%%%%%%%%%%%%%%%%%%%%%%%%%%%%%%%%%%%%%%%%%%%%%%%%%%%%
    %%%%%%%%%%%%%%%%%%%%%%%%%%%%%%%%%%%%%%%%%%%%%%%%%%%%%%%%%%%%%%%%%%%%%%%%%%%%%%%
%%%%%%%%%%%%%%%%%%%%%%%%%%%%%%%%%%%%%%%%%%%%%%%%%%%%%%%%%%%%%%%%%%%%%%%%%%
\section{Class of confidence sets}
\label{sec:classConfSet}
%%%%%%%%%%%%%%%%%%%%%%%%%%%%%%%%%%%%%%%%%%%%%%%%%%%%%%%%%%%%%%%%%%%%%%%%%%%%%%%
%%%%%%%%%%%%%%%%%%%%%%%%%%%%%%%%%%%%%%%%%%%%%%%%%%%%%%%%%%%%%%%%%%%%%%%%%%%%%%%
First let us introduce some generic notation that is used throughout this work.
For two numbers $a, a' \in \bbR$ we denote by $a \vee a'$ (resp. $a \wedge a'$) the maximum (resp. minimum) between $a$ and $a'$.
For a positive real number $a$ we denote by $\floor{a}$ (resp. $\ceil{a}$) the largest (resp. the smallest) non-negative integer that is less than or equal (resp. greater than or equal) to $a$.
The standard Euclidean norm of a vector $x \in \bbR^d$ is denoted by $\norm{x}$ and the standard Lebesgue measure is denoted by $\Leb(\cdot)$.
A Euclidean ball centered at $x \in \bbR^d$ of radius $r > 0$ is denoted by $\class{B}(x, r)$.
For an arbitrary Borel measure $\mu$ on $\bbR^d$ that is absolutely continuous \wrt the Lebesgue measure we denote by $\supp(\mu)$ its support, that is, the set where the Radon-Nikodym derivative of $\mu$ \wrt $\Leb$ is strictly positive.
For a vector function $p: \bbR^d \mapsto \bbR^K$ and a Borel measure $\mu$ on $\bbR^d$ we define the infinity norm of $p$ as
\begin{align*}
    \norm{p}_{\infty, \mu} \eqdef \inf\enscond{C \geq 0}{\max_{k \in [K]}\abs{p_k(x)} \leq C,\,\,\text{{\em a.e.} $x \in \bbR^d$ \wrt } \mu}\enspace.
\end{align*}
In this work $C$ or its lower-cased versions always refer to some constants which might different from line to line.
Importantly, all these constants are independent of $n, N$ but could depend on $K, d$ and other parameters which are assumed to be fixed.
Before introducing the families of distributions $\class{P}$ that are considered in this work we need the following definitions.
\begin{ass}[$\alpha$-margin assumption]
    \label{ass:margin_assumption}
    We say that the distribution $\Prob$ of the pair $(X, Y) \in \bbR^d \times [K]$ satisfies $\alpha$-margin assumption if there exists $C_1 > 0$ and $t_0 \in (0, 1)$ such that for every positive $t \leq t_0$
    \begin{align*}
        \Prob_X\pare{0 < \abs{p_k(X) - G^{-1}(\beta)} \leq t} \leq C_1t^\alpha \enspace.
    \end{align*}
\end{ass}
Let us point out an important consequence of Assumption~\ref{ass:continuity_cdf}. We have that the condition
\begin{align*}
    \quad \Prob_X\pare{\abs{p_k(X) - G^{-1}(\beta)} \leq t} \leq C_1t^\alpha\enspace,
\end{align*}
for all $t \in [0, t_0]$ is equivalent to Assumption~\ref{ass:margin_assumption}.
Indeed, since the random variables $p_k(X)$'s cannot concentrate at a constant level, in particular at $G^{-1}(\beta)$.
Moreover, again due to the continuity Assumption~\ref{ass:continuity_cdf} we have
\begin{align*}
    \lim_{t \rightarrow +0} \Prob_X\pare{\abs{p_k(X) - G^{-1}(\beta)} \leq t} = 0\enspace,
\end{align*}
thus the $\alpha$-margin Assumption~\ref{ass:margin_assumption} specifies the rate of this convergence.
Finally, the restriction of the range of $t$ to $[0, t_0]$ in $\alpha$-margin Assumption~\ref{ass:margin_assumption} does not affect its global behavior as for all $t \in [0,1]$
\begin{align*}
    \Prob_X\pare{0 < \abs{p_k(X) - G^{-1}(\beta)} \leq t} \leq c_1t^\alpha,\quad \text{with }c_1 = C_1 \vee t_0^{-\alpha}\enspace.
\end{align*}

Let $c_0$ and $r_0$ be two positive constants. We say that a Borel set $A\subset \bbR^d$ is a $(c_0, r_0)$-regular set if
\begin{align*}
    \Leb\pare{A \cap \class{B}(x,r)} \geq c_0 \Leb\pare{\class{B}(x,r)},\qquad \forall r\in(0,r_0], \ \forall x\in A\enspace.
\end{align*}
\begin{definition}[Strong density]
    We say that the probability measure $\Prob_X$ on $\bbR^d$ satisfies the $(\mu_{\min}, \mu_{\max}, c_0, r_0)$-strong density assumption if it is supported on a compact $(c_0,r_0)$-regular set $A \subset \bbR^d$ and has a density $\mu$ w.r.t. the Lebesgue measure such that $\mu(x) = 0$ for all $x \in \bbR^d \setminus A$ and
    \begin{align*}
        0<\mu_{\min} \leq \mu(x) \leq \mu_{\max}<\infty,\quad \forall x \in A\enspace.
    \end{align*}
\end{definition}
\begin{definition}[H\"{o}lder class,~\cite{Tsybakov08}]
    \label{def:holder_class}
    We say that a function $h: \bbR^d \rightarrow \bbR$ is $(\gamma, L)$-H\"{o}lder for $\gamma > 0$ and $L > 0$ if $h$ is $\floor{\gamma}$ times continuously differentiable and $\forall x, x' \in \bbR^d$ we have
    \begin{align*}
        \abs{h(x') - h_x(x')} \leq L \norm{x - x'}^\gamma\enspace,
    \end{align*}
    where $h_x(\cdot)$ is the Taylor polynomial of degree $\floor{\gamma}$ of $h(\cdot)$ at the point $x \in \bbR^d$.
    Consequently, the set of all functions from $\bbR^d$ to $\bbR$ satisfying the above conditions is called $(\gamma, L, \bbR^d)$-H\"{o}lder and is denoted by $\class{H}(\gamma, L, \bbR^d)$.
\end{definition}
\begin{definition}
    We denote by $\class{P}(L, \gamma, \alpha)$ a set of joint distributions on $\bbR^d \times [K]$ which satisfies the following conditions
    \begin{itemize}
        \item the marginal $\Prob_X$ satisfies the $(\mu_{\min}, \mu_{\max}, c_0, r_0)$-strong density,
        \item for all $k \in [K]$ the $k^{\text{th}}$ regression function $p_k(\cdot) = \Prob(Y = k | X = \cdot)$ belongs to the $(\gamma, L, \bbR^d)$-H\"{o}lder class, that is $p_k \in \class{H}(\gamma, L, \bbR^d)$ for all $k \in [K]$,
        \item for all $k \in [K]$ the regression function $p_k$ satisfy the $(C_1, \alpha, \beta)$-Margin assumption,
        \item for all $k \in [K]$, the cumulative distribution function $F_{p_k}$ of $p_k(X)$ is continuous.
    \end{itemize}
\end{definition}
The family of distributions $\class{P}(L, \gamma, \alpha)$ is similar to the one considered in~\citep{Audibert_Tsybakov07} in the context of binary classification.
The only major difference is the continuity Assumption~\ref{ass:continuity_cdf}, which does not allow to re-use in a straightforward way their construction for lower bounds.

%%%%%%%%%%%%%%%%%%%%%%%%%%%%%%%%%%%%%%%%%%%%%%%%%%%%%%%%%%%%%%%%%%%%%%%%%%%%%%%
    %%%%%%%%%%%%%%%%%%%%%%%%%%%%%%%%%%%%%%%%%%%%%%%%%%%%%%%%%%%%%%%%%%%%%%%%%%%%%%%
%%%%%%%%%%%%%%%%%%%%%%%%%%%%%%%%%%%%%%%%%%%%%%%%%%%%%%%%%%%%%%%%%%%%%%%%%%%%%%%
\section{Lower bounds}
\label{sec:lower_bounds}
%%%%%%%%%%%%%%%%%%%%%%%%%%%%%%%%%%%%%%%%%%%%%%%%%%%%%%%%%%%%%%%%%%%%%%%%%%%%%%%
%%%%%%%%%%%%%%%%%%%%%%%%%%%%%%%%%%%%%%%%%%%%%%%%%%%%%%%%%%%%%%%%%%%%%%%%%%%%%%%
The main results in the present work are the lower bounds we provide in this section.
In particular, we establish in Section~\ref{subsec:naive_approach} the inconsistency of top-$\beta$ procedures (see Eq.~\eqref{eq:IntroTopBetaEstimator} for a definition of the method).
Therefore more elaborate methods are required in this framework.
As pointed out in the introduction, we distinguish two types of estimators: \emph{supervised} and \emph{semi-supervised} for which we provide lower bounds in Section~\ref{subsec:supervised_estimators}.
The obtained rates highlight the benefit of the semi-supervised approach in the context of the confidence set classification.

Before considering the lower bounds, let us first display connection between the different minimax rates. Such links are used in the proofs of the lower bounds.

\begin{prop}
    \label{prop:distanceComparison}
    Let $\Gamma$ be a measurable function from $\bbR^d$ to $2^{[K]}$, $\beta \in [k]$ and assume that Assumption~\ref{ass:continuity_cdf} is fulfilled, then
    \begin{align*}
        \perf(\Gamma) - \perf(\Gamma^*_\beta)
        &=
        \risk_\beta(\Gamma) - \risk_\beta(\Gamma^*_\beta) + G^{-1}(\beta)\left(\beta - \info(\Gamma)\right)\enspace,\\
        \risk_\beta(\Gamma) - \risk_\beta(\Gamma^*_\beta)
        &=
        \sum_{k =1}^K \Exp_{\Prob_X}\left[|p_k(X)-G^{-1}(\beta)|\1_{\{k \in \Gamma(X) \triangle \Gamma^{*}_{\beta}(X)\}}\right]\enspace.
    \end{align*}
    Furthermore, if additionally Assumption~\ref{ass:margin_assumption} is satisfied with $\alpha > 0$, then there exist $C > 0$ which depends only on $K, \alpha, C_1$ such that for any pair of confidence set classifiers $\Gamma, \Gamma^{'}$ it holds that
    \begin{equation}
        \label{eq:H_D_comparison_lower_bound}
        \Exp_{\Prob_X}\abs{\Gamma(X) \triangle \Gamma^{'}(X)} \leq C \left(\risk_\beta(\Gamma) - \risk_\beta(\Gamma')\right)^{\alpha/(\alpha+1)}\enspace.
    \end{equation}
\end{prop}
\begin{prop}
    \label{prop:order_risks}
    For any $K \geq 2$, $\beta \in [K]$ and $n, N \in \bbN$ the following relation between minimax rates holds:
    \begin{align*}
        \excess^{\Ham}_{n, N}(\classGamma;\class{P}) \geq \excess^{\Disc}_{n, N}(\classGamma;\class{P}) \geq \excess^{\Exc}_{n, N}(\classGamma;\class{P})\enspace.
    \end{align*}
\end{prop}
Proposition~\ref{prop:distanceComparison}, and in particular Eq.~\eqref{eq:H_D_comparison_lower_bound} gives an easy way to establish a lower bound on $\excess^{\Exc}_{n, N}(\classGamma;\class{P})$ via a lower bound on the Hamming distance $\excess^{\Ham}_{n, N}(\classGamma;\class{P})$.
However, this approach does not allow to get $(N + n)^{-1/2}$ (resp. $n^{-1/2}$) part of the rate in the lower bound of $\excess^{\Exc}_{n, N}(\hat\Upsilon_{\SSE}, \class{P})$ (resp. $\excess^{\Exc}_{n, N}(\hat\Upsilon_{\SE}, \class{P})$).
Besides, Proposition~\ref{prop:order_risks} allows to prove a lower bound on the discrepancy $\excess^{\Disc}_{n, N}(\classGamma;\class{P})$ with the correct rate via the lower bound on the excess risk $\excess^{\Exc}_{n, N}(\classGamma;\class{P})$.

%%%%%%%%%%%%%%%%%%%%%%%%%%%%%%%%%%%%%%%%%%%%%%%%%%%%%%%%%%%%%%%%%%%%%%%%%%%%%%%
%%%%%%%%%%%%%%%%%%%%%%%%%%%%%%%%%%%%%%%%%%%%%%%%%%%%%%%%%%%%%%%%%%%%%%%%%%%%%%%
\subsection{Inconsistency of the top-\texorpdfstring{$\beta$}{Lg} procedure}
\label{subsec:naive_approach}
%%%%%%%%%%%%%%%%%%%%%%%%%%%%%%%%%%%%%%%%%%%%%%%%%%%%%%%%%%%%%%%%%%%%%%%%%%%%%%%
%%%%%%%%%%%%%%%%%%%%%%%%%%%%%%%%%%%%%%%%%%%%%%%%%%%%%%%%%%%%%%%%%%%%%%%%%%%%%%%
Before stating our results on the supervised and the semi-supervised estimators, we discuss another interesting class of confidence sets, which might be a natural choice at the first sight.
We consider estimators which consists of $\beta$ classes at every point $x \in \bbR^d$ since such estimators naturally satisfy $\info{(\hGamma)} = \beta$.
Let us denote by $\hat\Upsilon_\beta$ the set of all estimators $\hGamma$ such that $|{{\hGamma(x)}}| = \beta$ for all $x \in \bbR^d$, that is,
\begin{align*}
    \hat\Upsilon_\beta = \enscond{\hGamma \in \hat\Upsilon}{\lvert{\hGamma(x)}\rvert = \beta, \,\,\text{{\em a.e.} } x \in \bbR^d \text{ \wrt } \Leb}\enspace.
\end{align*}
Despite an obvious restriction on the cardinal of the confidence sets, the family of estimators $\hat\Upsilon_\beta$ is rather broad.
Indeed, every procedure which estimates the regression functions $p_k(\cdot)$'s and includes the top $\beta$ scores as the output are included in $ \hat\Upsilon_\beta$.
The nature of the estimator can also be different, that is, the estimates could be based on the ERM, non-parametric or parametric approaches.
Clearly, the family $\hat\Upsilon_\beta$ is neither included in $\hat\Upsilon_{\SE}$ nor in $\hat\Upsilon_{\SSE}$ and has a non-trivial intersection with both.
The next result states that there is no uniformly consistent estimator $\hGamma \in \hat\Upsilon_\beta$ over the family of distributions $\class{P}(L, \gamma, \alpha)$.
\begin{prop}
    \label{prop:no_consistency_top_beta}
    Assume that $K \geq 4$, $\beta \in [\floor{K / 2} - 1]$ and $\beta \geq 2$, then for all $n, N \in \bbN$ we have
    \begin{align*}
        \excess^{\Exc}_{n, N}\pare{\hat\Upsilon_\beta; \class{P}(L, \gamma, \alpha)} \geq \frac{\beta - 1}{4 K} \enspace.
    \end{align*}
\end{prop}
The proof builds an explicit construction of a distribution $\Prob$ whose $\beta$-Oracle satisfies $\absin{\Gamma^*_\beta(x)} > \beta$ for all $x$ in some $A \subset \bbR^d$ with $\Prob_X(A) > 0$.
Indeed, if such a distribution exists then there is no estimator in $\hat\Upsilon_\beta$ that would consistently estimate this $\beta$-Oracle.
The negative result established in Proposition~\ref{prop:no_consistency_top_beta} is rather instructive by itself as it advocates that a more involved estimation procedure ought to be constructed.
%%%%%%%%%%%%%%%%%%%%%%%%%%%%%%%%%%%%%%%%%%%%%%%%%%%%%%%%%%%%%%%%%%%%%%%%%%%%%%%
%%%%%%%%%%%%%%%%%%%%%%%%%%%%%%%%%%%%%%%%%%%%%%%%%%%%%%%%%%%%%%%%%%%%%%%%%%%%%%%
\subsection{Supervised vs semi-supervised estimation}
\label{subsec:supervised_estimators}
%%%%%%%%%%%%%%%%%%%%%%%%%%%%%%%%%%%%%%%%%%%%%%%%%%%%%%%%%%%%%%%%%%%%%%%%%%%%%%%
%%%%%%%%%%%%%%%%%%%%%%%%%%%%%%%%%%%%%%%%%%%%%%%%%%%%%%%%%%%%%%%%%%%%%%%%%%%%%%%
Clearly, estimators which achieve the infimum in the minimax rates are either supervised or semi-supervised, thus a lower bound on $\excess^{\square}_{n, N}(\hat\Upsilon_{\SE};\class{P})$ together with a lower bound on $\excess^{\square}_{n, N}(\hat\Upsilon_{\SSE};\class{P})$ yield a lower bound on $\excess^{\square}_{n, N}(\hat\Upsilon;\class{P})$.
However, a lower bound on $\excess^{\square}_{n, N}(\hat\Upsilon;\class{P})$ does not discriminate between the supervised and the semi-supervised estimators.
\begin{theo}[Supervised estimation]
\label{thm:lower_bound_supervised}
    Let $K \leq 3$, $\beta \in [\floor{K / 2} - 1]$. If $2 \alpha\ceil{\frac{\gamma}{2}}  \leq d$, then there exist constants $c, c', c'' > 0$ such that for all $n, N \in \bbN$
    \begin{align*}
        \excess^{\Ham}_{n, N}(\hat\Upsilon_{\SE};\class{P}(L, \gamma, \alpha))
        &\geq
        c\pare{{n}^{-\frac{\alpha \gamma}{2\gamma + d}}\bigvee n^{-1/2}}\enspace,\\
        \excess^{\Exc}_{n, N}(\hat\Upsilon_{\SE};\class{P}(L, \gamma, \alpha))
        &\geq
        c'\pare{n^{-\frac{(1 + \alpha) \gamma}{2\gamma + d}} \bigvee n^{-1/2}}\enspace,\\
        \excess^{\Disc}_{n, N}(\hat\Upsilon_{\SE};\class{P}(L, \gamma, \alpha))
        &\geq
        c''\pare{n^{-\frac{(1 + \alpha) \gamma}{2\gamma + d}} \bigvee n^{-1/2}}\enspace.
    \end{align*}
\end{theo}
Based on this results we observe that the lower bound for the Hamming risk $\excess^{\Ham}_{n, N}$ is slower than those for the other risks.
It is even more significant that the best rate that a supervised estimator can achieve for all of the risks is $n^{-1/2}$ even if the margin assumption holds.
This is the major difference with the classical settings where the value of threshold is known (such as classification and level set estimation).
Indeed, under the same assumptions on the family of distributions, besides the continuity Assumption~\ref{ass:continuity_cdf}, the minimax rate in those frameworks is $n^{-{(1 + \alpha) \gamma}/{(2\gamma + d)}}$ as proved for instance in~\citep{Audibert_Tsybakov07,Rigollet_Vert09}.
Next theorem deals with semi-supervised procedures and displays another behavior.
%%%%%%%%%%%%%%%%%%%%%%%%%%%%%%%%%%%%%%%%%%%%%%%%%%%%%%%%%%%%%%%%%%%%%%%%%%%%%%%
%%%%%%%%%%%%%%%%%%%%%%%%%%%%%%%%%%%%%%%%%%%%%%%%%%%%%%%%%%%%%%%%%%%%%%%%%%%%%%%
% \subsection{Semi-supervised estimation}
% \label{subsec:semi_supervised_estimators}
%%%%%%%%%%%%%%%%%%%%%%%%%%%%%%%%%%%%%%%%%%%%%%%%%%%%%%%%%%%%%%%%%%%%%%%%%%%%%%%
%%%%%%%%%%%%%%%%%%%%%%%%%%%%%%%%%%%%%%%%%%%%%%%%%%%%%%%%%%%%%%%%%%%%%%%%%%%%%%%
\begin{theo}[Semi-supervised estimation]
        \label{thm:lower_bound_semi_supervised}
    Let $K \geq 3$, $\beta \in [\floor{K / 2} - 1]$.  If $2 \alpha\ceil{\frac{\gamma}{2}}  \leq d$, then there exist constants $c, c', c'' > 0$ such that for all $n, N \in \bbN$
    \begin{align*}
        \excess^{\Ham}_{n, N}(\hat\Upsilon_{\SSE};\class{P}(L, \gamma, \alpha))
        &\geq
        c\pare{n^{-\frac{\alpha \gamma}{2\gamma + d}} \bigvee (n + N)^{-1/2}}\enspace,\\
        \excess^{\Exc}_{n, N}(\hat\Upsilon_{\SSE};\class{P}(L, \gamma, \alpha))
        &\geq
        c'\pare{n^{-\frac{(1 + \alpha) \gamma}{2\gamma + d}} \bigvee (n + N)^{-1/2}}\enspace,\\
        \excess^{\Disc}_{n, N}(\hat\Upsilon_{\SSE};\class{P}(L, \gamma, \alpha))
        &\geq
        c''\pare{n^{-\frac{(1 + \alpha) \gamma}{2\gamma + d}} \bigvee (n + N)^{-1/2}}\enspace.
    \end{align*}
\end{theo}
First, observe that the lower bound for the Hamming distance is, as in the supervised setting, worse than for the other measures of performance.
However there is a major difference with the supervised case:
as compared to Theorem~\ref{thm:lower_bound_supervised}, it is possible for a semi-supervised estimator to achieve rates that are faster than $n ^{-1/2}$ if the size of the unlabeled dataset $N \in \bbN$ is large enough.
In particular, when we consider $\excess^{\Exc}_{n, N}$ or $\excess^{\Disc}_{n, N}$ the following relations are necessary to get fast rates
\begin{align*}
    (n + N)^{-1/2}
    =
    \littleo\pare{n^{-{(1 + \alpha) \gamma}/{(2\gamma + d)}}},\quad
    n^{-{(1 + \alpha) \gamma}/{(2\gamma + d)}}
    =
    \littleo(n^{-1/2})\enspace.
\end{align*}
In this case, we recover the same fast rates as in the classical settings of classification and level set estimation.
It suggests that the lack of knowledge of the threshold $G^{-1}(\beta)$ does not alter the quality of estimation for the semi-supervised procedure, provided that $N$ is sufficiently large.
% In contrast, no supervised estimator can achieve rates that are faster than $n^{-1/2}$ even if the underlying margin assumption and the smoothness are extremely favorable.
Next corollary makes these observations clearer.
\begin{corollary}
    \label{cor:lower_bound}
    Assume that the rates in Theorem~\ref{thm:lower_bound_semi_supervised} (resp. Theorem~\ref{thm:lower_bound_supervised}) are minimax, that is, there exist a confidence set $\hGamma_{\SSE}$ (resp. $\hGamma_{\SE}$) that achieves these rates.
    Regarding $\excess^{\Exc}_{n, N}$ and $\excess^{\Disc}_{n, N}$ the following conclusions hold
    \begin{itemize}
        \item There is no semi-supervised
              estimator that achieves faster rate than $\hGamma_{\SE}$ if:
              \begin{align*}
                \begin{cases}
                    \frac{(1 + \alpha) \gamma}{2\gamma + d} \leq 1/2\\
                    N \in \bbN
                \end{cases}\qquad
                   \text{ or }\qquad
                \begin{cases}
                    \frac{(1 + \alpha) \gamma}{2\gamma + d} > 1/2\\
                    N = \bigO(n)
                \end{cases}\enspace.
              \end{align*}
        \item The rate of $\hGamma_{\SSE}$ is faster than the rate of any supervised estimator if:
        \begin{align*}
            \frac{(1 + \alpha) \gamma}{2\gamma + d} > 1/2\qquad \text{ and }\qquad n = \littleo(N)\enspace.
        \end{align*}
        Moreover, if there exists $\rho > 0$ such that $n^{1 + \rho} = \littleo(N)$, then the rate of $\hGamma_{\SSE}$ is polynomially faster than $n^{-1/2}$.
        \item The rate of $\hGamma_{\SSE}$ is fast similarly to the classical frameworks if
        \begin{align*}
            \frac{(1 + \alpha) \gamma}{2\gamma + d} > 1/2\qquad \text{ and }\qquad N = \Omega\pare{n^{\frac{2(1 + \alpha) \gamma}{2\gamma + d}}}\enspace.
        \end{align*}
    \end{itemize}
\end{corollary}
Clearly, similar observation is true for the Hamming risk $\excess^{\Ham}_{n, N}$; however the regime when improvement is possible thanks to semi-supervised approaches is narrowed as $n^{-{(1 + \alpha) \gamma}/{(2\gamma + d)}} = \littleo\pare{n^{-{\alpha \gamma}/{(2\gamma + d)}}}$.
We summarize Corollary~\ref{cor:lower_bound} in Table~\ref{tab:summary_rates}.
\begin{center}
    \label{tab:summary_rates}
    \begin{table}[t!]
        \begin{tabular}{|p{1cm}||c||ccc|}
            \hline
            $\hfil\frac{(1 + \alpha) \gamma}{2\gamma + d}$ & $N, n$ & $\SE$ rate & $\SSE$  rate & $\SSE > \SE$ \\
            % \noalign{\smallskip}
            \hline
            % \hline
            % \noalign{\smallskip}
            \hline
            $\hfil\leq \frac{1}{2}$ & $ N \in \bbN$, $n \in \bbN$ & $\hfil n^{-\frac{(1 + \alpha) \gamma}{2\gamma + d}}$ & $ n^{-\frac{(1 + \alpha) \gamma}{2\gamma + d}}$ & NO \\
            % \hline
            $\hfil> \frac{1}{2}$ & $N = \bigO(n)$ & $\hfil n^{-\frac{1}{2}}$ & $ n^{-\frac{1}{2}}$ & NO \\
            % \hline
            $\hfil> \frac{1}{2}$ & $n = \littleo(N)$ & $\hfil n^{-\frac{1}{2}}$ & $\hfil N^{-\frac{1}{2}} \bigvee n^{-\frac{(1 + \alpha) \gamma}{2\gamma + d}} $ & YES \\
            % \hline
            $\hfil> \frac{1}{2}$ & $\hfil N = \Omega\pare{n^{\frac{2(1 + \alpha) \gamma}{2\gamma + d}}}$ & $\hfil n^{-\frac{1}{2}}$ &$n^{-\frac{(1 + \alpha) \gamma}{2\gamma + d}}$ & YES \\ \hline
        \end{tabular}%
        \caption{This table summarizes observations of Corollary~\ref{cor:lower_bound} for $\excess^{\Exc}_{n, N}$ and $\excess^{\Disc}_{n, N}$. Depending on the relations between $\alpha, \gamma, d$ and $N, n$ the semi-supervised approach can significantly improve the rates of convergence.}
    \end{table}
\end{center}
Essentially, the above results suggest that the advantage of the semi-supervised approaches over the supervised ones depends not only on the underlying family of distributions $\class{P}$ but also on the metric that is considered.
Yet, necessary and sufficient conditions that must be imposed in general on the problem and the metric so that the semi-supervised estimation provably improve upon the supervised one remain an open problem.

A final remark we could make before going further concerns the assumption on the parameters $\alpha$ and $\gamma$.
The condition $2 \alpha\ceil{\frac{\gamma}{2}} \leq d$ in the lower bounds is slightly more restrictive than the conditions given in~\citep{Audibert_Tsybakov07} (they have $\alpha \gamma \leq d$).
We believe that this is an artifact of our proof and could be avoided with a finer choice of hypotheses.
Simple modifications of the lower bound of~\cite{Audibert_Tsybakov07} do not work in our settings because their hypotheses are not satisfying Assumption~\ref{ass:continuity_cdf}.
In contrast, the construction of~\cite{Rigollet_Vert09} satisfies\footnote{Modified properly to fit the classification framework.} Assumption~\ref{ass:continuity_cdf} but their lower bound is limited by the condition $\alpha \gamma \leq 1$, that is, it does not cover the fast rates as long as the dimension $d > 2$.

\subsection{Sketch of the proof}
\label{subsec:sketch_of_the_proof}

In order to prove the lower bounds of Theorems~\ref{thm:lower_bound_supervised},~\ref{thm:lower_bound_semi_supervised} we actually prove two separate lower bounds on the minimax rates.
The two lower bounds that we prove are naturally connected with the proposed two-step estimator in Eq.~\eqref{eq:IntroSemiSupEstimator}, that is, the first lower bound is connected with the problem of non-parametric estimation of $p_k$ for all $k \in [K]$ and the second describes the estimation of the unknown threshold $G^{-1}(\beta)$.

In particular, the first lower bound is closely related to the one provided in~\citep{Audibert_Tsybakov07,Rigollet_Vert09}, however, the continuity Assumption~\ref{ass:continuity_cdf} makes the proof more involved and results in a final construction of hypotheses that differs significantly.
This part of our lower bound relies on Fano's inequality in the form of~\cite{Birge05}.
The second lower bound is based on two hypotheses testing and is derived by constructing two different marginal distributions of $X \in \bbR^d$ which are sufficiently close and a fixed regression function $p(\cdot)$.
Crucially, these marginal distributions admit two different values of threshold $G^{-1}(\beta)$ and thus two different $\beta$-Oracle.
In this part we make use of Pinsker's inequality, see for instance~\citep{Tsybakov08}.

In order to discriminate the supervised and the semi-supervised procedures we make use of Definition~\ref{def:supervised_semi_supervised}.
Notice that every supervised procedure thanks to Definition~\ref{def:supervised_semi_supervised} is not ``sensitive'' to the expectation taken \wrt the unlabeled dataset $\data_N$, that is, randomness is only induced by the labeled dataset $\data_n$.
This strategy allows to eliminate the dependence of the lower bound on the size of the unlabeled dataset $\data_N$ for supervised procedures.
Informally, the lower bound on $\excess^{\square}_{n, N}(\hat\Upsilon_{\SE};\class{P})$ is obtained from the lower bound on $\excess^{\square}_{n, N}(\hat\Upsilon_{\SSE};\class{P})$ by setting $N = 0$.
%%%%%%%%%%%%%%%%%%%%%%%%%%%%%%%%%%%%%%%%%%%%%%%%%%%%%%%%%%%%%%%%%%%%%%%%%%%%%%%
    %%%%%%%%%%%%%%%%%%%%%%%%%%%%%%%%%%%%%%%%%%%%%%%%%%%%%%%%%%%%%%%%%%%%%%%%%%%%%%%
%%%%%%%%%%%%%%%%%%%%%%%%%%%%%%%%%%%%%%%%%%%%%%%%%%%%%%%%%%%%%%%%%%%%%%%%%%%%%%%
\section{Upper bounds}
\label{sec:upper_bounds}
%%%%%%%%%%%%%%%%%%%%%%%%%%%%%%%%%%%%%%%%%%%%%%%%%%%%%%%%%%%%%%%%%%%%%%%%%%%%%%%
%%%%%%%%%%%%%%%%%%%%%%%%%%%%%%%%%%%%%%%%%%%%%%%%%%%%%%%%%%%%%%%%%%%%%%%%%%%%%%%
In this section, we show that we can build confidence set estimators that achieve, up to a logarithmic factor, the lower bounds stated in Theorems~\ref{thm:lower_bound_supervised}-\ref{thm:lower_bound_semi_supervised}.
In other words, those estimators are {\it nearly} optimal in the minimax sense.
To come straight to the point, we delay the construction of the estimators to Section~\ref{subsec:construction_of_the_estimator} and their properties to Section~\ref{subsec:properties_estimators}, and focus right now on their upper bounds.
\begin{theo}[Supervised estimation]
    \label{thm:upper_bound_supervised}
    Let $K \in \bbN$, $\beta \in [K - 1]$, then there exists a supervised estimator $\hGamma_{\SE} \in \hat\Upsilon_{\SE}$ and a constant $C > 0$ such that for all $n, N \in \bbN$ we have
    \begin{align*}
        \excess^{\Ham}_{n, N}(\hGamma_{\SE};\class{P}(L, \gamma, \alpha))
        &\leq
        C\pare{n^{-\frac{\alpha \gamma}{2\gamma + d}} \bigvee n^{-1/2}}\enspace,\\
        \excess^{\Exc}_{n, N}(\hGamma_{\SE};\class{P}(L, \gamma, \alpha))
        &\leq
        C'\pare{\pare{\frac{n}{\log n}}^{-\frac{(1 + \alpha) \gamma}{2\gamma + d}} \bigvee n^{-1/2}}\enspace,\\
        \excess^{\Disc}_{n, N}(\hGamma_{\SE};\class{P}(L, \gamma, \alpha))
        &\leq
        C''\pare{\pare{\frac{n}{\log n}}^{-\frac{(1 + \alpha) \gamma}{2\gamma + d}} \bigvee n^{-1/2}}\enspace.
    \end{align*}
\end{theo}

\begin{theo}[Semi-supervised estimation]
    \label{thm:upper_bound_semi_supervised}
    Let $K \in \bbN$, $\beta \in [K - 1]$, then there exists a semi-supervised estimator $\hGamma_{\SSE} \in \hat\Upsilon_{\SSE}$ and constants $C, C', C'' > 0$ such that for all $n, N \in \bbN$ we have
    \begin{align*}
        \excess^{\Ham}_{n, N}(\hGamma_{\SSE};\class{P}(L, \gamma, \alpha))
        &\leq
        C\pare{n^{-\frac{\alpha \gamma}{2\gamma + d}} \bigvee (n + N)^{-1/2}}\enspace,\\
        \excess^{\Exc}_{n, N}(\hGamma_{\SSE};\class{P}(L, \gamma, \alpha))
        &\leq
        C'\pare{\pare{\frac{n}{\log n}}^{-\frac{(1 + \alpha) \gamma}{2\gamma + d}} \bigvee \pare{n + N}^{-1/2}}\enspace,\\
        \excess^{\Disc}_{n, N}(\hGamma_{\SSE};\class{P}(L, \gamma, \alpha))
        &\leq
        C''\pare{\pare{\frac{n}{\log n}}^{-\frac{(1 + \alpha) \gamma}{2\gamma + d}} \bigvee \pare{n + N}^{-1/2}}\enspace.
    \end{align*}
\end{theo}

We show here that the lower bounds of Theorems~\ref{thm:lower_bound_supervised}-\ref{thm:lower_bound_semi_supervised} are achievable.
In particular, in the case of Hamming risk, the upper bounds are optimal; whereas for the Excess risk and the Discrepancy, the upper bounds fit the lower bounds up to a logarithmic factor.
Thus, the comments we made in Corollary~\ref{cor:lower_bound} are correct.
Let us mention that the presence of the logarithmic factor in these upper bounds is due to $\ell_\infty$-norm estimation (see Lemma~\ref{lem:wasserstein_infty_bound_main}).

Hamming risk as a measure of performance was considered in the settings of~\cite{Mauricio_Jing_Wasserman18}.
They also establish upper bounds for this measure, though do not assess their optimality.
Besides, as we already mentioned,~\cite{Denis_Hebiri17} provide an upper bound on the excess risk in the context of ERM.
Let us point out, that the comparison with these two works is not fair as the assumptions and even frameworks under which we and they formulate results are different.

%%%%%%%%%%%%%%%%%%%%%%%%%%%%%%%%%%%%%%%%%%%%%%%%%%%%%%%%%%%%%%%%%%%%%%%%%%%%%%%
%%%%%%%%%%%%%%%%%%%%%%%%%%%%%%%%%%%%%%%%%%%%%%%%%%%%%%%%%%%%%%%%%%%%%%%%%%%%%%%
\subsection{Construction of the estimators}
\label{subsec:construction_of_the_estimator}
%%%%%%%%%%%%%%%%%%%%%%%%%%%%%%%%%%%%%%%%%%%%%%%%%%%%%%%%%%%%%%%%%%%%%%%%%%%%%%%
%%%%%%%%%%%%%%%%%%%%%%%%%%%%%%%%%%%%%%%%%%%%%%%%%%%%%%%%%%%%%%%%%%%%%%%%%%%%%%%

Building estimators $\hGamma_{\SE}$ and $\hGamma_{\SSE}$ that reach the rates in the former upper bounds involves a preliminary estimators $\hat p_k$ of the regression functions $p_k$, $k \in [K]$.
These estimators $\hat p_k$ are constructed using an arbitrary half $\data_{\floor{n/2}}$ of the labeled dataset $\data_n$ and they satisfy the following assumptions.
\begin{ass}[Exponential concentration]
    \label{ass:exponential_concentration_estimator}
    There exist estimators $\hat p_k$ for all $k \in [K]$ based on $\data_{\floor{n/2}}$ and positive constants $C_1, C_2$ such that for all $k \in [K]$ and all $n \geq 2$ we have for all $\delta > 0$
    \begin{align*}
        \sup_{\Prob \in \class{P}(L, \gamma, \alpha)}\Prob^{\otimes \floor{n/2}}\pare{\abs{\hat{p}_k(x) - p_k(x)} \geq \delta} \leq C_1 \exp\pare{-C_2n^{\frac{2\gamma}{2\gamma + d}}\delta^2}\enspace,
    \end{align*}
    for almost all $x \in \bbR^d$ w.r.t. $\Prob_X$.
    % a constant $C_2 > 0$ such that for all $n \geq 2$
    % \begin{equation*}
    %     \sup_{\Prob \in \class{P}(L, \gamma, \alpha)}\Exp_{\data_{\floor{n/2}}}\norm{p - \hat p}_{\infty, \Prob_X}^{1 + \alpha} \leq C_2\pare{\frac{n}{\log n}}^{-\frac{\gamma}{2\gamma + d}}\enspace.
    % \end{equation*}
    % Moreover, there exist constants $C'_1, C'_2$ such that for all $k \in [K]$ and all $n \geq 2$ we have for all $\delta > 0$
    % \begin{align*}
    %     \sup_{\Prob \in \class{P}(L, \gamma, \alpha)}\Prob^{\otimes \floor{n/2}}\pare{\abs{\hat{p}_k(x) - p_k(x)} \geq \delta} \leq C'_1 \exp\pare{-C'_2n^{\frac{2\gamma}{2\gamma + d}}\delta^2}\enspace,
    % \end{align*}
    %for almost all $x \in \bbR^d$ \wrt $\Prob_X$.
\end{ass}
\begin{ass}[Continuity of CDF]
    \label{ass:continuity_cdf_estimator}
    For all $k \in [K]$ the cumulative distribution function $F_{\hat{p}_k}(t) \eqdef \Prob_X(\hat{p}_k(X) \leq t)$ of $\hat{p}_k(X)$ is almost surely $\Prob^{\otimes \floor{n / 2}}$ continuous on $(0, 1)$.
\end{ass}
First let us point out that Assumption~\ref{ass:exponential_concentration_estimator} induces that there exists a constant  $C > 0$ such that for all $n \geq 2$ and all $\alpha > 0$
    \begin{equation*}
        \sup_{\Prob \in \class{P}(L, \gamma, \alpha)}\Exp_{\data_{\floor{n/2}}}\norm{p - \hat p}_{\infty, \Prob_X}^{1 + \alpha} \leq C\pare{\frac{n}{\log n}}^{-\frac{(1 + \alpha)\gamma}{2\gamma + d}}\enspace.
    \end{equation*}
Assumption~\ref{ass:exponential_concentration_estimator} is commonly used in the statistical community when we deal with rates of convergence in the classification settings~\citep{Audibert_Tsybakov07,Lei14,Mauricio_Jing_Wasserman18}.
It is for instance satisfied by the locally polynomial estimator~\citep{Stone77,Tsybakov86,Audibert_Tsybakov07}.
Assumption~\ref{ass:continuity_cdf_estimator} can always be satisfied by slightly processing any estimator $\hat p$.
Indeed, assume Assumption~\ref{ass:continuity_cdf_estimator} fails to be satisfied by some estimator $\hat p$.
It means that there exists a subset of $\bbR^d$ of non-zero measure such that at least one $\hat{p}_k$, with $k\in[K]$, is constant on this set.
Then, if we add a \emph{deterministic} continuous function of a sufficiently bounded variation\footnote{It is sufficient to make sure that adding the function preserves its statistical properties, that is, Assumption~\ref{ass:exponential_concentration_estimator}} to $\hat p$ such regions can no longer exist.

Since, the threshold level $G^{-1}(\beta)$ is not known beforehand, it ought to be estimated using data.
A straightforward estimator of this threshold can be constructed using the unlabeled dataset $\data_N$.
To make our presentation mathematically correct we introduce the following notation $\data_n = \data_{\floor{n / 2}} \bigcup \data_{\ceil{n / 2}}$, where $\data_{\floor{n / 2}}$ is the dataset used to build the estimators $\hat p_k$ for $k \in [K]$.
Now, all the labels are removed from $\data_{\ceil{n / 2}}$, that is it consists of $\ceil{n / 2}$ \iid samples from $\Prob_X$.
The supervised and semi-supervised estimators of $G(\cdot)$ are defined as
\begin{align*}
    \hat G_{\SE}(\cdot)
    &=
    \frac{1}{\ceil{n / 2}}\sum_{X \in \data_{\ceil{n / 2}}}\sum_{k = 1}^K\ind{\hat p_k(X) > \cdot}\enspace,\\
    \hat G_{\SSE}(\cdot)
    &=
    \frac{1}{\ceil{n / 2} + N}\sum_{X \in \data_N \bigcup \data_{\ceil{n / 2}}}\sum_{k = 1}^K\ind{\hat p_k(X) > \cdot}\enspace,
\end{align*}
respectively.
Finally, we are in position to define $\hGamma_{\SE}$ and $\hGamma_{\SSE}$ as
\begin{align*}
    \hGamma_{\SE}(x) &= \enscond{k \in [K]}{\hat p_k(x) \geq \hat G_{\SE}^{-1}(\beta)}\enspace,\\
    \hGamma_{\SSE}(x) &= \enscond{k \in [K]}{\hat p_k(x) \geq \hat G_{\SSE}^{-1}(\beta)}\enspace,
\end{align*}
for all $x\in\bbR^d$. Note that $\hGamma_{\SE}$ is clearly supervised in the sense of Definition~\ref{def:supervised_semi_supervised}, as it is independent of the unlabeled sample $\data_N$.
In contrast, $\hGamma_{\SSE}$ is semi-supervised, since we can find two samples $\data_N$ and $\data'_N$ which induce different confidence sets.
To show that the estimators introduced in this section satisfy the statements of Theorems~\ref{thm:upper_bound_supervised}-\ref{thm:upper_bound_semi_supervised} we refine the proof technique used in~\citep{Denis_Hebiri17}.
That is, we introduce an intermediate quantity
\begin{align*}
    \tilde{G}(\cdot) \eqdef \sum_{k = 1}^K\Prob_X\pare{\hat p_k(X) > \cdot}\enspace,
\end{align*}
and the associated confidence set, which we refer to as the pseudo Oracle confidence set given for all $x\in\bbR^d$ by
\begin{align*}
    \tilde{\Gamma}(x) \eqdef \enscond{k \in [K]}{\hat p_k(x) \geq \tilde{G}^{-1}(\beta)}\enspace.
\end{align*}
The confidence set $\tilde{\Gamma}$ assumes knowledge of the marginal distribution $\Prob_X$ and is seen as an idealized version of both $\hGamma_{\SE}$ and $\hGamma_{\SSE}$, note however, that the pseudo Oracle $\tilde{\Gamma}$ is not an estimator.

%%%%%%%%%%%%%%%%%%%%%%%%%%%%%%%%%%%%%%%%%%%%%%%%%%%%%%%%%%%%%%%%%%%%%%%%%%%%%%%
%%%%%%%%%%%%%%%%%%%%%%%%%%%%%%%%%%%%%%%%%%%%%%%%%%%%%%%%%%%%%%%%%%%%%%%%%%%%%%%
\subsection{Properties of the plug-in confidence sets}
\label{subsec:properties_estimators}
%%%%%%%%%%%%%%%%%%%%%%%%%%%%%%%%%%%%%%%%%%%%%%%%%%%%%%%%%%%%%%%%%%%%%%%%%%%%%%%
%%%%%%%%%%%%%%%%%%%%%%%%%%%%%%%%%%%%%%%%%%%%%%%%%%%%%%%%%%%%%%%%%%%%%%%%%%%%%%%

An important step of our analysis is the following lemma, that bounds the difference between $\tilde{G}^{-1}(\beta)$ and ${G}^{-1}(\beta)$.
\begin{lemme}[Upper bound on the thresholds]
    \label{lem:wasserstein_infty_bound_main}
    Let Assumption~\ref{ass:continuity_cdf} be satisfied, then for all $\beta \in [K]$
    \begin{align*}
        \abs{G^{-1}(\beta) - \tilde{G}^{-1}(\beta)} \leq \norm{p - \hat p}_{\infty, \Prob_X},\quad \text{almost surely }\Prob^{\otimes n} \otimes \Prob^{\otimes N}_X\enspace.
    \end{align*}
\end{lemme}
The proof of Lemma~\ref{lem:wasserstein_infty_bound_main} uses elementary properties of the generalized inverse functions which are provided in Appendix.
Besides, let us mention, that the difference $\absin{G^{-1}(\beta) - \tilde{G}^{-1}(\beta)}$ resembles the Wasserstein infinity distance which gives an alternative approach to prove Lemma~\ref{lem:wasserstein_infty_bound_main}, see~\citep{Bobkov_Ledoux14}.
Lemma~\ref{lem:wasserstein_infty_bound_main} explains the extra $\log n$ factor that appears in the upper bound, as the minimax estimation in sup norm contains the $\log n$ factor, see for instance~\citep{Stone82,Tsybakov08}.
Another important property of the introduced estimators $\hGamma_{\SE}$ and $\hGamma_{\SSE}$ is obtained via Assumption~\ref{ass:continuity_cdf_estimator}.
It describes the deviation of the information of $\hGamma_{\SE}$ and $\hGamma_{\SSE}$ from the desired level $\beta$.
\begin{prop}[\cite{Denis_Hebiri17}]
    \label{prop:level_of_info_estimators}
    Let $\hat p_k$ for all $k \in [K]$ be arbitrary estimators of the regression functions constructed using $\data_{\floor{n / 2}}$ that satisfies Assumption~\ref{ass:continuity_cdf_estimator}, then there exist constants $C, C' > 0$ such that for all $n, N \in \bbN$ it holds that
    \begin{align*}
        \Exp_{(\data_n, \data_N)}\abs{\beta - \info\pare{\hGamma_{\SE}}} &\leq Cn^{-1/2}\enspace,\\
        \Exp_{(\data_n, \data_N)}\abs{\beta - \info\pare{\hGamma_{\SSE}}} &\leq C'(N + n)^{-1/2}\enspace.
    \end{align*}
\end{prop}
Note that if $\hat p_k$ satisfies Assumption~\ref{ass:continuity_cdf_estimator} for all $k \in [K]$, then $\beta = \info(\tilde\Gamma)$.
This simple fact is a step in the proof of Proposition~\ref{prop:level_of_info_estimators}.
Finally, combination of Lemma~\ref{lem:wasserstein_infty_bound_main}, Proposition~\ref{prop:level_of_info_estimators}, Assumption~\ref{ass:exponential_concentration_estimator} with the peeling argument used in~\citep[Lemma 3.1]{Audibert_Tsybakov07} yields the results of Theorems~\ref{thm:upper_bound_supervised}-\ref{thm:upper_bound_semi_supervised}.

%%%%%%%%%%%%%%%%%%%%%%%%%%%%%%%%%%%%%%%%%%%%%%%%%%%%%%%%%%%%%%%%%%%%%%%%%%%%%%%
    \begin{table}[t!]
\begin{center}
\resizebox{0.49\textwidth}{!}{
\begin{tabular}{l || c || c}
\multicolumn{3}{c}{{$K = 10$}}\\
\hline
\multicolumn{1}{c}{$\beta$} &  \multicolumn{1}{c}{$\beta$-Oracle} & \multicolumn{1}{c}{$\topproc$-$\beta$ Oracle} \\
\hline\noalign{\smallskip}
 2  & 0.05 (0.01) \;\;&\;\; 0.09 (0.01) \\
 5  & 0.00 (0.00) \;\;&\;\;  0.01 (0.00) \\
\hline
\end{tabular}
}

\vspace*{0.25cm}

\resizebox{0.49\textwidth}{!}{
\begin{tabular}{l || c || c}
\multicolumn{3}{c}{$K = 100$} \\
\hline
\multicolumn{1}{c}{$\beta$} &  \multicolumn{1}{c}{$\beta$-Oracle} & \multicolumn{1}{c}{$\topproc$-$\beta$ Oracle} \\
\hline\noalign{\smallskip}
 2  & 0.39 (0.01)  \;\;&\;\;  0.42 (0.01) \\
 5  & 0.20 (0.01)  \;\;&\;\;  0.22 (0.01) \\
10  & 0.09 (0.01)  \;\;&\;\;  0.11 (0.01) \\
20  & 0.03 (0.01)  \;\;&\;\;  0.04 (0.01) \\
\hline
\end{tabular}
}
\caption{\label{table:tablePerfBayes}
For each of the $B = 100$ repetitions and each model, we derive the estimated
errors $\perf_M$ of the $\beta$-Oracle and of the $\topproc$-$\beta$ Oracle \wrt $\beta$.
We compute the means and standard deviations (between parentheses) over the $B = 100$ repetitions.
Top: the data are generated according to $K=10$ -- Bottom: the data are generated according to $K=100$.}

\end{center}
\end{table}
% In this section, we evaluate the plug-in $\beta$-set numerically.
% First, in Section~\ref{subsec:simStudy}, a simulation study is led to support our theory.
% Plug-in $\beta$-set are then applied in Section~\ref{subsec:realData} to real dataset.

%%%%%%%%%%%%%%%%%%%%%%%%%%%%%%%%%%%%%%%%%%
\subsection{Simulation study}
\label{subsec:simStudy}
%%%%%%%%%%%%%%%%%%%%%%%%%%%%%%%%%%%%%%%%%%
The goal of this part is to numerically address the following points.
\begin{itemize}
    \item[1)] Is it more advantageous to go outside of the \emph{classical} multi-class classification settings and consider the confidence set framework?
    To respond to this question we compute the Bayes optimal multi-class classifier and view it as a confidence set with one label.
    We compare this Bayes rule with the $\beta$-Oracle in terms of the error $\perf(\cdot)$ using various values of $\beta \in [K]$ and $K \in \bbN$.
    \item[2)] How does the $\beta$-Oracle confidence set compares to another "Oracle" ($\topproc$-$\beta$ Oracle) which simply includes classes corresponding to the largest values of $p_k(\cdot)$'s?
    \item[3)] Does the proposed plug-in approach indeed gives a good approximation of the $\beta$-Oracle through the error $\perf(\cdot)$ and the information $\info(\cdot)$?
    \item[4)] Despite demonstrating the minimax inconsistency of the top-$\beta$ approach, we wonder whether in some scenarios it can achieve a comparable performance against our semi-supervised plug-in procedure.
\end{itemize}
% We perform a simulation study in order to evaluate the performance of our plug-in procedure.
% We compare our method with the strategy which is based on the top-$\beta$ levels conditional probabilities.
% In the sequel this method is referred as $\topproc$ procedure.
% Additionally, we investigate the influence of the parameter $K$.
% Moreover, we display two measures of performance that are related to all of the risks we investigated in the theoretical part: the error $\perf\pare{\Gamma} = \Prob\pare{Y \notin \Gamma(X)}$ and the information $\info\pare{\Gamma} = \Exp_{\Prob_X}\abs{\Gamma(X)}$.

We consider two simulation schemes depending on the parameter $K \in \{10, 100\}$.
For each $K$, we generate $(X,Y)$ according to a mixture model.
More precisely,
\begin{itemize}
\item[i)] the label $Y$ follows uniform distribution on $[K]$;
\item[ii)] conditional on $Y = k$, the feature $X$ is generated according to a multivariate gaussian distribution with mean $\mu_k \in \mathbb{R}^{10}$ and identity covariance matrix.
\end{itemize}
For each $k \in [K]$, the vectors $\mu_k$ are \iid realizations of uniform distribution on $[0,4]^{10}$.
For this distribution, we have
\begin{equation*}
p_k(X) = \dfrac{f_k(X)}{\sum_{j = 1}^K f_j(X)},
\end{equation*}
where for each $k \in [K]$, $f_k(X)$ is the density function of a multivariate gaussian distribution with mean parameter $\mu_k$ and identity covariance matrix.

For each $K$, the \emph{missclassification} error of the \emph{classical} multi-class classification Bayes rule is evaluated based on a sufficiently large dataset.
It is valued at $0.22$ and at $0.60$ for $K = 10$ and for $K = 100$ respectively.
These values are relatively high, which suggests that confusion is induced by the large number of classes.
Hence, it is reasonable to apply the confidence set approach to this problem.
\begin{table}[t!]
\begin{center}
\resizebox{0.49\textwidth}{!}{
\begin{tabular}{l || c || c}
\multicolumn{1}{c}{$\beta$} &  \multicolumn{1}{c}{$K = 10$} & \multicolumn{1}{c}{$K = 100$} \\
\hline\noalign{\smallskip}
 2  & 2.00 (0.03)  \;\;&\;\;  2.00 (0.03) \\
 5  & 5.00 (0.08)  \;\;&\;\;  5.00 (0.06) \\
10  &  $\cdot$     \;\;&\;\;  10.00 (0.13) \\
20  &  $\cdot$     \;\;&\;\;  20.02 (0.31) \\
\hline
\end{tabular}
}
\caption{\label{table:tableCardBayes}
For each of the $B = 100$ repetitions and each model, we derive the estimated
information levels $\info_M$ of the $\beta$-Oracle set w.r.t. $\beta$.
We compute the means and standard deviations (in parentheses) over the $B = 100$ repetitions.
Left: the data are generated according to $K=10$ -- Right: the data are generated according to $K=100$.}

\end{center}
\end{table}
In the sequel, we aim at providing the estimation of the error of the $\beta$-Oracle.
To this end, for $\beta \in \{2, 5, 10, 20\}$ and each $K$, we repeat $B$ times the following steps.
\begin{enumerate}
\item[i)] simulate two datasets $\mathcal{D}_N$ and $\mathcal{D}_M$ with $N=10000$ and $M = 1000$;
\item[ii)] based on $\mathcal{D}_N$, we compute the empirical counterpart of $G$ and provide an approximation of the $\beta$-Oracle $\Gamma^*_{\beta}$ given in Eq.~\eqref{eq:bayesConfSet} (we recall that
this step requires a dataset which contains only unlabeled features);
\item[iii)] finally, over $\mathcal{D}_{M}$, we compute the empirical counterparts $\perf_M$ (of $\perf(\Gamma^*_{\beta})$) and $\info_M$ (of $\info
(\Gamma_{\beta}^*)$).
\end{enumerate}
From this estimates, we compute the mean and the standard deviation of $\perf_M$ and $\info_M$.
Tables~\ref{table:tablePerfBayes} and~\ref{table:tableCardBayes} present values of the error and of the information which are achieved by the $\beta$-\emph{Oracle} and by the $\topproc$-$\beta$ Oracle.

\begin{table}[t!]
\begin{center}
\resizebox{\textwidth}{!}{
\begin{tabular}{l || ccc || ccc}
\multicolumn{1}{c}{} & \multicolumn{6}{c}{{$K = 10$}}\\
\hline
\multicolumn{1}{c}{} &  \multicolumn{3}{c}{$\hGamma_{\SSE}$} & \multicolumn{3}{c}{$\topproc$-$\beta$} \\
\hline\noalign{\smallskip}
 $\beta$ \;\;\;  & \;\;   \texttt{rforest} \;\;\; &  \texttt{softmax reg} \;\;\;& \texttt{deep learn} & \;\;   \texttt{rforest} \;\;\; &  \texttt{softmax reg} \;\;\;&
 \texttt{deep learn} \\
\noalign{\smallskip}
\hline
\noalign{\smallskip}
 2  & 0.09 (0.01) & 0.06 (0.01) & 0.09 (0.01) \;\;&\;\; 0.13 (0.01) & 0.10 (0.01) & 0.13 (0.02) \\
 5  & 0.01 (0.00) & 0.00 (0.00) & 0.01 (0.00) \;\;&\;\;  0.02 (0.00) & 0.01 (0.00) & 0.02 (0.00) \\
\hline
\end{tabular}
}

% \vspace*{0.25cm}

\resizebox{\textwidth}{!}{
\begin{tabular}{l || ccc || ccc}
\multicolumn{1}{c}{} & \multicolumn{6}{c}{{$K = 100$}}\\
\hline
\multicolumn{1}{c}{} &  \multicolumn{3}{c}{$\hGamma_{\SSE}$} & \multicolumn{3}{c}{$\topproc$-$\beta$} \\
\hline\noalign{\smallskip}
  $\beta$ \;\;\;  & \;\;   \texttt{rforest} \;\;\; &  \texttt{softmax reg} \;\;\;& \texttt{deep learn} & \;\;   \texttt{rforest} \;\;\; &  \texttt{softmax reg} \;\;
  \;&   \texttt{deep learn} \\
\noalign{\smallskip}
\hline
\noalign{\smallskip}
 2  & 0.48 (0.02) & 0.93 (0.01) & 0.46 (0.02) \;\;&\;\;  0.51 (0.01) & 0.96 (0.01) & 0.48 (0.02) \\
 5  & 0.30 (0.02) & 0.85 (0.02) & 0.25 (0.02) \;\;&\;\;  0.31 (0.01) & 0.90 (0.01) & 0.27 (0.01) \\
10  & 0.17 (0.01) & 0.75 (0.02) & 0.12 (0.01) \;\;&\;\;  0.18 (0.01) & 0.80 (0.01) & 0.14 (0.01) \\
20  & 0.07 (0.01) & 0.59 (0.02) & 0.04 (0.01) \;\;&\;\;  0.09 (0.01) & 0.61 (0.02) & 0.06 (0.01) \\
\hline
\end{tabular}
}
\caption{\label{table:tablePerfPlug}
For each of the $B = 100$ repetitions and for each model, we derive the estimated errors~$\perf$ of three different $\hGamma_{\SSE}$'s \wrt $\beta$.
We compute the means and standard deviations (in parentheses) over the $B = 100$ repetitions.
For each $\beta$ and for each $N$, the $\hGamma_{\SSE}$'s, as well as the $\topproc$ procedures are based on, from left to right,
\texttt{rforest}, \texttt{softmax reg} and \texttt{deep learn},
which are respectively the random forest, the softmax regression and the deep learning methods.
Top: the data are generated according to $K=10$ -- Bottom: the data are generated according to $K=100$.}

\end{center}
\end{table}
We now move towards the construction of our semi-supervised plug-in estimators $\hGamma_{\SSE}$.
For each $K$ and each $\beta$, we evaluate the performance of $\hGamma_{\SSE}$ according to three different estimations of the regression function: the $\hat p_k$'s are based on random forests, softmax regression and deep learning procedures.
Let us point out, that for random forests and softmax regression algorithms, the random variables $\hat{p}_k(X)$ appear to be not continuous.
Hence Assumption~\ref{ass:continuity_cdf_estimator} is violated.
To alleviate this issue, we add to $\hat{p}_k(X)$ an independent small perturbation $|\mathcal{N}(0,e^{-10})|$ for simplicity.
The evaluation of the performance of $\hGamma_{\SSE}$ relies on the following steps
\begin{enumerate}
\item[i)] simulate three datasets $\mathcal{D}_n$, $\mathcal{D}_N$ and $\mathcal{D}_M$;
\item[ii)] based on $\mathcal{D}_n$, we compute the estimators $\hat{p}_k$ of $p_k$ according to the considered procedure;
\item[iii)] based on $\mathcal{D}_N$ and $\hat{p}_k$ we compute the function $\hat{G}$ and the estimator $\hGamma_{\SSE}$ as in Eq.~\eqref{eq:IntroSemiSupEstimator} (we recall that
this step requires a dataset which contains only unlabeled features);
\item[iv)] finally, we compute over $\mathcal{D}_M$ the empirical counterpart of $\perf$ and of $\info$ for the considered $\hGamma_{\SSE}$.
\end{enumerate}

%%%%%%%%%%%%%%%%%%%%%%%%%%%%%%%%%%%%%%%%%%%%%%%%%%%%%%%%%%%%%%%%%%%%%%%%%%%%%%%%%%%%%%%%%%%%%%%%%%%%%%%%%%%%%%%%%%%%%%%
%%%%%%%%%%%%%%%%%%%%%%%%%%%%%%%%%%%%%%%%%%%%%%%%%%%%%%%%%%%%%%%%%%%%%%%%%%%%%%%%%%%%%%%%%%%%%%%%%%%%%%%%%%%%%%%%%%%%%%%%%

%%%%%%%%%%%%%%%%%%%%%%%%%%%%%%%%%%%%%%%%%%%%%%%%%%%%%%%%%%%%%%%%%%%%%%%%%%%%%%%%%%%%%%%%%%%%%%%%%%%%%%%%%%%%%%%%%%%%%%%
%%%%%%%%%%%%%%%%%%%%%%%%%%%%%%%%%%%%%%%%%%%%%%%%%%%%%%%%%%%%%%%%%%%%%%%%%%%%%%%%%%%%%%%%%%%%%%%%%%%%%%%%%%%%%%%%%%%%%%%%%

%%%%%%%%%%%%%%%%%%%%%%%%%%%%%%%%%%%%%%%%%%%%%%%%%%%%%%%%%%%%%%%%%%%%%%%%%%%%%%%%%%%%%%%%%%%%%%%%%%%%%%%%%%%%%%%%%%%%%
%%%%%%%%%%%%%%%%%%%%%%%%%%%%%%%%%%%%%%%%%%%%%%%%%%%%%%%%%%%%%%%%%%%%%%%%%%%%%%%%%%%%%%%%%%%%%%%%%%%%%%%%%%%%%%%%%%%%%

\begin{table}[t!]
\begin{center}
\resizebox{\textwidth}{!}{
\begin{tabular}{l || ccc || ccc}
\multicolumn{1}{c}{} & \multicolumn{6}{c}{{$K = 10$}}\\
\hline
\multicolumn{1}{c}{} &  \multicolumn{3}{c}{$N = 100$} & \multicolumn{3}{c}{$N = 10000$} \\
\hline\noalign{\smallskip}
 $\beta$ \;\;\;  & \;\;   \texttt{rforest} \;\;\; &  \texttt{softmax reg} \;\;\;& \texttt{deep learn} & \;\;   \texttt{rforest} \;\;\; &  \texttt{softmax reg} \;\;\;&
 \texttt{deep learn} \\
\noalign{\smallskip}
\hline
\noalign{\smallskip}
 2  & 2.01 (0.09) & 2.01 (0.10) & 2.02 (0.11) \;\;&\;\; 2.00 (0.02) & 2.00 (0.03) & 2.00 (0.03) \\
 5  & 5.02 (0.18) & 4.99 (0.20) & 5.00 (0.21) \;\;&\;\;  5.00 (0.06) & 5.00 (0.08) & 5.00 (0.07) \\
\hline
\end{tabular}
}

\vspace*{0.25cm}

\resizebox{\textwidth}{!}{
\begin{tabular}{l || ccc || ccc}
\multicolumn{1}{c}{} & \multicolumn{6}{c}{{$K = 100$}}\\
\hline
\multicolumn{1}{c}{} &  \multicolumn{3}{c}{$N = 100$} & \multicolumn{3}{c}{$N = 10000$} \\
\hline\noalign{\smallskip}
  $\beta$ \;\;\;  & \;\;   \texttt{rforest} \;\;\; &  \texttt{softmax reg} \;\;\;& \texttt{deep learn} & \;\;   \texttt{rforest} \;\;\; &  \texttt{softmax reg} \;\;
  \;&   \texttt{deep learn} \\
\noalign{\smallskip}
\hline
\noalign{\smallskip}
 2  & 2.02 (0.10) & 2.09 (0.43) & 2.01 (0.09) \;\;&\;\;  2.00 (0.03) & 2.02 (0.15) & 2.00 (0.02) \\
 5  & 4.97 (0.15) & 5.27 (0.70) & 5.01 (0.24) \;\;&\;\;  5.00 (0.04) & 5.01 (0.27) & 5.00 (0.07) \\
10  & 9.98 (0.24) & 10.02 (1.00) & 10.02 (0.42) \;\;&\;\;  10.01 (0.09) & 10.05 (0.32) & 10.00 (0.16) \\
20  & 20.08 (0.48) & 19.74 (0.98) & 20.11 (0.85) \;\;&\;\;  20.00 (0.16) & 20.01 (0.36) & 20.01 (0.28) \\
\hline
\end{tabular}
}
\caption{\label{table:tableCardPlug}
For each of the $B = 100$ repetitions and for each model, we derive the estimated information levels $\info$ of three different $\hGamma_{\SSE}$'s \wrt $\beta$
and the sample size $N$.
We compute the means and standard deviations (in parentheses) over the $B = 100$ repetitions.
For each $\beta$ and each $N$, the $\hGamma_{\SSE}$'s are based on, from left to right,
\texttt{rforest}, \texttt{softmax reg} and \texttt{deep learn},
which are respectively the random forest, the softmax regression and the deep learning procedures.
Top: the data are generated according to $K=10$ -- Bottom: the data are generated according to $K=100$.}

\end{center}
\end{table}

%%%%%%%%%%%%%%%%%%%%%%%%%%%%%%%%%%%%%%%%%%%%%%%%%%%%%%%%%%%%%%%%%%%%%%%%%%%%%%%%%%%%%%%%%%%%%%%%%%%%%%%%%%%%%%%%
%%%%%%%%%%%%%%%%%%%%%%%%%%%%%%%%%%%%%%%%%%%%%%%%%%%%%%%%%%%%%%%%%%%%%%%%%%%%%%%%%%%%%%%%%%%%%%%%%%%%%%%%%%%%%%%%%%%

Again, during these experiments, we compute means and standard deviations.
The parameters $K, n, N$ are fixed as follows: for $K = 10$, we fix $n=1000$ and $N \in \{100,10000\}$; for $K = 100$ we fix $n=10000$ and $N \in \{100,10000\}$.
Finally, the size of $\data_M$ is fixed to $M = 1000$.
The results are illustrated in Tables~\ref{table:tablePerfPlug} and~\ref{table:tableCardPlug}.

As benchmark for the continuation of our experiments, the classical \emph{missclassification} errors of the multi-class classifiers based on random forests, softmax regression and deep learning methods are valued respectively at $0.28$, $0.24$, $0.29$ for $K = 10$, and at $0.65$, $0.98$ $0.63$ for $K = 100$.

Turning to Table~\ref{table:tablePerfBayes} we confirm the intuition that the error of the $\beta$-Oracle decreases as the value of the parameter $\beta$ increases.
Nevertheless, for moderate values of $\beta$, compared to $K$, we obtain a
satisfactory improvement compared to \emph{standard} multi-class classification Bayes rule.
For instance, when $K = 10$ and $\beta = 2$ the error of the $2$-Oracle confidence set is $0.05$, whereas the Bayes classifier has $0.22$; likewise, when $K = 100$ and $\beta = 5$ the the classification error decreases from $0.60$ to $0.20$.
Table~\ref{table:tablePerfBayes} shows that the $\topproc$-$\beta$ Oracle is slightly outperformed by the $\beta$-Oracle in terms of the error, but still performs well.

From Tables~\ref{table:tableCardBayes} and~\ref{table:tableCardPlug}, we observe that the approximation of the information is reasonably good and it gets better with $N$ the number of unlabeled data.
Besides, Tables~\ref{table:tablePerfBayes} and~\ref{table:tablePerfPlug} demonstrate that our algorithm is sensitive to the choice of the underlying estimator $\hat p_k$.
Indeed, when $\hat p_k$ is estimated via the softmax regression, our algorithm fails to give a good approximation to the error of the $\beta$-Oracle.

Table~\ref{table:tablePerfPlug} provides similar conclusions regarding $\hGamma_{\SSE}$, though, unlike the theoretical quantities, there are more scenarios where our method is better than its $\topproc$-$\beta$ counterpart.
Let us point out, that for $K =100$ methods that are based on the softmax regression perform poorly in this setup.
%%%%%%%%%%%%%%%%%%%%%%%%%%%%%%%%%%%%%%%%%%%%%%%%%%%%%%%%%%%%%%%%%%%%%%%%%%%%%%%
    %%%%%%%%%%%%%%%%%%%%%%%%%%%%%%%%%%%%%%%%%%%%%%%%%%%%%%%%%%%%%%%%%%%%%%%%%%%%%%%
%%%%%%%%%%%%%%%%%%%%%%%%%%%%%%%%%%%%%%%%%%%%%%%%%%%%%%%%%%%%%%%%%%%%%%%%%%%%%%%
\section{Discussions}
\label{sec:discussion}
%%%%%%%%%%%%%%%%%%%%%%%%%%%%%%%%%%%%%%%%%%%%%%%%%%%%%%%%%%%%%%%%%%%%%%%%%%%%%%%
%%%%%%%%%%%%%%%%%%%%%%%%%%%%%%%%%%%%%%%%%%%%%%%%%%%%%%%%%%%%%%%%%%%%%%%%%%%%%%%

%%%%%%%%%%%%%%%%%%%%%%%%%%%%%%%%%%%%%%%%%%%%%%%%%%%%%%%%%%%%%%%%%%%%%%%%%%%%%%%
%%%%%%%%%%%%%%%%%%%%%%%%%%%%%%%%%%%%%%%%%%%%%%%%%%%%%%%%%%%%%%%%%%%%%%%%%%%%%%%
\subsection{Around continuity Assumption~\ref{ass:continuity_cdf}}
\label{sub:around_assumption}
%%%%%%%%%%%%%%%%%%%%%%%%%%%%%%%%%%%%%%%%%%%%%%%%%%%%%%%%%%%%%%%%%%%%%%%%%%%%%%%
%%%%%%%%%%%%%%%%%%%%%%%%%%%%%%%%%%%%%%%%%%%%%%%%%%%%%%%%%%%%%%%%%%%%%%%%%%%%%%%
The bedrock of this paper is Assumption~\ref{ass:continuity_cdf}.
Based on it, we ensure that the $\beta$-Oracle confidence set given by Eq.~\eqref{eq:bayesConfSet} is indeed of information $\beta$.
On top of that, the explicit formulation of excess risk in Proposition~\ref{prop:distanceComparison} relies on the continuity of function $G(\cdot)$.
Should Assumption~\ref{ass:continuity_cdf} fail to be satisfied, then there might be no $\beta$-Oracle given by thresholding on some level $\theta \in (0,1)$.
Indeed, assume Assumption~\ref{ass:continuity_cdf} is not satisfied but one can build a $\beta$-Oracle having the form $\Gamma^*_\beta(\cdot) = \enscond{k\in [K]}{p_k(\cdot) > \theta}$ with some $\theta$, then
\begin{align*}
    \beta = \info(\Gamma^*_\beta) = G(\theta)\enspace.
\end{align*}
However, without the continuity, the function $G(\cdot)$ is not surjective and therefore, the equation $G(\theta) = \beta $ may have no solution, which contradicts the fact that $ \info(\Gamma^*_\beta)= \beta $.
Therefore, the settings without the continuity of $G(\cdot)$ deserve a separate study.
Let us also point out that the continuity assumption implies that the $\beta$-Oracle can also be defined as
\begin{align*}
    \Gamma^*_\beta \in \argmin\enscond{\perf\pare{\Gamma}}{\Gamma \in \Upsilon \text { \st } \info(\Gamma) \leq \beta}\enspace,
\end{align*}
where the inequality used in place of the equality.
Indeed, under continuity assumption thanks to Propositions~\ref{prop:minimizer_of_penalized_risk} and~\ref{prop:distanceComparison} we have for all confidence sets $\Gamma$ such that $\info(\Gamma) \leq \beta$
\begin{align*}
    \perf(\Gamma) - \perf(\Gamma^*_\beta) =
    \underbrace{\risk_\beta(\Gamma) - \risk_\beta(\Gamma^*_\beta)}_{\geq 0} + \underbrace{G^{-1}(\beta)\left(\beta - \info(\Gamma)\right)}_{\geq 0}\enspace,
\end{align*}
which implies that the $\beta$-Oracle $\Gamma^*_\beta$ is a minimizer.

%%%%%%%%%%%%%%%%%%%%%%%%%%%%%%%%%%%%%%%%%%%%%%%%%%%%%%%%%%%%%%%%%%%%%%%%%%%%%%%
%%%%%%%%%%%%%%%%%%%%%%%%%%%%%%%%%%%%%%%%%%%%%%%%%%%%%%%%%%%%%%%%%%%%%%%%%%%%%%%
\subsection{Around Lipschitz continuity of \texorpdfstring{$G^{-1}(\cdot)$}{Lg}}
\label{sub:aound_lipschitz}
%%%%%%%%%%%%%%%%%%%%%%%%%%%%%%%%%%%%%%%%%%%%%%%%%%%%%%%%%%%%%%%%%%%%%%%%%%%%%%%
%%%%%%%%%%%%%%%%%%%%%%%%%%%%%%%%%%%%%%%%%%%%%%%%%%%%%%%%%%%%%%%%%%%%%%%%%%%%%%%
Under the assumptions needed in this work, and in particular the continuity
assumption we showed two important facts:
i) no supervised approach can achieve fast rates, that is, faster than $n^{-1/2}$;
ii) some semi-supervised approaches can achieve fast rate.

One might wonder whether extra assumptions on the problem allow a supervised method to get faster rates than $n^{-1/2}$.
We give to this question a partial answer following the recent work of~\cite{Bobkov_Ledoux14} and more precisely their Theorem~5.11.
Applying this result to our framework, we can state that there exists a positive constant $c$ such that
\begin{align*}
    \Exp_{\data_{N}}\abs{G^{-1}(\beta) - G^{-1}_N(\beta)} \leq c \; \Lip(G^{-1})N^{-1/2}\enspace,
\end{align*}
where $\Lip(G^{-1})$ is the Lipschitz constant of $G^{-1}(\cdot)$ and $G^{-1}_N(\cdot)$ is the generalized inverse of
\begin{align*}
    G_N(\cdot) = \frac{1}{N}\sum_{X \in \data_N}\sum_{k = 1}^K \ind{p_k(X) > \cdot}\enspace.
\end{align*}
If, on top of the above, one can show that for any $\alpha > 0$ and for some positive constant~$c'$
\begin{align*}
    \Exp_{\data_{N}}\abs{G^{-1}(\beta) - G^{-1}_N(\beta)}^{1+\alpha} \leq c' \; \Lip^{1+\alpha}(G^{-1})N^{-(1+\alpha)/2}\enspace,
\end{align*}
then under Lipschitz continuity of $G^{-1}(\cdot)$, we can prove that
\begin{align*}
    \excess^{\Exc}_{n, N}(\hGamma_{\square};\class{P}(L, \gamma, \alpha))
    &\lesssim
    {\pare{\frac{n}{\log n}}^{-\frac{(1 + \alpha) \gamma}{2\gamma + d}} }\enspace,
\end{align*}
where $\square$ stands for $\SE$ or $\SSE$.
This would illustrate that both $\hGamma_{\SE}$ and $\hGamma_{\SSE}$ are statistically equivalent under Lipschitz condition on $G^{-1}(\cdot)$, that is, both reach the same rate and the impact of the unlabeled data $\data_N$ is negligible.
We plan to further investigate the influence of this Lipschitz condition on the minimax rates of convergence in our future works.
Since in the present contribution we do not impose this assumption on $G^{-1}(\cdot)$, the upper bound of~\cite{Bobkov_Ledoux14} is not applicable and we had to rely on a different approach.

%%%%%%%%%%%%%%%%%%%%%%%%%%%%%%%%%%%%%%%%%%%%%%%%%%%%%%%%%%%%%%%%%%%%%%%%%%%%%%%
%%%%%%%%%%%%%%%%%%%%%%%%%%%%%%%%%%%%%%%%%%%%%%%%%%%%%%%%%%%%%%%%%%%%%%%%%%%%%%%
\subsection{Around extra logarithm}
\label{sub:around_extra_log}
%%%%%%%%%%%%%%%%%%%%%%%%%%%%%%%%%%%%%%%%%%%%%%%%%%%%%%%%%%%%%%%%%%%%%%%%%%%%%%%
%%%%%%%%%%%%%%%%%%%%%%%%%%%%%%%%%%%%%%%%%%%%%%%%%%%%%%%%%%%%%%%%%%%%%%%%%%%%%%%
Theorems~\ref{thm:lower_bound_semi_supervised} and~\ref{thm:upper_bound_semi_supervised} demonstrate that for the excess risk and the Discrepancy, the upper and the lower bounds differ by a logarithmic factor.
As we have already pointed out, this factor appears in the upper bounds due to Lemma~\ref{lem:wasserstein_infty_bound_main} which relates the difference between two thresholds to the infinity norm.
One might hope that if we manage to replace the infinity norm by any other $\ell_q$-norm on the right hand side of the inequality in Lemma~\ref{lem:wasserstein_infty_bound_main} this logarithm can be eliminated.
Unfortunately, it appears that this bound is actually tight, in a sense that one can construct a distribution $\Prob$ and an estimator $\hat p_k$ for all $k \in [K]$ such that an equality is achieved in Lemma~\ref{lem:wasserstein_infty_bound_main}.
These arguments suggest that the obtained upper bound should be optimal.
They also imply that the lower bounds could be further refined to get an extra logarithmic factor.
%while the lower might be made more accurate.
%     the problem appears in the lower bound.
%We believe that our way of providing the lower bounds is sub-optimal and could be further refined to get the extra logarithmic factor.
Let us also mention that the continuity Assumption~\ref{ass:continuity_cdf} in combination with the margin Assumption~\ref{ass:margin_assumption} are main obstacles that did not allow us to provide better lower bounds.
Nevertheless, our proofs are already involved and our results allow to make non trivial conclusions even without going into the details concerning the logarithms.

%%%%%%%%%%%%%%%%%%%%%%%%%%%%%%%%%%%%%%%%%%%%%%%%%%%%%%%%%%%%%%%%%%%%%%%%%%%%%%%
%%%%%%%%%%%%%%%%%%%%%%%%%%%%%%%%%%%%%%%%%%%%%%%%%%%%%%%%%%%%%%%%%%%%%%%%%%%%%%%
\section{Conclusion}
\label{sec:conclusion}
%%%%%%%%%%%%%%%%%%%%%%%%%%%%%%%%%%%%%%%%%%%%%%%%%%%%%%%%%%%%%%%%%%%%%%%%%%%%%%%
%%%%%%%%%%%%%%%%%%%%%%%%%%%%%%%%%%%%%%%%%%%%%%%%%%%%%%%%%%%%%%%%%%%%%%%%%%%%%%%
In this work we have studied the minimax settings of confidence set multi-class classification.
First of all, following previous works we have shown that a top-$\beta$ type procedure is inconsistent in our settings and more involved estimators should be proposed.
Besides, we have demonstrated that no supervised estimator can achieve rates that are faster than $n^{-1/2}$, which stays in contrast with other classical settings.
Additionally, we have shown that fast rates are achievable for semi-supervised techniques provided that the size of the unlabeled sample is large enough.
Consecutively, we have established that our lower bounds are either optimal or nearly optimal by providing a supervised and a semi-supervised estimators which are tractable in practice.
Our future works shall be focused on the Lipschitz condition of $G^{-1}(\cdot)$ discussed in Section~\ref{sub:aound_lipschitz}, in particular, we want to understand how this extra assumption affects our lower bounds.

%%%%%%%%%%%%%%%%%%%%%%%%%%%%%%%%%%%%%%%%%%%%%%%%%%%%%%%%%%%%%%%%%%%%%%%%%%%%%%%
    \vskip 0.2in
    \bibliographystyle{imsart-nameyear}
    \bibliography{bibli}
%%%%%%%%%%%%%%%%%%%%%%%%%%%%%%%%%%%%%%%%%%%%%%%%%%%%%%%%%%%%%%%%%%%%%%%%%%%%%%%
    \appendix

%%%%%%%%%%%%%%%%%%%%%%%%%%%%%%%%%%%%%%%%%%%%%%%%%%%%%%%%%%%%%%%%%%%%%%%%%%%%%%%
%%%%%%%%%%%%%%%%%%%%%%%%%%%%%%%%%%%%%%%%%%%%%%%%%%%%%%%%%%%%%%%%%%%%%%%%%%%%%%%
\section{Technical results}
\label{app:technical_results}
%%%%%%%%%%%%%%%%%%%%%%%%%%%%%%%%%%%%%%%%%%%%%%%%%%%%%%%%%%%%%%%%%%%%%%%%%%%%%%%
%%%%%%%%%%%%%%%%%%%%%%%%%%%%%%%%%%%%%%%%%%%%%%%%%%%%%%%%%%%%%%%%%%%%%%%%%%%%%%%
Here we provide proofs for our result.
This Appendix is composed of the following part: in Appendix~\ref{app:technical_results} we introduce some technical results used for our proofs;
Appendix~\ref{sec:upper_bounds_appendix} is devoted to the proofs of the upper bounds; Appendix~\ref{sec:proof_of_the_lower_bounds} provides with the proofs our our main lower bounds; finally, in Appendix~\ref{app:fail_naive_top_k} we prove the inconsistency of top-$\beta$ approaches.

In this section we gather several technical results which are used to derive the contributions of this work.
Let us start by introducing notation used in the appendix.
Given any two probability measures $\Prob_1, \Prob_2$ on some space measurable space $(\class{X}, \class{A})$ the Kullback–Leibler divergence between $\Prob_1$ and $\Prob_2$ is defined as
\begin{align}
    \label{eq:KL_definitions}
    \KL(\Prob_1, \Prob_2) \eqdef \begin{cases}
                                \int_{\class{X}} \log\pare{\frac{d\Prob_1}{d\Prob_2}}d\Prob_1,&\quad \supp(\Prob_1) \subset \supp(\Prob_2)\\
                                +\infty, &\quad \text{otherwise}
                            \end{cases}\enspace,
\end{align}
and the total variation distance is defined as
\begin{align}
    \label{eq:TV_definition}
    \TV(\Prob_1, \Prob_2) \eqdef \sup_{A \in \class{A}} \abs{\Prob_1(A) - \Prob_2(A)}\enspace.
\end{align}
We start with Fano's inequality in the form proved by~\citep{Birge05}.
\begin{lemme}[Fano's inequality~\citep{Birge05}]
    \label{lem:birges_for_lower_bound}
    Let $\ens{\Prob_i}_{i = 0}^m$ be a finite family of probability measures on $(\class{X}, \class{A})$ and let $\ens{A_i}_{i = 0}^m$ be a finite family of disjoint events such that $A_i \in \class{A}$ for each $i = 0, \ldots, m$.
    Then,
    \begin{align*}
        \min_{i \in \{0, 1, \ldots, m\}} \Prob_i(A_i) \leq \pare{0.71 \bigvee \frac{\frac{1}{m}\sum_{i = 1}^m \KL(\Prob_i, \Prob_0)}{\log(m + 1)}}\enspace.
    \end{align*}
\end{lemme}
\begin{lemme}[Pinsker's inequality]
    \label{lem:pinsker_inequality}
    Given any two probability measures $\Prob_1, \Prob_2$ on some measurable space $(\class{X}, \class{A})$ we have
    \begin{align*}
        \TV(\Prob_1, \Prob_2) \leq \sqrt{\frac{1}{2}\KL(\Prob_1, \Prob_2)}\enspace.
    \end{align*}
\end{lemme}
\begin{lemme}[Hoeffding's inequality~\citep{Hoeffding63}]
    Let $b >0$ be a real number, and $N$ be a positive integer. Let $X_1,\ldots, X_N$ be $N$ random variables having values in $[0,b]$, then
\begin{align*}
     \Prob \left( \left|  \frac{1}{N}  \sum_{i=1}^N     \left( X_i - \mathbb{E}\left[ X_i \right]  \right) \right|  \geq t   \right)
 \leq
 2 \exp\left(- \frac{2 N t^2 }{b^2}   \right),\qquad \forall t > 0\enspace.
\end{align*}
\end{lemme}
\begin{prop}[Properties of the generalized inverse]
    \label{prop:inverse_G_properties}
    Let $X \in \bbR^d$ and $\Prob_X$ be a Borel measure on $\bbR^d$, let $p: \bbR^d \rightarrow [0, 1]^K$ be a vector function, we define for all $t \in [0, 1]$ and all $\beta \in (0, K)$
    \begin{align*}
        G(t) &\eqdef \sum_{k = 1}^K\Prob_X(p_k(X) > t),\quad
        G^{-1}(\beta) \eqdef \inf\enscond{t \in [0, 1]}{G(t) \leq \beta}\enspace.
    \end{align*}
    Then,
    \begin{itemize}
        \item  for all $t \in (0, 1)$ and $\beta \in (0, K)$ we have
        $G^{-1}(\beta) \leq t \iff G(t) \leq \beta$.
        \item if for all $k \in [K]$ the mappings $t \mapsto \Prob_X(p_k(X) > t)$ are continuous on $(0, 1)$, then
        \begin{itemize}
            \item for all $\beta \in (0, K)$ we have $G(G^{-1}(\beta)) = \beta$\enspace.
        \end{itemize}
    \end{itemize}
\end{prop}
The next result is an analogue of the classical inverse transform theorem~\citep[Lemma 21.1]{Vaart98} and was already established by~\cite{Denis_Hebiri17}.
\begin{lemme}
\label{lem:uniform_inverse_cdf_analogue}
Let $\varepsilon$ distributed from a uniform distribution on $[K]$
and $Z_1, \ldots, Z_K$, $K$ real valued random variables independent from $\varepsilon$, such that the function $t \mapsto H(t)$ defined as
\begin{align*}
    H(t) \eqdef \frac{1}{K}\sum_{k = 1}^K \Prob(Z_k \leq t)\enspace,
\end{align*}
is continuous.
Consider random variable $Z = \sum_{k = 1}^K Z_k \1_{\{\varepsilon = k\}}$ and let $U$ be distributed according to the uniform distribution on $[0,1]$.
Then
\begin{align*}
H(Z)  \overset{\mathcal{L}}{=} U \;\; {\rm and} \;\; H^{-1}(U)  \overset{\mathcal{L}}{=} Z\enspace,
\end{align*}
where $H^{-1}$ denotes the generalized inverse of $H$.
\end{lemme}
\begin{proof}
First we note that for every $t \in [0,1]$,
$\Prob\left(H(Z) \leq t\right) = \Prob\left(Z \leq H^{-1}(t)\right)$.
Moreover, we have
\begin{align*}
\Prob\left(H(Z) \leq t\right) & = \sum_{k = 1}^K \Prob(Z \leq H^{-1}(t), \varepsilon = k) \\
                                & = \dfrac{1}{K} \sum_{k = 1}^K \Prob(Z_k \leq H^{-1}(t)) \;\; &({\rm with} \;\;
                                                                \varepsilon \;\; {\rm independent \;\; of} \;\; Z)\enspace\\
                                   & =H(H^{-1}(t)) \\
                                   & =t \;\; &({\rm with} \;\; H \;\; {\rm continuous})\enspace.
\end{align*}
To conclude the proof, we observe that
\begin{align*}
\Prob\left(H^{-1}(U) \leq t\right) & = \Prob\left(U \geq H(t)\right)
                                        = \dfrac{1}{K} \sum_{k = 1}^K \Prob(Z_k \leq t) \\
                                        & = \sum_{k = 1}^K \Prob\left(Z_k \leq t, \varepsilon  = k\right)
                                         =\Prob(Z\leq t)\enspace.
\end{align*}
\end{proof}

%%%%%%%%%%%%%%%%%%%%%%%%%%%%%%%%%%%%%%%%%%%%%%%%%%%%%%%%%%%%%%%%%%%%%%%%%%%%%%%
%%%%%%%%%%%%%%%%%%%%%%%%%%%%%%%%%%%%%%%%%%%%%%%%%%%%%%%%%%%%%%%%%%%%%%%%%%%%%%%
\section{Upper bounds}
\label{sec:upper_bounds_appendix}
%%%%%%%%%%%%%%%%%%%%%%%%%%%%%%%%%%%%%%%%%%%%%%%%%%%%%%%%%%%%%%%%%%%%%%%%%%%%%%%
%%%%%%%%%%%%%%%%%%%%%%%%%%%%%%%%%%%%%%%%%%%%%%%%%%%%%%%%%%%%%%%%%%%%%%%%%%%%%%%
In this section we prove Theorems~\ref{thm:upper_bound_supervised} and~\ref{thm:upper_bound_semi_supervised}.
It will be clear from our analysis that the proof of Theorem~\ref{thm:upper_bound_supervised} follows directly from Theorem~\ref{thm:upper_bound_semi_supervised} by setting $N = 0$ in the statement of Theorem~\ref{thm:upper_bound_semi_supervised}.
Thus, in this section for simplicity we omit the subscript $\SSE$ from $\hGamma_{\SSE}$.
Recall that our dataset consists of three parts $\data_{\floor{n / 2}}, \data_{\ceil{n / 2}}, \data_N$.
The set $\data_{\floor{n / 2}}$ is used to construct an estimator $\hat p$ of the regression function $p$, that is, $\hat p$ is independent from both $\data_{\ceil{n / 2}}, \data_N$.
The other two sets $\data_{\ceil{n / 2}}, \data_N$ are used in a semi-supervised manner to estimate the threshold, that is, we erase the labels from $\data_{\ceil{n / 2}}$.
Let $\beta\in[K-1]$, and also recall the definition of the proposed semi-supervised estimator for a given $x\in\bbR^d$
\begin{align*}
    \hGamma(x) = \enscond{k \in [K]}{\hat p_k(x) \geq \hat G^{-1}(\beta)}\enspace,
\end{align*}
with $\hat p_k(x)$ satisfying Assumptions~\ref{ass:continuity_cdf_estimator},~\ref{ass:exponential_concentration_estimator} for all $k \in [K]$.
Moreover, $\hat G^{-1}(\beta)$ defined as the generalized inverse of
\begin{align*}
    \hat G(t)
    =
    \frac{1}{\ceil{n / 2} + N}\sum_{X \in \data_N \bigcup \data_{\ceil{n / 2}}}\sum_{k = 1}^K\1_{\hat p_k(X) > t}\enspace,
\end{align*}
where $t\in[0,1]$.
Additionally, recall that the $\beta$-Oracle is given as
\begin{align}
    \Gamma_\beta^* (x)  =  \enscond{k \in [K]}{p_{k}(x)\geq G^{-1}(\beta)}\enspace,
\end{align}
where $G^{-1}(\cdot)$ is the generalized inverse of
\begin{align*}
    G(t) \eqdef \sum_{k = 1}^K \Prob(p_k(X) \geq t)\enspace.
\end{align*}
Lastly, let us re-introduce an idealized version $\tGamma$ of the proposed estimator $\hGamma$ which 'knows' the marginal distribution $\Prob_X$ of the feature vector $X \in \bbR^d$ as
\begin{align*}
    \tGamma(x) = \enscond{k \in [K]}{\hat p_k(x) \geq \tilde G^{-1}(\beta)}\enspace,
\end{align*}
with $\tilde{G} \eqdef \sum_{k = 1}^K\Prob_X(\hat p_k(X) > t)$, conditionally on the data.
The following result is needed to relate the threshold $\tilde G^{-1}(\beta)$ of $\tGamma$ to the true value of the threshold $G^{-1}(\beta)$.
\begin{lemme}[Upper-bound on the thresholds]
    \label{lem:wasserstein_infty_bound_appendix}
    Let $X \in \bbR^d$ and $\Prob_X$ be a Borel measure on $\bbR^d$.
    For two vector functions $p, \hat p: \bbR^d \rightarrow [0, 1]^K$, we define
    \begin{align*}
        G(\cdot) \eqdef \sum_{k = 1}^K\Prob_X(p_k(X) > \cdot), \quad \tilde{G}(\cdot) \eqdef \sum_{k = 1}^K\Prob_X(\hat p_k(X) > \cdot)\enspace.
    \end{align*}
    If for all $k \in [K]$ the mapping $t \mapsto \Prob_X(p_k(X) > t)$ is continuous on $(0, 1)$, then for every $\beta \in (0, K)$
    \begin{align*}
        \abs{G^{-1}(\beta) - \tilde{G}^{-1}(\beta)} \leq \norm{\hat p - p}_{\infty, \Prob_X}\enspace.
    \end{align*}
\end{lemme}
\begin{proof}
    The proof of this result is very similar to the proof of~\cite[Theorem 2.12]{Bobkov_Ledoux14}.
    We start by defining the following quantity
    \begin{align*}
        h^* = \inf\enscond{h \geq 0}{\forall t \in [0, 1]\,\, \tilde{G}(t + h) \leq G(t) \leq \tilde{G}(t - h)}\enspace.
    \end{align*}
    Due to the definition of $h^*$ we have that for all $t \in [0, 1]$
    \begin{align*}
        \tilde{G}(t + h^*) \leq G(t) \leq \tilde{G}(t - h^*)\enspace,
    \end{align*}
    that is, applying Proposition~\ref{prop:inverse_G_properties} to the second inequality we get for all $t \in [0, 1]$
    \begin{align*}
        t - h^* \leq \tilde{G}^{-1}(G(t))\enspace,
    \end{align*}
    thus, for $t = G^{-1}(\beta)$ with $\beta \in (0, K)$ thanks to Proposition~\ref{prop:inverse_G_properties} we get
    \begin{align*}
        G^{-1}(\beta) -  \tilde{G}^{-1}(\beta) \leq h^*\enspace.
    \end{align*}
    The inequality $\tilde{G}^{-1}(\beta) - G^{-1}(\beta) \leq h^*$ is obtained in the same way.
    Thus, we have proved that
    \begin{align*}
        \abs{G^{-1}(\beta) - \tilde{G}^{-1}(\beta)} \leq h^*\enspace.
    \end{align*}
    Finally, notice that for all $t \in [0, 1]$
    \begin{align*}
       \underbrace{\sum_{k = 1}^K\Prob_X\pare{\hat p_k(X) > t + \norm{\hat p - p}_{\infty, \Prob_X}}}_{\tilde{G}\pare{t + \norm{\hat p - p}_{\infty, \Prob_X}}}
       &\leq
       \underbrace{\sum_{k = 1}^K\Prob_X(p_k(X) > t)}_{G(t)}\\
       &\leq
       \underbrace{\sum_{k = 1}^K\Prob_X\pare{\hat p_k(X) > t - \norm{\hat p - p}_{\infty, \Prob_X}}}_{\tilde{G}\pare{t - \norm{\hat p - p}_{\infty, \Prob_X}}}\enspace,
    \end{align*}
    where we used the fact that for all $k \in [K]$
    \begin{align*}
       \Prob_X(\hat p_k(X) > t + \abs{\hat p_k(X) - p_k(X)})
       &\leq
       \Prob_X(p_k(X) > t)\\
       &\leq
       \Prob_X(\hat p_k(X) > t - \abs{\hat p_k(X) - p_k(X)})\enspace,
    \end{align*}
    and $\Prob_X\pare{\abs{\hat p_k(X) - p_k(X)} \leq \norm{\hat p - p}_{\infty, \Prob_X}} = 1$.
    Therefore by definition of $h^*$, we can write $h^* \leq \norm{\hat p - p}_{\infty, \Prob_X}$ and we conclude.
\end{proof}

We are in position to prove Theorem~\ref{thm:upper_bound_semi_supervised}, let us point out that the most difficult part in Theorem~\ref{thm:upper_bound_semi_supervised} is the upper-bound on the excess risk.
The upper-bound on the discrepancy follows the same arguments as the ones we use for the excess-risk.

{\bf Excess risk and discrepancy:} to upper-bound the excess risk we first separate it into two parts as
\begin{align*}
    \risk_\beta(\hGamma) - \risk_\beta(\Gamma^*_\beta)
    &=
    \underbrace{\pare{\risk_\beta(\tGamma) - \risk_\beta(\Gamma^*_\beta)}}_{R_1} + \underbrace{\pare{\risk_\beta(\hGamma) - \risk_\beta(\tGamma)}}_{R_2}\enspace.
\end{align*}
Recall that thanks to Proposition~\ref{prop:distanceComparison} we have
\begin{align*}
    R_1
    &=
    \sum_{k =1}^K \Exp\left[|p_k(X)-G^{-1}(\beta)|\1_{\{k \in \tGamma(X) \triangle \Gamma^{*}_{\beta}(X)\}}\right]\enspace.
\end{align*}
Moreover, let us point out that if some $k \in \tGamma(X) \triangle \Gamma^{*}_{\beta}(X)$ then either
\begin{align*}
    \begin{cases}
        p_k(X) - G^{-1}(\beta) \geq 0\\
        \hat p_k(X) - \tilde G^{-1}(\beta) < 0
    \end{cases}\quad\text{or}\quad
    \begin{cases}
        p_k(X) - G^{-1}(\beta) < 0\\
        \hat p_k(X) - \tilde G^{-1}(\beta) \geq 0
    \end{cases}\enspace,
\end{align*}
holds.
Thus on the event $k \in \tGamma(X) \triangle \Gamma^{*}_{\beta}(X)$ we have
\begin{align*}
    \abs{p_k(X) - G^{-1}(\beta)}
    &\leq
    \abs{\hat p_k(X) -  p_k(X) + G^{-1}(\beta) - \tilde G^{-1}(\beta)}\\
    &\leq
    \abs{\hat p_k(X) -  p_k(X)} + \abs{G^{-1}(\beta) - \tilde G^{-1}(\beta)}\enspace.
\end{align*}
Therefore, for $R_1$ using Lemma~\ref{lem:wasserstein_infty_bound_appendix} and the observations above we can write
\begin{align*}
    R_1
    &\leq
    \sum_{k =1}^K \Exp\left[|p_k(X)-G^{-1}(\beta)|\ind{\abs{p_k(X) - G^{-1}(\beta)} \leq \abs{\hat p_k(X) -  p_k(X)} + \abs{G^{-1}(\beta) - \tilde G^{-1}(\beta)}}\right]\\
    &\leq
    \sum_{k =1}^K \Exp\left[|p_k(X)-G^{-1}(\beta)|\ind{\abs{p_k(X) - G^{-1}(\beta)} \leq 2\norm{\hat p -  p}_{\infty, \Prob_X}}\right]\\
    &\leq
    \sum_{k =1}^K \Exp\left[2\norm{\hat p - p}_{\infty, \Prob_X}\ind{\abs{p_k(X) - G^{-1}(\beta)} \leq 2\norm{\hat p - p}_{\infty, \Prob_X}}\right]\\
    &=
    2\norm{\hat p -  p}_{\infty, \Prob_X}\sum_{k =1}^K\Prob_X\pare{\abs{p_k(X) - G^{-1}(\beta)} \leq 2\norm{\hat p -  p}_{\infty, \Prob_X}}\enspace,
\end{align*}
finally, using the margin Assumption~\ref{ass:margin_assumption} we get almost surely data
\begin{align*}
    R_1 \leq c_12^{1 + \alpha}K\norm{\hat p - p}^{1 + \alpha}_{\infty, \Prob_X}\enspace.
\end{align*}
Integrating over the data from both sides and using Assumption~\ref{ass:exponential_concentration_estimator} we get
\begin{align*}
    \Exp_{(\data_n, \data_N)}R_1 \leq c_1C_22^{1 + \alpha}K\pare{{\frac{n}{\log n}}}^{-\frac{(1 + \alpha)\gamma}{2\gamma + d}}\enspace.
\end{align*}
For $R_2$ the following trivial upper-bound holds
\begin{align}
    \label{eq:R_2_as_hamming}
    R_2
    &=
    \pare{\perf(\hGamma) - \perf(\tGamma)} + G^{-1}(\beta)\pare{\info(\hGamma) - \info(\tGamma)}\nonumber\\
    &=
    \sum_{k = 1}^K \Exp\pare{p_k(X) - G^{-1}(\beta)}\pare{\ind{k \in \tGamma(X)} - \ind{k \in \hGamma(X)}}\nonumber\\
    &\leq
    \underbrace{\sum_{k = 1}^K \Exp\abs{\ind{k \in \tGamma(X)} - \ind{k \in \hGamma(X)}}}_{\Exp\abs{\tGamma(X) \triangle \hGamma(X)}}\\
    &=
    \sum_{k = 1}^K \Exp\abs{\ind{\hat p_k(X) \geq \hat G^{-1}(\beta)} - \ind{\hat p_k(X) \geq \tilde G^{-1}(\beta)}}\enspace,\nonumber
\end{align}
now, thanks to the first property of Proposition~\ref{prop:inverse_G_properties} we can write
\begin{align*}
    R_2
    &\leq
    \sum_{k = 1}^K \Exp\abs{\ind{\hat G(\hat p_k(X)) \leq \beta} - \ind{\tilde G(\hat p_k(X)) \leq \beta}}\\
    &\leq
    \sum_{k = 1}^K\Prob_X\pare{\abs{\hat G(\hat p_k(X)) - \tilde G(\hat p_k(X))} \geq \abs{\tilde G(\hat p_k(X)) - \beta}}
\end{align*}
To finish our proof we make use of the peeling technique of~\citep[Lemma 3.1]{Audibert_Tsybakov07}.
That is, we define for $\delta > 0$ and $k \in [K]$
\begin{align*}
    {A}^k_0
    &=
    \ens{\abs{\tilde{G}(\hat{p}_k(X)) -\beta} \leq \delta}\\
    {A}^k_j
    &=
    \ens{2^{j-1}\delta  < \abs{\tilde{G}(\hat{p}_k(X)) -\beta} \leq 2^{j} \delta } , \;\; j \geq 1.
\end{align*}
Since, for every $k \in [K]$, the events $({A}^k_j)_{j\geq 0}$ are mutually exclusive, we deduce
\begin{align}
\label{eq:eqDec}
    \sum_{k = 1}^K &\Prob_X\pare{\abs{\hat{G}(\hat{p}_k(X)) - \tilde{G}(\hat{p}_k(X))} \geq \abs{\tilde{G}(\hat{p}_k(X)) -\beta}}
    =
    \\
    &\sum_{k = 1}^K \sum_{j\geq 0} \Prob_X\pare{\abs{\hat{G}(\hat{p}_k(X)) - \tilde{G}(\hat{p}_k(X))} \geq \abs{\tilde{G} (\hat{p}_k(X)) -\beta} , A_j^k}\nonumber\enspace.
\end{align}
Now, we consider $\varepsilon$ uniformly distributed on $[K]$ independent of the data and $X$.
Conditional on the data and under Assumption~\ref{ass:continuity_cdf_estimator}, we apply Lemma~\ref{lem:uniform_inverse_cdf_analogue} with $Z_k = \hat{p}_k(X)$, $Z = \sum_{k = 1}^K Z_k \1_{\{\varepsilon = k\}}$ and then obtain that
$\tilde{G}(Z)$ is uniformly distributed on $[0,K]$.
Therefore, for all $j \geq 0$ and $\delta > 0$, we deduce
\begin{align*}
    \dfrac{1}{K} \sum_{k=1}^K \Prob_X\left(|\tilde{G}(\hat{p}_k(X)) - \beta | \leq 2^j \delta\right)
    &=
    \Prob_X\left(|\tilde{G}(Z) - \beta | \leq 2^j \delta\right)
    \leq
    \dfrac{2^{j+1} \delta}{K}\enspace.
\end{align*}
Hence, for all $j \geq 0$, we obtain
\begin{equation}
    \label{eq:eqEns}
    \sum_{k = 1}^K \Prob_X(A_j^k) \leq \sum_{k=1}^K \Prob_X\left(|\tilde{G}(\hat{p}_k(X)) - \beta | \leq 2^j \delta\right) \leq 2^{j+1} \delta\enspace.
\end{equation}
Next, we observe that for all $j \geq 1$
\begin{align}
    \label{eq:eqAux}
    \sum_{k = 1}^K \Prob_X\left(|\hat{G}(\hat{p}_k(X)) - \tilde{G}(\hat{p}_k(X))| \geq |\tilde{G} (\hat{p}_k(X)) -\beta| , A_j^k\right)
    \leq
    \\
    \sum_{k = 1}^K \Prob_X\left(|\hat{G}(\hat{p}_k(X)) - \tilde{G}(\hat{p}_k(X))| \geq 2^{j - 1}\delta , A_j^k\right)\nonumber\enspace.
\end{align}
Thus, we obtain that
\begin{align*}
    R_2
    \leq
    \sum_{k = 1}^K \sum_{j \geq 0}\Prob_X\left(|\hat{G}(\hat{p}_k(X)) - \tilde{G}(\hat{p}_k(X))| \geq 2^{j - 1}\delta , A_j^k\right)\enspace,
\end{align*}
almost surely data.
Integrating from both sides with respect to the data we get
\begin{align*}
    \Exp_{(\data_n, \data_N)}R_2
    &\leq
    \sum_{k = 1}^K \sum_{j \geq 0}\Exp_{(\data_n, \data_N)}\Prob_X\left(|\hat{G}(\hat{p}_k(X)) - \tilde{G}(\hat{p}_k(X))| \geq 2^{j - 1}\delta , A_j^k\right)\\
    &=
    \sum_{k = 1}^K \sum_{j \geq 0}\Exp_{(\data_{\floor{n / 2}}, \data_{\ceil{n / 2}}, \data_N, X\sim \Prob_X)}\ind{|\hat{G}(\hat{p}_k(X)) - \tilde{G}(\hat{p}_k(X))| \geq 2^{j - 1}\delta}\ind{A_j^k}\enspace.
\end{align*}
recall that the function $\ind{A_j^k}$ for all $j \geq 0$ and $k \in [K]$ is independent from $\data_{\ceil{n / 2}}, \data_N$, thus we can write
\begin{align*}
    &\Exp_{(\data_{\floor{n / 2}}, \data_{\ceil{n / 2}}, \data_N, X\sim \Prob_X)}\ind{|\hat{G}(\hat{p}_k(X)) - \tilde{G}(\hat{p}_k(X))| \geq 2^{j - 1}\delta}\ind{A_j^k}
    =
    \\
    &\Exp_{(\data_{\floor{n / 2}}, X \sim \Prob_X)}\Exp_{(\data_{\ceil{n / 2}}, \data_N)}\left[\ind{|\hat{G}(\hat{p}_k(X)) - \tilde{G}(\hat{p}_k(X))| \geq 2^{j - 1}\delta} \Big\rvert \data_{\floor{n / 2}}, X\right]\ind{A_j^k}]
\end{align*}

Now, since conditional on $(\data_{\floor{n / 2}}, X)$,  $\hat{G}(\hat{p}_k(X))$ is an empirical mean of \iid random variables of common mean $\tilde{G}(\hat{p}_k(X)) \in [0,K]$, we deduce from Hoeffding's inequality that
\begin{align*}
    \Exp_{(\data_{\ceil{n / 2}}, \data_N)}\left[\ind{|\hat{G}(\hat{p}_k(X)) - \tilde{G}(\hat{p}_k(X))| \geq 2^{j - 1}\delta} \Big\rvert \data_{\floor{n / 2}}, X\right]
    \leq
    2e^{-\dfrac{(N + \ceil{n / 2})\delta^2 2^{2j-1}}{K^2}}\enspace.
\end{align*}
Therefore, treating $A_0^k$ separately, we get from inequalities of Eqs.~\eqref{eq:eqDec},~\eqref{eq:eqEns}, and~\eqref{eq:eqAux}
\begin{align*}
    \Exp_{(\data_n, \data_N)}R_2
    \leq
    2\delta + \delta \sum_{j \geq 1} 2^{j+2} \exp\left(-\dfrac{(N + \ceil{n / 2})\delta^2 2^{2j-1}}{K^2}\right)\enspace.
\end{align*}
Finally, choosing $\delta = \dfrac{K}{\sqrt{N + \ceil{n / 2}}}$ in the above inequality, we finish the proof.

{\bf Hamming risk:} here we provide an upper bound on the Hamming risk.
First, by the triangle inequality we can write for the proposed estimator $\hGamma$ and the pseudo Oracle $\beta$ set $\tGamma$
\begin{align*}
    \Exp_{(\data_n, \data_N)} \Exp_{X \sim \Prob_X}\abs{\hGamma(X) \triangle \Gamma^*_\beta(X)}
    \leq
    &\Exp_{(\data_n, \data_N)} \Exp_{X \sim \Prob_X}\abs{\tGamma(X) \triangle \Gamma^*_\beta(X)}\\
    &+
    \Exp_{(\data_n, \data_N)} \Exp_{X \sim \Prob_X}\abs{\hGamma(X) \triangle \tGamma(X)}\enspace.
\end{align*}
Notice that for the term $\Exp_{(\data_n, \data_N)} \Exp_{X \sim \Prob_X}\abs{\hGamma(X) \triangle \tGamma(X)}$ we can re-use the proof technique used for the term $R_2$ in Eq.~\eqref{eq:R_2_as_hamming}.
Thus, it remain to upper-bound the term $\Exp_{(\data_n, \data_N)} \Exp_{X \sim \Prob_X}\abs{\tGamma(X) \triangle \Gamma^*_\beta(X)}$.
The proof on this part closely follows the machinery used in~\cite{Denis_Hebiri17}, however, let us mention that they used this method to obtain a bound on the Discrepancy which leads to a sub-optimal rate.
Nevertheless, their approach gives a correct rate if instead of the Discrepancy we bound the Hamming distance.
For the sake of completeness we write the principal parts of the proof here.

First of all, by the definition of sets $\Gamma^*_\beta$ and $\tGamma$ we can write for $(*) = \Exp_{X \sim \Prob_X}\abs{\tGamma(X) \triangle \Gamma^*_\beta(X)}$
\begin{align*}
    (*) = \sum_{k = 1}^K\Exp_{X \sim \Prob_X}\abs{\ind{\hat{p}_k(X) \geq \tilde{G}^{-1}(\beta)} - \ind{{p}_k(X) \geq {G}^{-1}(\beta)}}\enspace,
\end{align*}
Now if $\hat{p}_k(X) \geq \tilde{G}^{-1}(\beta)$ and ${p}_k(X) < {G}^{-1}(\beta)$ we can have the following situations
\begin{itemize}
    \item if $\tilde{G}^{-1}(\beta) > {G}^{-1}(\beta)$, then $\abs{{p}_k(X) - {G}^{-1}(\beta)} \leq \abs{\hat{p}_k(X) - p_k(X)}$;
    \item if $\tilde{G}^{-1}(\beta)\leq {G}^{-1}(\beta)$, then either $\abs{{p}_k(X) - {G}^{-1}(\beta)} \leq \abs{\hat{p}_k(X) - p_k(X)}$ or $\hat{p}_k(X) \in \pare{\tilde{G}^{-1}(\beta), {G}^{-1}(\beta)}$;
\end{itemize}
Similar conditions are satisfied if $\hat{p}_k(X) < \tilde{G}^{-1}(\beta)$ and ${p}_k(X) \geq {G}^{-1}(\beta)$.
Using the above arguments we can upper-bound $(*)$ as
\begin{align*}
    (*)
    \leq
    &\sum_{k = 1}^K\Prob_X\pare{\abs{{p}_k(X) - {G}^{-1}(\beta)} \leq \abs{\hat{p}_k(X) - p_k(X)}}\\
    &+
    \ind{\tilde{G}^{-1}(\beta)\leq {G}^{-1}(\beta)}\sum_{k = 1}^K\Prob_X\pare{\tilde{G}^{-1}(\beta) < \hat{p}_k(X) < G^{-1}(\beta)}\\
    &+
    \ind{{G}^{-1}(\beta) < \tilde{G}^{-1}(\beta)}\sum_{k = 1}^K\Prob_X\pare{{G}^{-1}(\beta) < \hat{p}_k(X) < \tilde{G}^{-1}(\beta)}\\
    =
    &\sum_{k = 1}^K\Prob_X\pare{\abs{{p}_k(X) - {G}^{-1}(\beta)} \leq \abs{\hat{p}_k(X) - p_k(X)}}\\
    &+
    \abs{\tilde{G}\pare{\tilde{G}^{-1}(\beta)} - \tilde{G}\pare{G^{-1}(\beta)}}\enspace.
\end{align*}
Thanks to the continuity Assumption~\ref{ass:continuity_cdf_estimator} on the estimator and the continuity Assumption~\ref{ass:continuity_cdf} on the distribution we clearly have $\tilde{G}\pare{\tilde{G}^{-1}(\beta)} = \beta = {G}\pare{{G}^{-1}(\beta)}$.
Moreover, we can write
\begin{align*}
    \lvert\tilde{G}\pare{\tilde{G}^{-1}(\beta)}- &\tilde{G}\pare{G^{-1}(\beta)}\rvert =
    \abs{{G}\pare{{G}^{-1}(\beta)}- \tilde{G}\pare{G^{-1}(\beta)}}\\
    &\leq
    \sum_{k = 1}^K \Exp_{X \sim \Prob_X}\abs{\ind{\hat{p}_k(X) \geq G^{-1}(\beta)} - \ind{p_k(X) \geq G^{-1}(\beta)}}\\
    &\leq
    \sum_{k = 1}^K\Prob_X\pare{\abs{{p}_k(X) - {G}^{-1}(\beta)} \leq \abs{\hat{p}_k(X) - p_k(X)}}\enspace.
\end{align*}
Thus, our bound reads as
\begin{align*}
    (*) \leq 2\sum_{k = 1}^K\Prob_X\pare{\abs{{p}_k(X) - {G}^{-1}(\beta)} \leq \abs{\hat{p}_k(X) - p_k(X)}}\enspace.
\end{align*}
Finally, in order to upper-bound the term above one can use the peeling argument of~\cite{Audibert_Tsybakov07} applied with the exponential concentration inequality provided by Assumption~\ref{ass:exponential_concentration_estimator}.
This part of the proof we omit here and refer the reader to~\cite{Denis_Hebiri17} or to~\cite{Audibert_Tsybakov07} for a complete result.

Let us emphasize that the argument above is only possible due to the continuity Assumptions~\ref{ass:continuity_cdf},~\ref{ass:continuity_cdf_estimator} on the distribution and the estimator respectively.

%%%%%%%%%%%%%%%%%%%%%%%%%%%%%%%%%%%%%%%%%%%%%%%%%%%%%%%%%%%%%%%%%%%%%%%%%%%%%%%
%%%%%%%%%%%%%%%%%%%%%%%%%%%%%%%%%%%%%%%%%%%%%%%%%%%%%%%%%%%%%%%%%%%%%%%%%%%%%%%
\section{Proof of the lower bounds}
\label{sec:proof_of_the_lower_bounds}
%%%%%%%%%%%%%%%%%%%%%%%%%%%%%%%%%%%%%%%%%%%%%%%%%%%%%%%%%%%%%%%%%%%%%%%%%%%%%%%
%%%%%%%%%%%%%%%%%%%%%%%%%%%%%%%%%%%%%%%%%%%%%%%%%%%%%%%%%%%%%%%%%%%%%%%%%%%%%%%
This section is devoted to the proof of the lower bounds provided by Theorems~\ref{thm:lower_bound_supervised}-\ref{thm:lower_bound_semi_supervised}.
Before proceeding to the proofs let us briefly sketch the high-level strategy used in this work.
In order to prove the lower bounds of Theorems~\ref{thm:lower_bound_supervised}-\ref{thm:lower_bound_semi_supervised} we actually prove to separate lower bounds on the minimax risk.
Clearly, if some non-negative quantity is lower-bounded by two different values, therefore it is lower-bounded by the maximum between the two.
The two lower bounds that we prove are naturally connected with the proposed two-steps estimator, that is, the first lower bound is connected with the problem of non-parametric estimation of $p_k$ for all $k \in [K]$ and the second describes the estimation of the unknown threshold $G^{-1}(\beta)$.

In particular, the first lower bound is closely related to the one provided in~\citep{Audibert_Tsybakov07,Rigollet_Vert09}, though, crucially the continuity Assumption~\ref{ass:continuity_cdf} makes the proof more involved.
The second lower bound is based on two hypotheses testing and is derived by constructing two different marginal distributions of $X \in \bbR^d$ and a fixed regression vector $p(\cdot)$.
In this part we make use of Pinsker's inequality recalled in Lemma~\ref{lem:pinsker_inequality}.

In order to discriminate the supervised and the semi-supervised procedures we make use of Definition~\ref{def:supervised_semi_supervised}.
Notice that every supervised procedure thanks to Definition~\ref{def:supervised_semi_supervised} is not 'sensitive' to the expectation taken \wrt the unlabeled dataset $\data_N$, that is, randomness is only induced by the labeled dataset $\data_n$.
This strategy allows to eliminate the dependence of the lower bound on the size of the unlabeled dataset $\data_N$ for supervised procedures.
Indeed, let $\hGamma$ be any supervised estimator in the sense of Definition~\ref{def:supervised_semi_supervised}, then for any real valued function of confidence sets $Z$ we have
\begin{align*}
    \Exp_{(\data_n, \data_N)}[\Exp_{\Prob_X}Z(\hGamma(X; \data_n, \data_N))] = \Exp_{\data_n}[\Exp_{\Prob_X}Z(\hGamma(X; \data_n, \data'_N))]\enspace,
\end{align*}
with $\data'_N$ being an arbitrary set of $N$ points in $\bbR^d$.

%%%%%%%%%%%%%%%%%%%%%%%%%%%%%%%%%%%%%%%%%%%%%%%%%%%%%%%%%%%%%%%%%%%%%%%%%%%%%%%
%%%%%%%%%%%%%%%%%%%%%%%%%%%%%%%%%%%%%%%%%%%%%%%%%%%%%%%%%%%%%%%%%%%%%%%%%%%%%%%
\subsection{Part I: \texorpdfstring{$(N + n)^{-1/2}$}{Lg}}
\label{subsec:semi_supervised_part_of_the_rate}
%%%%%%%%%%%%%%%%%%%%%%%%%%%%%%%%%%%%%%%%%%%%%%%%%%%%%%%%%%%%%%%%%%%%%%%%%%%%%%%
%%%%%%%%%%%%%%%%%%%%%%%%%%%%%%%%%%%%%%%%%%%%%%%%%%%%%%%%%%%%%%%%%%%%%%%%%%%%%%%
Here we prove that the rate $(N + n)^{-1/2}$ is optimal for semi-supervised methods, as already mentioned the rate for the supervised methods can be obtained by formally setting $N = 0$.
The constant $C', C, c$ are always assumed to be independent of $N, n$ and can differ from line to line.
Let us fix $\beta \in \{1, \ldots, \floor{K / 2} - 1\}$ and $K \geq 5$.
For a positive constant $C < 1/2$ we define the following sequence
\begin{align*}
    \kappa_{N, n} = C(N + n)^{-\frac{1}{2}} < 0.5\enspace.
\end{align*}
To prove the lower bound we construct two distribution $\Prob_0$ and $\Prob_1$ on $\bbR^d$ sharing the same regression function $p(\cdot) = (p_1(\cdot), \ldots, p_K(\cdot))$ and with different marginals admitting densities $\mu_0, \mu_1$.
First, for a fixed parameter $0< \rho <1$ and fixed constants $0 < r_0 < r_1 < r_2 < r_3 < r_4$ to be specified we define the following sets
\begin{align*}
% r_i(3+\rho)/4, r_i(3\rho + 1)/4
    \class{X}_0 &= \enscond{x \in \bbR^d}{\norm{x} \leq r_0}\enspace,\\
    \class{X}_1 &= \enscond{x \in \bbR^d}{\norm{x - (\underbrace{r_1 + \rho, 0, \ldots, 0}_{\in \bbR^d})^\top} \leq \rho/2}\enspace,\\
    \class{X}_2 &= \enscond{x \in \bbR^d}{\norm{x - (\underbrace{r_2 + \rho, 0, \ldots, 0}_{\in \bbR^d})^\top} \leq \rho)}\enspace,\\
    \class{X}_3 &= \enscond{x \in \bbR^d}{\norm{x - (\underbrace{r_3 + \rho, 0, \ldots, 0}_{\in \bbR^d})^\top} \leq \rho/2}\enspace,\\
    \class{X}_4 &= \enscond{x \in \bbR^d}{r_4 \leq \norm{x} \leq 2 r_4}\enspace.
\end{align*}

\begin{figure}[!t]

  \centering
  \includegraphics[width=\linewidth]{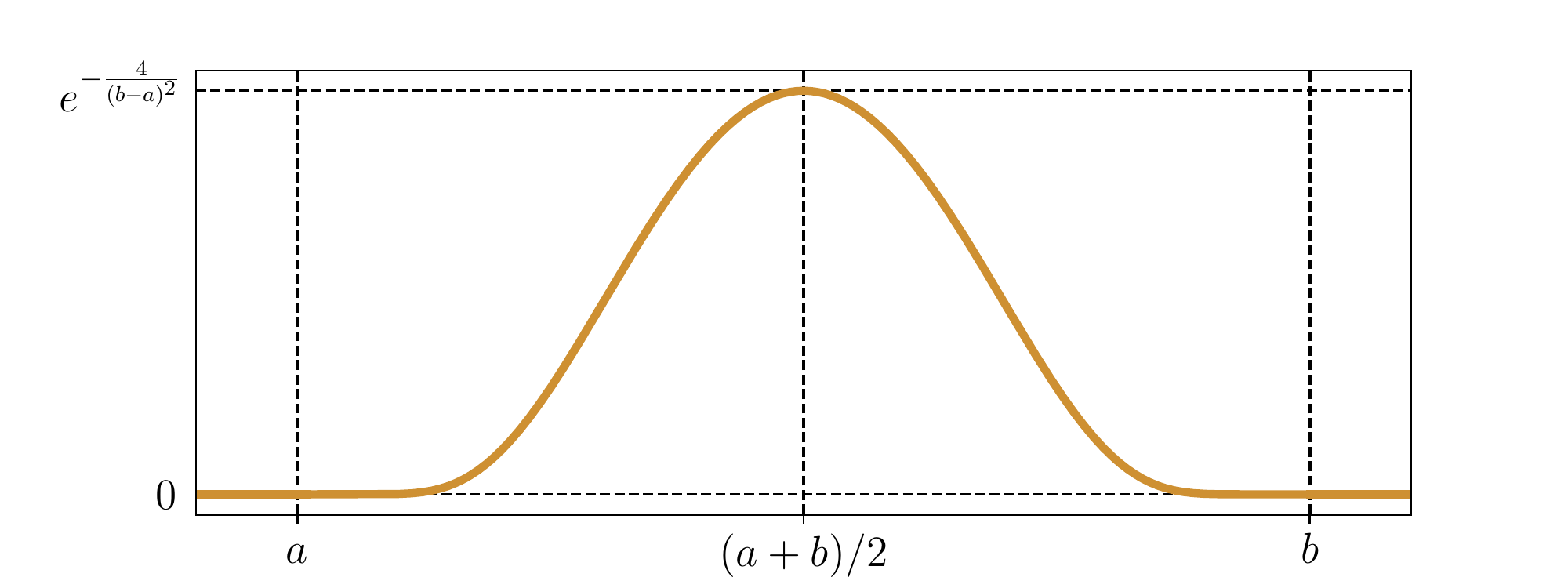}
  \captionof{figure}{Bump function: $x \mapsto \psi_{a, b}(x)$. Importantly, this function is supported on $(a, b)$ and is infinitely smooth.}
  \label{fig:psi_function}

\end{figure}

\begin{figure}[!t]
          \centering
          \includegraphics[width=\linewidth]{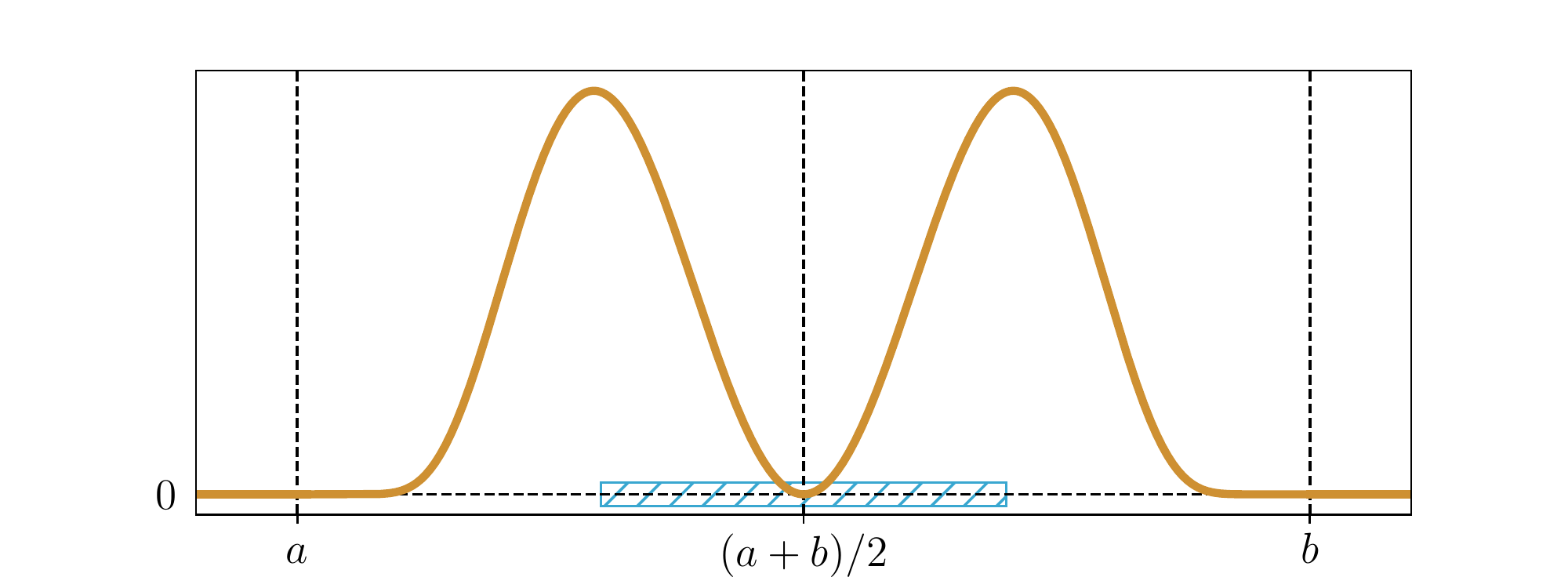}
          \captionof{figure}{Dumped bump function: $x \mapsto \pare{x - \frac{a + b}{2}}^{2\ceil{\frac{\gamma}{2}}}\psi_{a, b}(x)$. Importantly, this function behaves as polynomial of even degree $2\ceil{\frac{\gamma}{2}}$ in the affinity of $\tfrac{a + b}{2}$, while being infinitely smooth and supported on $(a, b)$.
          It means that if we select a measure which is supported in the affinity of $\tfrac{a + b}{2}$ (light-blue hatched region) the function on the plot is essentially polynomial \wrt such a measure.}
          \label{fig:dumped_function}
    \end{figure}

Let us denote by $o_i = (r_i + \rho, 0, \ldots, 0)^\top$ for $i=1,2,3$ the centers of $\class{X}_1$, $\class{X}_2$ and $\class{X}_3$.
Using these sets we define the regression vector as

\begin{align*}
    p_1(x) = \ldots = p_{2\beta}(x) &= \begin{cases}
                                    \frac{1}{2\beta} - \frac{\varphi_0(x)}{2\beta}, &\quad x \in \class{X}_0\\
                                    \frac{3K + 2\beta}{8K\beta} - \frac{\varphi_1(x)}{2\beta}, &\quad x \in \class{X}_1\\
                                    \frac{2K + 4\beta}{8K\beta} - \frac{\varphi_2(x)}{2\beta}, &\quad x \in \class{X}_2\\
                                    \frac{K + 6\beta}{8K\beta} - \frac{\varphi_3(x)}{2\beta}, &\quad x \in \class{X}_3\\
                                    \frac{1}{K} - \frac{\varphi_4(x)}{2\beta}, &\quad x \in \class{X}_4
                               \end{cases}\\
    p_{2\beta + 1}(x) = \ldots = p_{K}(x)
    &= \begin{cases}
                                    \frac{\varphi_0(x)}{K - 2\beta}, &\quad x \in \class{X}_0\\
                                    \frac{1}{4K} + \frac{\varphi_1(x)}{K - 2\beta}, &\quad x \in \class{X}_1\\
                                    \frac{2K - 4\beta}{4K(K - 2\beta)} + \frac{\varphi_2(x)}{K - 2\beta}, &\quad x \in \class{X}_2\\
                                    \frac{3K - 6\beta}{4K(K - 2\beta)} + \frac{\varphi_3(x)}{K - 2\beta}, &\quad x \in \class{X}_3\\
                                    \frac{1}{K} + \frac{\varphi_4(x)}{K - 2\beta}, &\quad x \in \class{X}_4
                               \end{cases}\enspace,
\end{align*}
In order to define the functions $\varphi_i$ for $i = 0, \ldots, 4$ we first define a one dimensional function of two real-valued parameters $a < b$
\begin{align*}
     \psi_{a, b}(x) = \begin{cases}
                        \exp\pare{-\frac{1}{(b - x)(x - a)}}, &\quad x \in (a, b)\\
                        0, &\quad \text{otherwise}
                      \end{cases}\enspace.
 \end{align*}
Figure~\ref{fig:psi_function} illustrates the behavior of $\psi_{a, b}$ function in one dimension.
Note that for every $a, b \in \bbR$ the function above is infinitely smooth.
Using the definition of $\psi_{a, b}$ we define the functions $\varphi_i$ for $i = 0, \ldots, 4$ as
\begin{align*}
    \varphi_0(x) &= \frac{C'}{2}\pare{\frac{K - 2\beta}{8K\beta} \bigwedge \frac{1}{4K}}\psi_{-1, r_0}(\norm{x})\enspace,\\
    \varphi_i(x) &= \frac{C'}{2}\rho^{\gamma}\pare{\frac{K - 2\beta}{8K\beta} \bigwedge \frac{1}{4K}}\pare{\frac{\norm{x - o_i}}{\rho}}^{2 \ceil{\frac{\gamma}{2}}}\psi_{-1, 1}\pare{\frac{\norm{x - o_i}}{\rho}},\quad i = 1, 3\enspace,\\
    \varphi_2(x) &= \frac{C'}{2}\rho^{\gamma}\pare{\frac{K - 2\beta}{8K\beta} \bigwedge \frac{1}{4K}}\psi_{-1, 1}\pare{\frac{\norm{x - o_2}}{\rho}}\enspace,\\
    \varphi_4(x) &= \frac{C'}{2}\pare{\frac{K - 2\beta}{8K\beta} \bigwedge \frac{1}{4K}}\psi_{r_4, 2 r_4}(\norm{x})\enspace,
\end{align*}
and the constant $C' \leq 1$ is chosen small enough so that each function $\varphi_i$ for $i = 0, \ldots, 4$ is $(\gamma, L)$-H\"older.
Let us point out that such value $C'$ exists and is independent of $n, N$, indeed, the mapping
\begin{align*}
    x \mapsto C'\norm{x}^{2 \ceil{\frac{\gamma}{2}}} \psi_{-1, 1}(\norm{x})\enspace,
\end{align*}
is infinitely smooth, thus it is $(\gamma, L)$-H\"older for a properly chosen $C'$.
Figure~\ref{fig:dumped_function} demonstrates the behavior of the considered construction in one dimension.
Note that $\varphi_i(x)$ for $i = 1,3$ are obtained from the previous mapping by re-scaling which preserves the H\"older constant $L$.
Same reasoning applies to $\varphi_i$ for $i = 0, 2, 4$.

Since $\beta < K / 2$ one can check that the following relations hold true
\begin{align*}
    0 &< \frac{1}{4K} < \frac{2K - 4\beta}{4K(K - 2\beta)} < \frac{3K - 4\beta}{4K(K - 2\beta)} < \frac{3K - 6\beta}{4K(K - 2\beta)} < \frac{1}{K} \\
    &< \frac{K + 6\beta}{8K\beta} < \frac{2K + 4\beta}{8K\beta} < \frac{3K + 2\beta}{8K\beta} < \frac{1}{2\beta}\enspace,
\end{align*}
which will help us to ensure that the thresholds under $\Prob_0, \Prob_1$ are $\frac{3K + 2\beta}{8K\beta}$ and $\frac{K + 6\beta}{8K\beta}$ respectively.
Now, we define two marginal distributions $\mu_0, \mu_1$ by their densities as
\begin{align*}
    \mu_0(x) = \begin{cases}
                    \frac{1/2}{\Leb(\class{X}_0)}, &\quad x \in \class{X}_0\\
                    \frac{\kappa_{N, n}}{\Leb(\class{X}_1)}, &\quad x \in \class{X}_1\\
                    \frac{\kappa_{N, n}}{\Leb(\class{X}_2)}, &\quad x \in \class{X}_2\\
                    \frac{\kappa_{N, n}}{\Leb(\class{X}_3)}, &\quad x \in \class{X}_3\\
                    \frac{1/2  - 3\kappa_{N, n}}{\Leb(\class{X}_4)}, &\quad x \in \class{X}_4
               \end{cases},\quad
    \mu_1(x) = \begin{cases}
                    \frac{1/2 - 3\kappa_{N, n}}{\Leb(\class{X}_0)}, &\quad x \in \class{X}_0\\
                    \frac{\kappa_{N, n}}{\Leb(\class{X}_1)}, &\quad x \in \class{X}_1\\
                    \frac{\kappa_{N, n}}{\Leb(\class{X}_2)}, &\quad x \in \class{X}_2\\
                    \frac{\kappa_{N, n}}{\Leb(\class{X}_3)}, &\quad x \in \class{X}_3\\
                    \frac{1/2}{\Leb(\class{X}_4)}, &\quad x \in \class{X}_4
               \end{cases}\enspace,
\end{align*}
and both $\mu_0, \mu_1$ are equal to zero in unspecified regions.
Clearly, the strong density assumption is satisfied on $\class{X}_0$ and $\class{X}_4$ since the density is lower and upper-bounded by a constant independent of both $N, n$.
The parameter $\rho$ is chosen such that the strong density assumption on $\class{X}_i$ for $i = 1, 2, 3$ is satisfied.
Notice that
\begin{align*}
    \Leb(\class{X}_i) = c \rho^d\enspace,
\end{align*}
for some constant $c > 0$ independent of $N, n$, thus we set $\rho = C(N + n)^{-1/2d}$.
For these hypotheses one can easily check that the thresholds $G^{-1}_0(\beta), G^{-1}_1(\beta)$ and the optimal $\beta$-sets $\Gamma^*_0, \Gamma^*_1$ are given as
\begin{align*}
     G^{-1}_0(\beta) = &\frac{3K + 2\beta}{8K\beta}, \quad G^{-1}_1(\beta) = \frac{K + 6\beta}{8K\beta}\enspace,\\
     \Gamma^*_0(x) &= \begin{cases}
                        \{1, \ldots, 2\beta\}, &\quad x \in \class{X}_0\\
                        \emptyset, &\quad \text{otherwise}
                     \end{cases}, \\
     \Gamma^*_1(x) &= \begin{cases}
                        \{1, \ldots, 2\beta\}, &\quad x \in \class{X}_0 \bigcup \class{X}_1 \bigcup \class{X}_2, \\
                        \emptyset, &\quad \text{otherwise}
                     \end{cases}\enspace.
 \end{align*}

{\bf The margin assumption:} we are in position to check the margin Assumption~\ref{ass:margin_assumption}.
Let $t_0 = \frac{1}{2}\pare{\frac{K - 2\beta}{8K\beta} \bigwedge \frac{1}{4K}}$, thus for every $k \in \{2\beta + 1, \ldots, K\}$ and every $t \leq t_0$ we have
\begin{align*}
    \Prob_0\pare{\abs{p_k(X) - G^{-1}_0(\beta)} \leq t} = 0,\quad
    \Prob_1\pare{\abs{p_k(X) - G^{-1}_1(\beta)} \leq t} = 0\enspace,
\end{align*}
moreover for every $k \in \{1, \ldots, 2\beta\}$ and every $t \leq t_0$ we can write
\begin{small}
\begin{align*}
    \Prob_0\Bigg(&\abs{p_k(X) - G^{-1}_0(\beta)} \leq t\Bigg)\\
    =
    &\Prob_0\pare{\frac{C'}{2}\rho^{\gamma}\pare{\frac{K - 2\beta}{8K\beta} \bigwedge \frac{1}{4K}}\pare{\frac{\norm{x - o_1}}{\rho}}^{2 \ceil{\frac{\gamma}{2}} }\psi_{-1, 1}\pare{\frac{\norm{x - o_1}}{\rho}} \leq 2\beta t, X\in \class{X}_1}\enspace,\\
    \Prob_1\Bigg(&\abs{p_k(X) - G^{-1}_1(\beta)} \leq t\Bigg)\\
    =
    &\Prob_1\pare{\frac{C'}{2}\rho^{\gamma}\pare{\frac{K - 2\beta}{8K\beta} \bigwedge \frac{1}{4K}}\pare{\frac{\norm{x - o_3}}{\rho}}^{2 \ceil{\frac{\gamma}{2}} }\psi_{-1, 1}\pare{\frac{\norm{x - o_3}}{\rho}}\leq 2\beta t,  X\in \class{X}_3}\enspace.
\end{align*}
\end{small}
Hence, for the $0$ hypothesis there exists $c$ independent of $N, n$ such that
\begin{small}
\begin{align*}
    \Prob_0\pare{\abs{p_k(X) - G^{-1}_0(\beta)} \leq t}
    &\leq
    \Prob_0\pare{ \pare{\frac{\norm{x - o_1}}{\rho}}^{2 \ceil{\frac{\gamma}{2}} }\psi_{-1, 1}\pare{\frac{\norm{x - o_1}}{\rho}} \leq c^2 \rho^{-\gamma}t, X\in \class{X}_1}
\end{align*}
\end{small}
Therefore we can write using the strong density assumption
\begin{small}
\begin{align*}
    \Prob_0\pare{\abs{p_k(X) - G^{-1}_0(\beta)} \leq t}
    &\leq
    \int_{\norm{x - o_1} \leq \rho/2}\ind{\pare{\frac{\norm{x - o_1}}{\rho}}^{2 \ceil{\frac{\gamma}{2}} }\psi_{-1, 1}\pare{\frac{\norm{x - o_1}}{\rho}} \leq c^2 \rho^{-\gamma}t}d\mu_0(x)\\
    &\leq
    C\int_{\norm{x - o_1} \leq \rho/2}\ind{\pare{\frac{\norm{x - o_1}}{\rho}}^{2 \ceil{\frac{\gamma}{2}} }\psi_{-1, 1}\pare{\frac{\norm{x - o_1}}{\rho}} \leq c^2 \rho^{-\gamma}t}dx\\
    &=
    C\int_{\norm{x} \leq \rho/2}\ind{\pare{\frac{\norm{x}}{\rho}}^{2 \ceil{\frac{\gamma}{2}} }\psi_{-1, 1}\pare{\frac{\norm{x}}{\rho}} \leq c^2 \rho^{-\gamma}t}dx\\
    &=
    C\rho^d\int_{\norm{x} \leq 1/2}\ind{{{\norm{x}}}^{2 \ceil{\frac{\gamma}{2}} }\psi_{-1, 1}\pare{{\norm{x}}} \leq c^2 \rho^{-\gamma}t}dx\enspace,
\end{align*}
\end{small}
Finally notice that for every $x \in \bbR^d$ such that $\norm{x} \leq 1/2$ we have for some $C > 0$
\begin{align*}
    \psi_{-1, 1}\pare{{\norm{x}}} \geq \psi_{-1, 1}\pare{1/2} \geq C\enspace,
\end{align*}
which implies that for some positive $C, C'$ independent of $N, n$ we can write
\begin{small}
\begin{align*}
    \Prob_0\pare{\abs{p_k(X) - G^{-1}_0(\beta)} \leq t}
    &\leq
    C\rho^d\int_{\norm{x} \leq 1/2}\ind{{{\norm{x}}}^{2\ceil{\frac{\gamma}{2}}} \leq C' \rho^{-\gamma}t}dx\\
    &=
    C\rho^d\int_{\norm{x} \leq 1/2}\ind{{{\norm{x}}} \leq C'^{1/(2\ceil{\frac{\gamma}{2}})} \rho^{-\gamma/(2 \ceil{\frac{\gamma}{2}})}t^{1/ (2 \ceil{\frac{\gamma}{2}})}}dx\\
    &\leq
    C\rho^{d(1 - \gamma/( 2 \ceil{\frac{\gamma}{2}})}t^{d/(2 \ceil{\frac{\gamma}{2}} )}\enspace.
\end{align*}
\end{small}
This implies that for as long as $\alpha \leq d/(2 \ceil{\frac{\gamma}{2}})$ (and since we have $\gamma \leq 2 \ceil{\frac{\gamma}{2}}$) the margin assumption is satisfied.
Moreover, these conditions imply that $\alpha \gamma \leq d$, which we will also require while proving the supervised part of the rate.
Same reasoning can be carried out for the case of the first hypothesis $\Prob_1$ on the set $\class{X}_3$.

Finally, the parameters $r_0, r_1, r_2, r_3$ are chosen as constants independent of $n, N$ such that there exists a smooth connection between the parts of the regression functions $p_k(\cdot)$ which are defined on $\class{X}_0, \class{X}_1, \class{X}_2, \class{X}_3, \class{X}_4$.
Notice that such a choice is possible since by the construction of functions $\varphi_i$ for $i = 0, 1, 2, 3, 4$ they are zeroed-out on the boundaries of $\class{X}_0, \class{X}_1, \class{X}_2, \class{X}_3, \class{X}_4$.
Thus in the region $\bbR^d \setminus\bigcup_{i = 0}^4 \class{X}_i$ it is sufficient to construct a function which connects four different constants smoothly.
We avoid this over complication on this part and hope that the guidelines provided above are sufficient for the understanding.

Notice that the constructed distributions are satisfying Assumption~\ref{ass:continuity_cdf} since the measures are only defined on $\class{X}_0, \class{X}_1, \class{X}_2, \class{X}_3, \class{X}_4$ and the regression functions on these sets are not concentrated around any constant.
% now if $2\beta t \geq \frac{C'}{2}\tau^{-\gamma}_{n, N} > \max_{x \in \bbR^d}\varphi_1(x)$, then
% \begin{align*}
%     \Prob_0\pare{ \varphi_1(X) \leq 2\beta t, X\in \class{X}_1} \leq \1_{2\beta t \geq \frac{C'}{2}\tau^{-\gamma}_{n, N}}\Prob_0\pare{X\in \class{X}_1} = \kappa_{N, n}\1_{2\beta t \geq \frac{C'}{2}\tau^{-\gamma}_{n, N}} \leq \pare{\frac{4\beta}{C'}}^{\alpha}\kappa_{N, n}\tau^{\alpha\gamma}_{n, N}t^{\alpha}
% \end{align*}
% if $2\beta t < \frac{C'}{2}\tau^{-\gamma}_{n, N}$, then
% \begin{align*}
%     \Prob_0\pare{ \varphi_1(X) \leq 2\beta t, X\in \class{X}_1} \leq \1_{2\beta t < \frac{C'}{2}\tau^{-\gamma}_{n, N}}\Prob_0\pare{X\in \class{X}_1} = \1_{2\beta t < \frac{C'}{2}\tau^{-\gamma}_{n, N}}\kappa_{N, n}
% \end{align*}

Before proceeding to the final stage of the proof let us mention that in what follows we use the de Finetti~\citep{Finetti72,Finetti74} notation which is common in probability.
That is, given a probability measure $\Prob$ on some measurable space $(\Omega_0, \class{A}_0)$ and a measurable function $X: (\Omega_0, \class{A}_0) \to (\bbR, \text{Borel}(\bbR))$ we write
\begin{align*}
    \Prob[X] \eqdef \Exp[X]\enspace.
\end{align*}

{\bf Bound on the KL-divergence:} we start by computing the KL-divergence between $\mu_0$ and $\mu_0$
\begin{align*}
    \text{KL}&(\mu_0, \mu_1)
    \eqdef
    \int_{\bbR^d} \mu_0(x)\log\pare{\frac{\mu_0}{\mu_1}} dx
    =
    \sum_{i = 0}^4 \int_{x \in \class{X}_i}\mu_0(x)\log\pare{\frac{\mu_0(x)}{\mu_1(x)}} dx\\
    &=
    \frac{1}{\Leb(\class{X}_0)}\int_{x \in \class{X}_0}\frac{1}{2}\log\pare{\frac{1/2}{1/2 - 3\kappa_{N, n}}} dx\\
    &\phantom{=}
    +
    \frac{1}{\Leb(\class{X}_4)}\int_{x \in \class{X}_4}\pare{\frac{1}{2} - 3\kappa_{N, n}}\log\pare{\frac{1/2 - 3\kappa_{N, n}}{1/2}} dx\\
    &=
    \frac{1}{2}\log\pare{\frac{1/2}{1/2 - 3\kappa_{N, n}}}\\
    &\phantom{=}+
    \pare{\frac{1}{2} - 3\kappa_{N, n}}\log\pare{\frac{1/2 - 3\kappa_{N, n}}{1/2}}\\
    &=
    -3\kappa_{N, n}\log\pare{{1 - 6\kappa_{N, n}}}
    \leq 36\kappa_{N, n}^2\enspace.
\end{align*}

{\bf Lower bound for the Hamming risk:} first of all let us introduce the following notation for $i = 0, 1$
\begin{align*}
    \Ham(\hGamma, \Gamma^*_i) \eqdef \mu_i\abs{\hGamma(X) \triangle \Gamma^*_i(X)} \enspace.
\end{align*}
Recall that we are interested in the following quantity
\begin{align*}
    \inf_{\hGamma}\sup_{\Prob \in \class{P}} \Exp_{(\data_n, \data_N)} \Exp_{X \sim \Prob_X} \abs{\hGamma(X) \triangle \Gamma^*_\beta(X)}\enspace,
\end{align*}
since the hypotheses $\Prob_0, \Prob_1 \in \class{P}$ we can write
\begin{align*}
    2\sup_{\Prob \in \class{P}} \Exp_{(\data_n, \data_N)} \Exp_{X \sim \Prob_X} \abs{\hGamma(X) \triangle \Gamma^*_\beta(X)} \geq (*)\enspace,
\end{align*}
where $(*)$ is defined as
\begin{align*}
    (*) = \mu_0^{\otimes (n + N)}\otimes \Prob_{Y | X}^{\otimes n} \Ham(\hGamma, \Gamma^*_0) + \mu_1^{\otimes (n + N)}\otimes \Prob_{Y | X}^{\otimes n} \Ham(\hGamma, \Gamma^*_1)\enspace,
\end{align*}
 thus, for the Hamming risk we can write
 \begin{align*}
    (*)
    &\geq
   \mu_0^{\otimes (n + N)}\otimes \Prob_{Y | X}^{\otimes n}\pare{\frac{d\mu_1^{\otimes (n + N)}\otimes \Prob_{Y | X}^{\otimes n}}{d\mu_0^{\otimes (n + N)}\otimes \Prob_{Y | X}^{\otimes n}} \bigwedge 1}\pare{\Ham(\hGamma, \Gamma^*_0) + \Ham(\hGamma, \Gamma^*_1)}\enspace.
 \end{align*}
 Now we focus our attention to the sum of two Hamming differences which appearing on the right hand side of the above inequality
 \begin{align*}
    \Ham(\hGamma, \Gamma^*_0) + &\Ham(\hGamma, \Gamma^*_1)
     = \mu_0\sum_{k = 1}^K\1_{k \in \hGamma(X) \triangle \Gamma^*_0(X)}
      + \mu_1\sum_{k = 1}^K\1_{k \in \hGamma(X) \triangle \Gamma^*_1(X)}\\
     &\geq
     \mu_0\pare{\frac{d\mu_1}{d\mu_0} \bigwedge 1}\sum_{k = 1}^K\1_{k \in \hGamma(X) \triangle \Gamma^*_0(X)} \\
     &\phantom{=}+ \mu_0\pare{\frac{d\mu_1}{d\mu_0} \bigwedge 1}\sum_{k = 1}^K\1_{k \in \hGamma(X) \triangle \Gamma^*_1(X)}\\
     &\geq
     {\mu_0\pare{\frac{d\mu_1}{d\mu_0} \bigwedge 1}\sum_{k = 1}^K\1_{k \in \Gamma_1^*(X) \triangle \Gamma^*_0(X)}} \quad\text{(Triangle inequality)}\\
     &=
     2\beta\mu_0\pare{\frac{d\mu_1}{d\mu_0} \bigwedge 1}\pare{\1_{\class{X}_1} + \1_{\class{X}_2}}\\
     &=
     2\beta\int_{\bbR^d}\pare{\frac{\mu_1(x)}{\mu_0(x)} \bigwedge 1}\pare{\1_{\class{X}_1} + \1_{\class{X}_2}} d\mu_0(x)\\
     &=
     2\beta\int_{\class{X}_1}\pare{\frac{\mu_1(x)}{\mu_0(x)} \bigwedge 1} d\mu_0(x)
     +
     2\beta\int_{\class{X}_2}\pare{\frac{\mu_1(x)}{\mu_0(x)} \bigwedge 1} d\mu_0(x)\\
     &=
     2\beta\Prob_0(\class{X}_1 \cup \class{X}_2) \geq 2\beta\kappa_{n, N}\enspace.
 \end{align*}
 Substituting this lower bound into the initial inequality we arrive at
 \begin{align*}
     (*)
     &\geq
     2\beta\kappa_{n, N}\mu_0^{\otimes (n + N)}\otimes \Prob_{Y | X}^{\otimes n}\pare{\frac{d\mu_1^{\otimes (n + N)}\otimes \Prob_{Y | X}^{\otimes n}}{d\mu_0^{\otimes (n + N)}\otimes \Prob_{Y | X}^{\otimes n}} \bigwedge 1}\\
     &=
     2\beta\kappa_{n, N}\pare{1 - \text{TV}\pare{\mu_0^{\otimes (n + N)}\otimes \Prob_{Y | X}^{\otimes n}, \mu_1^{\otimes (n + N)}\otimes \Prob_{Y | X}^{\otimes n}}}\\
     &=
     2\beta\kappa_{n, N}\pare{1 - \text{TV}\pare{\mu_0^{\otimes (n + N)}, \mu_1^{\otimes (n + N)}}}\\
     &\geq
     2\beta\kappa_{n, N}\pare{1 - \sqrt{\frac{1}{2}\text{KL}\pare{\mu_0^{\otimes (n + N)}, \mu_1^{\otimes (n + N)}}}}\quad\text{(Pinsker's inequality)}\\
     &\geq
     2\beta\kappa_{n, N}\pare{1 - \kappa_{n , N}\sqrt{n + N}}\enspace,
 \end{align*}
which implies the desired lower bound on the Hamming risk.

{\bf Lower bound for the $\beta$ excess risk:} this part is analogues to the case of the Hamming distance.
Let us recall that for every $\hGamma$ we have the following expression for $i = 0, 1$
\begin{align*}
    \Disc(\hGamma, \Gamma^*_i) \eqdef \risk_{\beta}(\hGamma) - \risk_\beta(\Gamma^*_i) = \mu_i\sum_{k = 1}^K\abs{p_k(X) - G^{-1}(\beta)}\ind{k \in \hGamma(X) \triangle \Gamma^*_i(X)}\enspace.
\end{align*}
Again, recall that we are interested in
\begin{align*}
    \inf_{\hGamma}\sup_{\Prob \in \class{P}}\Exp_{(\data_n, \data_N)}[\risk_{\beta}(\hGamma)] - \risk(\Gamma^*_\beta)
\end{align*}
similarly to the previous case, since the hypotheses $\Prob_0, \Prob_1 \in \class{P}$ we can write
\begin{align*}
    2\sup_{\Prob \in \class{P}}\Exp_{(\data_n, \data_N)}[\risk_{\beta}(\hGamma)] - \risk(\Gamma^*_\beta) \geq (**)\enspace,
\end{align*}
where $(**)$ is defined as
 \begin{align*}
     (**) = \mu_0^{\otimes (n + N)}\otimes \Prob_{Y | X}^{\otimes n} \Disc(\hGamma, \Gamma^*_0) + \mu_1^{\otimes (n + N)}\otimes \Prob_{Y | X}^{\otimes n} \Disc(\hGamma, \Gamma^*_1) \enspace,
 \end{align*}
 we can write
 \begin{small}
 \begin{align*}
     &(**)
    \geq
    \mu_0^{\otimes (n + N)}\otimes \Prob_{Y | X}^{\otimes n}\pare{\frac{d\mu_1^{\otimes (n + N)}\otimes \Prob_{Y | X}^{\otimes n}}{d\mu_0^{\otimes (n + N)}\otimes \Prob_{Y | X}^{\otimes n}} \bigwedge 1}\pare{\Disc(\hGamma, \Gamma^*_0) + \Disc(\hGamma, \Gamma^*_1)}
 \end{align*}
 and we continue in a similar fashion
 \begin{align*}
     \Disc(\hGamma, \Gamma^*_0) + \Disc(\hGamma, \Gamma^*_1)
     &=
     \mu_0\sum_{k = 1}^K\abs{p_k(X) - \frac{3K + 2\beta}{8K\beta}}\1_{k\in \hGamma(X) \triangle \Gamma^*_0(X)}\\
     &\phantom{=}+
     \mu_1\sum_{k = 1}^K\abs{p_k(X) - \frac{K + 6\beta}{8K\beta}}\1_{k\in \hGamma(X) \triangle \Gamma^*_1(X)}\\
     &\geq
     \mu_0\sum_{k = 1}^{2\beta}\abs{p_k(X) - \frac{3K + 2\beta}{8K\beta}}\1_{k\in \hGamma(X) \triangle \Gamma^*_0(X)}\1_{X \in \class{X}_2}\\
     &\phantom{=}+
     \mu_1\sum_{k = 1}^{2\beta}\abs{p_k(X) - \frac{K + 6\beta}{8K\beta}}\1_{k\in \hGamma(X) \triangle \Gamma^*_1(X)}\1_{X \in \class{X}_2}\enspace,
 \end{align*}
 since $\mu_0(x) = \mu_1(x)$ for all $x \in \class{X}_2$ we obtain
 \begin{small}
 \begin{align*}
     \Disc(\hGamma, &\Gamma^*_0) + \Disc(\hGamma, \Gamma^*_1)\\
     &\geq
     \mu_0\Bigg(\sum_{k = 1}^{2\beta}\abs{p_k(X) - \frac{3K + 2\beta}{8K\beta}}\1_{k\in \hGamma(X) \triangle \Gamma^*_0(X)}\1_{X \in \class{X}_2}\\
     &\phantom{=\mu_0\Bigg(}+
     \sum_{k = 1}^{2\beta}\abs{p_k(X) - \frac{K + 6\beta}{8K\beta}}\1_{k\in \hGamma(X) \triangle \Gamma^*_1(X)}\1_{X \in \class{X}_2}\Bigg)\\
     % &\geq
     % \mu_0\pare{\sum_{k = 1}^{2\beta}\pare{\abs{p_k(X) - \frac{3K + 2\beta}{8K\beta}} \bigwedge \abs{p_k(X) - \frac{K + 6\beta}{8K\beta}}}\pare{\1_{k\in \hGamma(X) \triangle \Gamma^*_0(X)} + \1_{k\in \hGamma(X) \triangle \Gamma^*_1(X)}}{\1_{X \in \class{X}_2}}}\\
     &\geq
     \mu_0\pare{\sum_{k = 1}^{2\beta}\pare{\abs{p_k(X) - \frac{3K + 2\beta}{8K\beta}} \bigwedge \abs{p_k(X) - \frac{K + 6\beta}{8K\beta}}}\1_{k\in \Gamma^*_1(X) \triangle \Gamma^*_0(X)}{\1_{X \in \class{X}_2}}}\\
     &=
     \mu_0\pare{\sum_{k = 1}^{2\beta}\pare{\abs{p_k(X) - \frac{3K + 2\beta}{8K\beta}} \bigwedge \abs{p_k(X) - \frac{K + 6\beta}{8K\beta}}}{\1_{X \in \class{X}_2}}}\\
     &=
     \mu_0\pare{2\beta\pare{\abs{ \frac{2K + 4\beta}{8K\beta} - \frac{\varphi_2(X)}{2\beta} - \frac{3K + 2\beta}{8K\beta}} \bigwedge \abs{ \frac{2K + 4\beta}{8K\beta} - \frac{\varphi_2(X)}{2\beta} - \frac{K + 6\beta}{8K\beta}}}{\1_{X \in \class{X}_2}}}\\
     &=
     \mu_0\pare{2\beta\pare{\abs{ \frac{K - 2\beta}{8K\beta} + \frac{\varphi_2(X)}{2\beta}} \bigwedge \abs{ \frac{K - 2\beta}{8K\beta} - \frac{\varphi_2(X)}{2\beta}}}{\1_{X \in \class{X}_2}}}\\
     &=
     \mu_0\pare{2\beta\abs{ \frac{K - 2\beta}{8K\beta} - \frac{\varphi_2(X)}{2\beta}}{\1_{X \in \class{X}_2}}}\enspace,
 \end{align*}
 \end{small}
 then, since $\frac{\varphi_2(x)}{\beta} \leq \frac{K - 2\beta}{8K\beta}$ for all $x \in \class{X}_2$, we have
 \begin{align*}
     \Disc(\hGamma, \Gamma^*_0) + \Disc(\hGamma, \Gamma^*_1) \geq \frac{2\beta (K - 2\beta)}{16K\beta}\mu_0(\class{X}_2) = \frac{K - 2\beta}{8K}\kappa_{n, N}\enspace.
 \end{align*}
 Thus,
 \begin{align*}
     (**)
     &\geq
     \frac{K - 2\beta}{8K}\kappa_{n, N}\pare{1 - \text{TV}\pare{\mu_0^{\otimes (n + N)}\otimes \Prob_{Y | X}^{\otimes n}, \mu_1^{\otimes (n + N)}\otimes \Prob_{Y | X}^{\otimes n}}}\\
     &\geq
     \frac{K - 2\beta}{8K}\kappa_{n, N}\pare{1 - \sqrt{\frac{1}{2}\text{KL}\pare{\mu_0^{\otimes (n + N)}\otimes \Prob_{Y | X}^{\otimes n}, \mu_1^{\otimes (n + N)}\otimes \Prob_{Y | X}^{\otimes n}}}}\\
     &\geq
     \frac{K - 2\beta}{8K}\kappa_{n, N}\pare{1 - \kappa_{n, N}\sqrt{n + N}}
 \end{align*}
Which concludes the first part of the lower bounds.
 \end{small}

%%%%%%%%%%%%%%%%%%%%%%%%%%%%%%%%%%%%%%%%%%%%%%%%%%%%%%%%%%%%%%%%%%%%%%%%%%%%%%%
%%%%%%%%%%%%%%%%%%%%%%%%%%%%%%%%%%%%%%%%%%%%%%%%%%%%%%%%%%%%%%%%%%%%%%%%%%%%%%%
\subsection{Part II: \texorpdfstring{$n^{-\alpha \gamma / (2\gamma + d)}$}{Lg}}
\label{subsec:supervised_part_of_the_rate}
%%%%%%%%%%%%%%%%%%%%%%%%%%%%%%%%%%%%%%%%%%%%%%%%%%%%%%%%%%%%%%%%%%%%%%%%%%%%%%%
%%%%%%%%%%%%%%%%%%%%%%%%%%%%%%%%%%%%%%%%%%%%%%%%%%%%%%%%%%%%%%%%%%%%%%%%%%%%%%%
In this section we prove that in case of the Hamming risk $\excess^{\Ham}$ the rate $n^{-\alpha \gamma/(2\gamma + d)}$ is minimax optimal.
Notice, that thanks to Proposition~\ref{prop:distanceComparison} a lower bound of order $n^{-\alpha \gamma/(2\gamma + d)}$ on the Hamming risk $\excess^{\Ham}$ immediately implies a lower bound of order $n^{-(\alpha + 1) \gamma/(2\gamma + d)}$ on both $\excess^{\Exc}$ and $\excess^{\Disc}$.

The proof is based on the reduction of the Hamming risk to a multiple hypotheses testing problem and an application of Fano's inequality provided by~\cite{Birge05} recalled in Lemma~\ref{lem:birges_for_lower_bound}.

Assume that $K \geq 5$ and fix some $\beta \in \{2, \ldots, (K - 2) \wedge \floor{K / 2}\}$, define the regular grid on $[0, 1]^d$ as
    \begin{align*}
        G_q \coloneqq \left\{\left(\frac{2k_1 + 1}{2q}, \ldots, \frac{2k_d + 1}{2q}\right)^\top \,:\, k_i \in \{0, \ldots, q - 1\}, i = 1, \ldots, d\right\}\enspace,
    \end{align*}
    and denote by $n_q(x) \in G_q$ as the closest point to of the grid $G_q$ to the point $x \in \bbR^d$.
    Such a grid defines a partition of the unit cube $[0, 1]^d \subset \bbR^d$ denoted by $\class{X}'_1, \ldots, \class{X}'_{q^d}$.
    Besides, denote by $\class{X}'_{-j} \coloneqq \{x \in \bbR^d\, : \, -x \in \class{X}_j'\}$ for all $j = 1, \ldots, q^d$.
    For a fixed integer $m \leq q^d$ and for any $j \in \{1, \ldots, m\}$ define $\class{X}_i \coloneqq \class{X}_i'$, $\class{X}_{-i} \coloneqq \class{X}_{-i}'$.
    Additionally we introduce the following set $\class{X}_0 = \class{B}(0, (4q)^{-1})$.
    For every $w \in W \eqdef \{-1, 1\}^m$ we build the distribution $\Prob_w \in \class{P}_W$, such that, the marginal distribution $\Prob_{w, X}$ is independent of $w \in \{-1, 1\}^m$ and the regression vector $(p^w_1(x), \ldots, p^w_K(x))$ is constructed as
    \begin{figure}[!t]
          \centering
          \includegraphics[width=\linewidth]{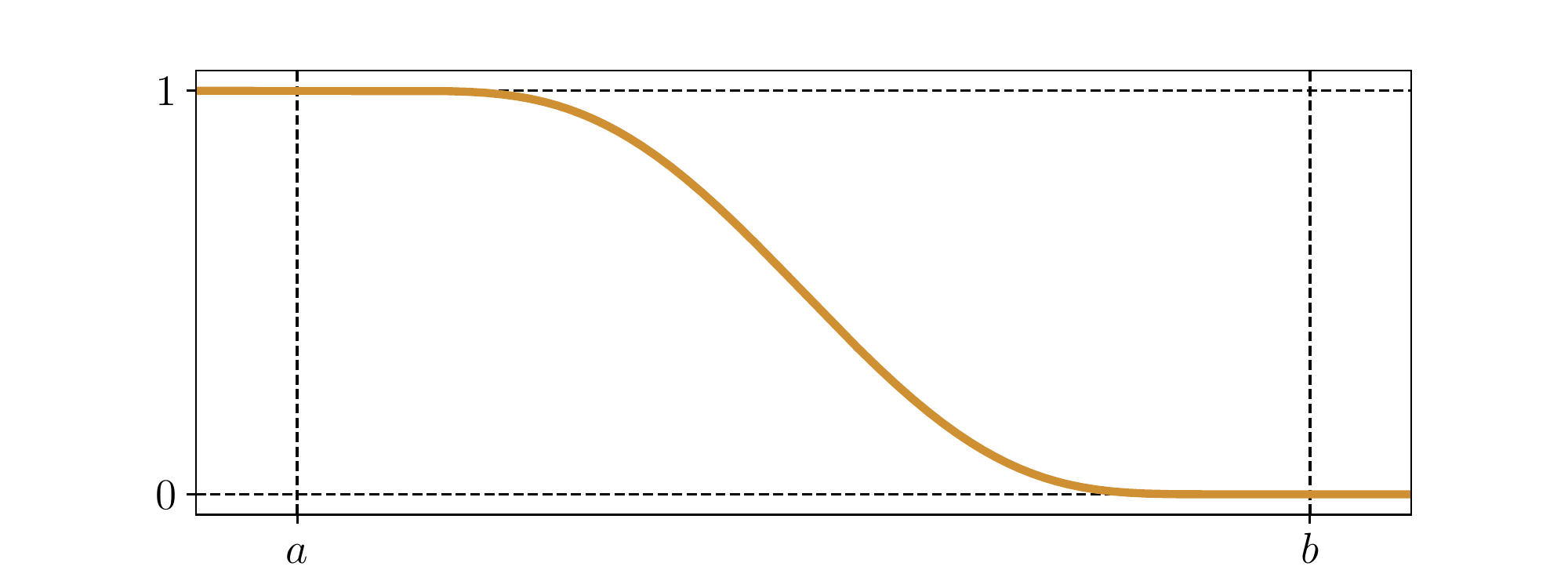}
          \captionof{figure}{Integrated bump: $x \mapsto \frac{\int_x^\infty \psi_{a, b}(t) dt}{\int_a^b \psi_{a, b}(t) dt}$. Importantly, this function is infinitely smooth and is equal to one or zero only outside of the interval $(a, b)$.}
          \label{fig:integrated_psi_function}
    \end{figure}
    \begin{align*}
        p^w_1(x) &= \ldots = p^w_{\beta - 1}(x) = v + \frac{c'}{\beta - 1} + \frac{g(x)}{\beta - 1}\enspace,\\
        p_{\beta}^w(x) &= \begin{cases}
                            v + \phi(x) , &\text{ if } x \in \class{X}_0\\
                            v + w_i \varphi(x - n_q(x)), &\text{ if } x \in \class{X}_i\\
                            v - w_i \varphi(x - n_q(x)), &\text{ if } x \in \class{X}_{-i}\\
                            \frac{1}{K}, &\text{ if } x \in \class{B}(0, \sqrt{d}) \setminus \left(\bigcup_{i = -m, i \neq 0}^m \class{X}_i\right)\\
                            \frac{3v}{2} + g(x), &\text{ if } x \in \bbR^d \setminus \class{B}(0, \sqrt{d} + \rho)\\
                            v + \xi(x), &\text{ if } x \in \class{B}(0, \sqrt{d} + \rho) \setminus \class{B}(0, \sqrt{d})
                         \end{cases}\enspace,\\
        p_{\beta + 1}^w(x) &= \begin{cases}
                            v - \phi(x) , &\text{ if } x \in \class{X}_0\\
                            v - w_i \varphi(x - n_q(x)), &\text{ if } x \in \class{X}_i\\
                            v + w_i \varphi(x - n_q(x)), &\text{ if } x \in \class{X}_{-i}\\
                            v, &\text{ if } x \in \class{B}(0, \sqrt{d}) \setminus \left(\bigcup_{i = -m, i \neq 0}^m \class{X}_i\right)\\
                            \frac{v}{2} - g(x), &\text{ if } x \in \bbR^d \setminus \class{B}(0, \sqrt{d} + \rho)\\
                            v - \xi(x), &\text{ if } x \in \class{B}(0, \sqrt{d} + \rho) \setminus \class{B}(0, \sqrt{d})
                         \end{cases}\enspace,\\
        p^w_{\beta + 2}(x) &= \ldots = p^w_K(x) = v - \frac{c'}{K - \beta - 1} - \frac{g(x)}{K - \beta - 1} \enspace,\\
    \end{align*}
    where $v \in [0, 1]$, $\varphi: \bbR^d \mapsto \bbR_+$, and $\xi: \bbR^d \mapsto \bbR_+$ are to be specified.
    The constants $v, c'$ are set as
    \begin{align*}
        v = \frac{1}{K},\quad c' = \frac{(\beta - 1)(K - \beta - 1)}{K^2}
    \end{align*}
    The function $\xi$ is constructed as
    \begin{align*}
        \xi(x) &= \frac{v}{2}\bar{u}\pare{\frac{\norm{x}_2 - \sqrt{d}}{\rho}},\,
        \bar{u}(x) = 1 - \frac{\int_{x}^{\infty}\psi_{0, 1}(t) dt}{\int_{0}^{1}\psi_{0, 1}(t) dt}\enspace,
    \end{align*}
    the function $\bar{u}$ is infinitely many times differentialble, is equal to zero on $(-\infty, 0]$ and to one on $[1, +\infty)$.
    Figure~\ref{fig:integrated_psi_function} shows the behavior of $1 - \bar{u}$.
    Taking the constant $\rho > 0$ big enough independently of $N, n$ we can ensure that the function $\xi$ is $(\gamma, L)$-H\"older.

    \begin{figure}[!t]
          \centering
          \includegraphics[width=\linewidth]{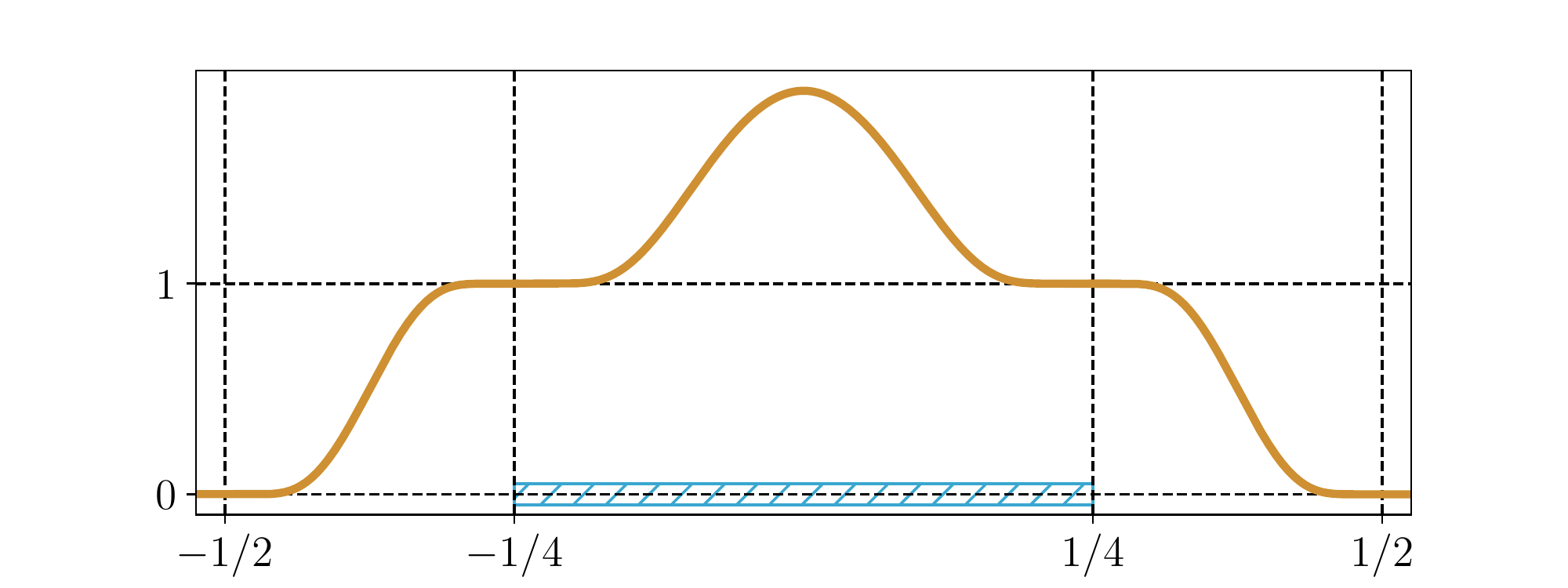}
          \captionof{figure}{The function $x \mapsto u_2(\abs{x}) + \psi_{-\frac{1}{4}, \frac{1}{4}}(x)$. Importantly, this function is infinitely smooth nowhere concentrates at any constant on $(-1/4, 1/4)$. If the marginal measure is supported on the light-blue hatched region the function is lower- and upper-bounded almost surely.}
          \label{fig:fansy_psi_function}
    \end{figure}
    The function $\phi$ is constructed similarly to the previous part of the rate, that is, for $\phi$ we choose
    \begin{align*}
        \phi(x) = C_{\phi} (2q)^{-\gamma} \pare{\frac{\norm{x}}{(2q)^{-1}}}^{2 \ceil{\frac{\gamma}{2}}}\psi_{-1, 1}\pare{\frac{\norm{x}}{(2q)^{-1}}}\enspace,
    \end{align*}
    with $C_{\phi}$ being sufficiently small such that $\phi(\cdot)$ is $(\gamma, L)$-H\"older and upper-bounded by $c'/2\wedge v/4$.
    For the function $\varphi$ we consider the following construction
    \begin{align*}
        \varphi(x) = C_{\varphi}q^{-\gamma}\pare{u_2\pare{\frac{\norm{x}}{q^{-1}}} + \psi_{-\frac{1}{4}, \frac{1}{4}}\pare{\frac{\norm{x}}{q^{-1}}}}\enspace,
    \end{align*}
    where $u_2(\cdot)$ is defined as
    \begin{align*}
        u_2(x) = \frac{\int_{x}^\infty \psi_{\frac{1}{4}, \frac{1}{2}}(t) dt}{\int_{1/4}^{1/2} \psi_{\frac{1}{4}, \frac{1}{2}}(t) dt}\enspace.
    \end{align*}
    Figure~\ref{fig:fansy_psi_function} explains the behavior of this function and helps for better understanding of our results.
    The constant $C_\varphi$ is chosen in such a way that the constructed function $\varphi(\cdot)$ is $(\gamma, L)$-H\"older and and upper-bounded by $c'/2\wedge v/4$.
    Notice that the function $\varphi(x)$ for all $x \in \class{B}(0, (4q)^{-1})$ satisfies
    \begin{align*}
       C_{\varphi}q^{-\gamma}\leq \varphi(x) \leq C_{\varphi}q^{-\gamma}\pare{1 + \psi_{-\frac{1}{4}, \frac{1}{4}}(0)} \leq 2C_{\varphi}q^{-\gamma}\enspace.
    \end{align*}
    Finally, the function $g$ is any $(\gamma, L)$-H\"older function with sufficiently bounded variation which is not concentrated around any constant, for example
    \begin{align*}
        g(x) = C_g\bar{u}\pare{\norm{x}-\sqrt{d} - \rho}\cos\pare{\norm{x} - \sqrt{d} - \rho}\enspace,
    \end{align*}
    For $C_g$ chosen small enough to ensure that it is $(\gamma, L)$-H\"older and has a bounded by $c'/2 \wedge v/4$ variation.

    It remains to define the marginal distribution of the vector $X \in \bbR^d$.
    We select a Euclidean ball in $\bbR^d$ denoted by $A_0$ that has an empty intersection with $\class{B}(0, \sqrt{d} + \rho)$ and whose Lebesgue measure is $\Leb(A_0) = 1 - mq^{-d}$.
    The density $\mu$ of the marginal distribution of $X \in \bbR^d$ is constructed as
    \begin{itemize}
        \item $\mu(x) = \frac{\tau}{\Leb(\class{B}(0, (4q)^{-1}))}$ for every $z \in G_q \cup \{0\}$ and every $x \in \class{B}(z, (4q)^{-1}))$ or $x \in \class{B}(-z, (4q)^{-1}))$,
        \item $\mu(x) = \frac{1 - 2m\tau}{\Leb(A_0)}$ for every $x \in A_0$,
        \item $\mu(x) = 0$ for every other $x \in \bbR^d$,
    \end{itemize}
    for some $\tau$ to be specified.
    Now, we check that the distributions constructed above belong to the set $\class{P}$ for every $w \in W$.
    Namely, we check the following list of assumption
    \begin{itemize}
        \item The functions $p_1^w, \ldots, p_K^w$ are defining some
        regression function for every $w \in W$. That is, for each $x \in \bbR^d$ we have $\sum_{k = 1}^K p_k^w(x) = 1$ and $0 \leq p_k^w(x) \leq 1$,
        \item the functions $p_1^w, \ldots, p_K^w$ are $(\gamma, L)$-H\"older,
        \item the function $G_w(t) \eqdef \sum_{k = 1}^K\int_{\bbR^d}\1_{\{p_k^w(x) \geq t\}} \mu(x) dx$ is continuous,
        \item the threshold $G^{-1}(\beta)$ is equal to $v$ for every $w \in W$,
        \item the marginal distribution satisfies the strong density assumption,
        \item the regression function satisfies $\alpha$-margin assumption.
    \end{itemize}
    {\bf{The regression function is well defined:}}
    to see this, notice that for every $w \in W$ and every $x \in \bbR^d$ we have by construction
    \begin{align*}
        p_{\beta + 1}^w(x) + p_{\beta}^w(x) &= 2v\enspace,\\
        \sum_{k = 1}^{\beta - 1} p^w_k(x) + \sum_{k = \beta + 2}^{K} p^w_k(x) &= (K - 2) v\enspace,
    \end{align*}
    and the combination of both with $v = 1/K$ implies that $\sum_{k = 1}^K p_k^w(x) = 1$.
    Moreover, as long as $\sup_{x \in \class{X}_i} \varphi(x) \leq v/2$ for every $i = -m, \ldots, -1, 1, \ldots, m$ we have for every $x \in \bbR^d$
    \begin{align*}
         0 < v/2 \leq p_{\beta + 1}^w(x) \leq 3v/2 \leq 1,\quad  0 < v/2 \leq p_{\beta}^w(x) \leq 3v/2 \leq 1\enspace,
    \end{align*}
    and by construction of the function $g$ we have for every $k = 1, \ldots, \beta - 1$, every $x \in \bbR^d$ and every $w \in W$
    \begin{align*}
       0 \leq p_k^w(x) \leq v + \frac{3c'}{2(\beta - 1)}\enspace,
    \end{align*}
    due to the choice of $c', v$ we have
    \begin{align*}
         v + \frac{3c'}{2(\beta - 1)} = \frac{1}{K} + \frac{3(K - \beta - 2)}{2K^2} \leq \frac{2}{K} \leq 1\enspace.
     \end{align*}
     Similarly, for every $k = \beta + 2, \ldots, K$, every $x \in \bbR^d$ and every $w \in W$
     \begin{align*}
         v - \frac{3c'}{2(K - \beta - 1)\wedge(\beta - 1)} \leq p_k^w(x) \leq 1\enspace,
     \end{align*}
     and with the choice of $v, c'$ specified above and the constraint $\beta \leq \floor{K / 2}$ we have
     \begin{align*}
         v - \frac{3c'}{2(K - \beta - 1)} = \frac{1}{K} - \frac{3(\beta - 1)}{2K^2} \geq \frac{1}{K} - \frac{3(K/2 - 1)}{2K^2} = \frac{1}{4K} + \frac{3}{2K^2} \geq 0\enspace.
     \end{align*}
     Thus, the construction above defines some regression function for every $w \in W$.

     {\bf{The regression function is $(\gamma, L)$-H\"older:}} this implication follows immediately from the construction of $\varphi, \xi, g$.

    {\bf{Continuity of $G(t)$:}} first let us show that $\int_{\bbR^d} \1_{p_k^w(x) \geq t} \mu(x) dx$ is continuous for every $k \in [K]$.
    For $k = 1, \ldots, \beta-1, \beta+2, \ldots, K$ the continuity follows from the fact that $g$ is not concentrated around any constant.
    For $k = \beta, \beta + 1$ we first write
    \begin{align*}
        \int_{\bbR^d} \1_{p_k^w(x) \geq t} \mu(x) dx
        &=
        \sum_{c \in G_q \cup - G_q}^m\frac{\tau}{\Leb\pare{\class{B}(c, (4q)^{-1})}}\int_{\class{B}(c, (4q)^{-1})} \1_{p_k^w(x) \geq t} dx\\
        &\phantom{=}
        +
        \frac{1 - 2m\tau}{\Leb(A_0)}\int_{A_0}\1_{p_k^w(x) \geq t} dx\enspace,
    \end{align*}
    thus for this choice of $k$ the continuity follows from the fact that $\varphi$ and $g$ are not concentrated around any constant.

    {\bf{Threshold $G^{-1}(\beta) = v$:}} to see this notice that for every $w \in W$,
    \begin{align*}
        \sum_{k = 1}^K\1_{p_k^w(x) \geq v} = \beta,\quad \text{a.e. } \mu\enspace,
    \end{align*}
    and the condition on the threshold follows from the continuity of $G(\cdot)$.
    Besides, the corresponding $\beta$-Oracle sets $\Gamma^*_w$ are given for every $w \in W$ as
    \begin{align*}
        \Gamma^*_w(x) = \begin{cases}
                        \{1, \ldots, \beta - 1, \beta\}, &\quad x\in\class{X}_i,\,w_i = 1,\\
                        \{1, \ldots, \beta - 1, \beta + 1\}, &\quad x\in\class{X}_i,\,w_i = -1,\\
                        \{1, \ldots, \beta - 1, \beta\}, &\quad x\in\class{X}_{-i},\,w_i = -1,\\
                        \{1, \ldots, \beta - 1, \beta + 1\}, &\quad x\in\class{X}_{-i},\,w_i = 1,\\
                        \{1, \ldots, \beta - 1, \beta\}, &\quad x \in \bbR^d \setminus (\bigcup_{i = -m}^m \class{X}_i)
                     \end{cases}
    \end{align*}

    {\bf{The strong density assumption:}} the strong density assumption can be checked following the proof of~\citep[Theorem 3.5]{Audibert_Tsybakov07} where an analogous construction of the marginal distribution was considered.

    {\bf{$\alpha$-margin assumption:}} for all $t \leq t_0 \eqdef v/4$, all $k \in [K] \setminus \{\beta, \beta + 1\}$ and all $w \in W$ we have
    \begin{align*}
        \mu\pare{\abs{p^w_k(X) - v} \leq t} = 0\enspace,
    \end{align*}
    thus for $k \in [K] \setminus \{\beta, \beta + 1\}$ the margin assumption is satisfied.
    It remains to check that the margin assumption is satisfied for $k \in \{\beta, \beta + 1\}$.
    Fix an arbitrary $w \in W$ and $k = \beta$, then for all $t \leq t_0$ we can write
    \begin{align*}
        \mu\Bigg(&\abs{p^w_k(X) - v} \leq t\Bigg)
        =
        \sum_{i = -m}^m\mu\pare{\abs{p^w_k(X) - v} \leq t, X \in \class{X}_i}\\
        &=
        \sum_{i = -m, i\neq 0}^m\mu\pare{\varphi(X - n_q(X)) \leq t, X \in \class{X}_i} + \mu\pare{\phi(X) \leq t, X \in \class{X}_0}\enspace.
    \end{align*}
    We separately upper-bound both terms which appear on the right hand side of the equality.
    \begin{align*}
        \mu&\pare{\phi(X) \leq t, X \in \class{X}_0}
        =
        \frac{\tau}{\Leb(\class{B}(0, (4q)^{-1}))}\int_{\class{B}(0, (4q)^{-1})}\ind{\phi(X) \leq t} dx\\
        &=
        \frac{\tau}{\Leb(\class{B}(0, (4q)^{-1}))}\int_{\class{B}(0, (4q)^{-1})}\ind{C_\phi (2q)^{-\gamma} \pare{\frac{\norm{x}}{(2q)^{-1}}}^{2 \ceil{\frac{\gamma}{2}}}\psi_{-1, 1}\pare{\frac{\norm{x}}{(2q)^{-1}}} \leq t} dx\\
        &=
        \frac{C\tau q^{-d}}{\Leb(\class{B}(0, (4q)^{-1}))}\int_{\class{B}(0, 1/2)}\ind{ \pare{{\norm{x}}}^{2 \ceil{\frac{\gamma}{2}}}\psi_{-1, 1}\pare{{\norm{x}}} \leq C_\phi^{-1}(2q)^{\gamma}t} dx\enspace,
    \end{align*}
    clearly there exists a constant $C$ such that for all $x \in \class{B}(0, 1/2)$ we have
    \begin{align*}
        \psi_{-1, 1}\pare{{\norm{x}}} \geq C\enspace,
    \end{align*}
    Therefore for some constant $C > 0$ we can write
    \begin{align*}
        \mu\pare{\phi(X) \leq t, X \in \class{X}_0}
        &\leq
        \frac{C\tau q^{-d}}{\Leb(\class{B}(0, (4q)^{-1}))}\int_{\class{B}(0, 1/2)}\ind{ {{\norm{x}}} \leq C(q)^{\gamma / 2 \ceil{\frac{\gamma}{2}}}t^{1 / 2 \ceil{\frac{\gamma}{2}}}} dx\\
        &\leq
        \frac{C\tau q^{-d\pare{1 - \gamma / 2 \ceil{\frac{\gamma}{2}}}}}{\Leb(\class{B}(0, (4q)^{-1}))}t^{d / 2 \ceil{\frac{\gamma}{2}}}\enspace,
    \end{align*}
    thanks to the strong density assumption we can write for some $C > 0$
    \begin{align*}
        \mu\pare{\phi(X) \leq t, X \in \class{X}_0} \leq Cq^{-d\pare{1 - \gamma / 2 \ceil{\frac{\gamma}{2}}}}t^{d / 2 \ceil{\frac{\gamma}{2}}}\enspace.
    \end{align*}
    Thus since $1 - \gamma / 2 \ceil{\frac{\gamma}{2}} \geq 0$ and $d / 2 \ceil{\frac{\gamma}{2}} \geq \alpha$ we can write for some $C > 0$
    \begin{align*}
        \mu\pare{\phi(X) \leq t, X \in \class{X}_0} \leq Ct^{\alpha}\enspace.
    \end{align*}
    To finish this part it remains to upper-bound the other term in the margin assumption
    \begin{small}
    \begin{align*}
        \sum_{i = -m, i\neq 0}^m\mu\pare{\varphi(X - n_q(X)) \leq t, X \in \class{X}_i}
        =
        \frac{2m\tau}{\Leb(\class{B}(0, (4q)^{-1}))}\int_{\class{B}(0, (4q)^{-1})}\ind{\varphi(X) \leq t} dx\enspace,
    \end{align*}
     \end{small}
    using the fact that the function $\varphi(x)$ for all $x \in \class{B}(0, (4q)^{-1})$ satisfies
    \begin{align*}
       C_{\varphi}q^{-\gamma}\leq \varphi(x) \leq C_{\varphi}q^{-\gamma}\pare{1 + \psi_{-\frac{1}{4}, \frac{1}{4}}(0)} \leq 2C_{\varphi}q^{-\gamma}\enspace,
    \end{align*}
    we can write for all $t \leq C_{\varphi}q^{-\gamma}$
    \begin{align*}
        \sum_{i = -m, i\neq 0}^m\mu\pare{\varphi(X - n_q(X)) \leq t, X \in \class{X}_i} = 0\enspace,
    \end{align*}
    moreover, for all $t \geq 2C_{\varphi}q^{-\gamma}$ we can write
    \begin{align*}
        \sum_{i = -m, i\neq 0}^m\mu\pare{\varphi(X - n_q(X)) \leq t, X \in \class{X}_i} \leq 2m\tau\enspace,
    \end{align*}
    and finally for $t \in (C_{\varphi}q^{-\gamma}, 2C_{\varphi}q^{-\gamma})$ we can write
    \begin{small}
    \begin{align*}
        \sum_{i = -m, i\neq 0}^m\mu\pare{\varphi(X - n_q(X)) \leq t, X \in \class{X}_i}
        &=
        \frac{2m\tau}{\Leb(\class{B}(0, (4q)^{-1}))}\int_{\class{B}(0, (4q)^{-1})}\ind{\varphi(X) \leq t} dx\\
        &\leq
        \frac{2m\tau}{\Leb(\class{B}(0, (4q)^{-1}))}\int_{\class{B}(0, (4q)^{-1})}\ind{C_{\varphi}q^{-\gamma} \leq t} dx\\
        &=2m\tau\enspace.
    \end{align*}
        \end{small}
    The above implies that for some constant $C > 0$ we have for all $t \leq t_0$
    \begin{align*}
        \sum_{i = -m, i\neq 0}^m\mu\pare{\varphi(X - n_q(X)) \leq t, X \in \class{X}_i}
        &\leq
        2\tau m\ind{t \leq 2C_{\varphi}q^{-\gamma}}\\
        &\leq
        C\tau m q^{\gamma \alpha} t^\alpha\enspace.
    \end{align*}
    Thus the margin assumption is satisfied as long as
    \begin{itemize}
        \item $\tau m = \bigO(q^{-\gamma \alpha})$;
        \item $2 \ceil{\frac{\gamma}{2}} \alpha \leq d$.
    \end{itemize}
    Similarly one can check that the margin assumption is satisfied for $k = \beta + 1$
    {\bf Bound on the KL-divergence:} we are in position to upper-bound the KL divergence between any two hypotheses.
    Fix some $w, w' \in W$, then using the upper bound on $\varphi(\cdot)$ we can write for some $C > 0$
    \begin{align*}
        \KL(\Prob_w, \Prob_{w'})
        &\leq
        2\sum_{i = -m, i\neq 0}^m\mu\pare{\varphi(X - n_q(X))\log\pare{\frac{1 + \varphi(X - n_q(X))}{1 - \varphi(X - n_q(X))}}, X \in \class{X}_i}\\
        &\leq
        Cm\tau q^{-2\gamma}
    \end{align*}

    {\bf How many hypotheses to take:} let us recall the following result which is a version of Varshamov-Gilbert bound~\citep{Gilbert52,Varshamov57}.
    \begin{lemme}
        \label{lem:gilbert_varshamov}
        Let $\delta(w, w')$ denote the Hamming distance between $w, w' \in W$ given by
        \begin{align*}
            \delta(w, w') \eqdef \sum_{i = 1}^m \ind{w_i \neq w'_i}\enspace.
        \end{align*}
        There exists $\class{W} \subset W$ such that for all $w \neq w' \in \class{W}$ we have
        \begin{align*}
            \delta(w, w') \geq \frac{m}{4}\enspace,
        \end{align*}
        and $\log\abs{\class{W}} \geq \frac{m}{8}$.
    \end{lemme}
    Denote $\class{W} \subset W$ the set provided by Lemma~\ref{lem:gilbert_varshamov} and by $\class{P}_{\class{W}}$ the set of distributions $\Prob^w$ with $w \in \class{W}$.
    Taking into account all the above we conclude that $\class{P}_{\class{W}}$ satisfies the assumptions of our result.

    {\bf Lower bound on the Hamming risk (applying Birg\'e's Lemma~\ref{lem:birges_for_lower_bound}):} finally, we are in position to lower bound the hamming risk.
    Recall that we are interested in the following quantity
    \begin{align*}
        \inf_{\hGamma}\sup_{\Prob \in \class{P}} \Exp_{(\data_n, \data_N)}\Exp_{\Prob_X} \abs{\hGamma(X) \triangle \Gamma^*_\beta(X)}\enspace.
    \end{align*}
    The rest of the proof follows standard arguments, which again using the de Finetti notation read as
    \begin{align*}
        \inf_{\hGamma}\sup_{\Prob \in \class{P}} \Exp_{(\data_n, \data_N)}\Exp_{\Prob_X} \abs{\hGamma(X) \triangle \Gamma^*_\beta(X)}
        \geq
        \inf_{\hGamma}\sup_{w \in \class{W}} \mu^{\otimes N}\otimes\Prob_w^{\otimes n}\mu\pare{\abs{\hGamma(X) \triangle \Gamma^*_w(X)}}\enspace.
    \end{align*}
    Denote by $\hat w$ the following minimizer
    \begin{align*}
        \hat w \in \argmin_{w \in \class{W}}\mu\pare{\abs{\hGamma(X) \triangle \Gamma^*_w(X)}}\enspace,
    \end{align*}
    thus if $w \neq \hat w$ we can write using the definition of $\hat w$ and the triangle inequality
    \begin{align*}
        2 \mu\pare{\abs{\hGamma(X) \triangle \Gamma^*_w(X)}}
        &\geq
        \mu\pare{\abs{\hGamma(X) \triangle \Gamma^*_w(X)}} + \mu\pare{\abs{\hGamma(X) \triangle \Gamma^*_{\hat w}(X)}}\\
        &\geq \mu\pare{\abs{\Gamma^*_{\hat w}(X) \triangle \Gamma^*_w(X)}} \geq 2\delta(w, \hat w) \mu(\class{X}_0)\\
        &= 2\delta(w, \hat w) \tau \geq \frac{m\tau}{2}\enspace.
    \end{align*}
    These arguments and Birge's lemma~\ref{lem:birges_for_lower_bound} imply that
    \begin{align*}
        \sup_{\Prob \in \class{P}} &\Exp_{(\data_n, \data_N)}\Exp_{\Prob_X} \abs{\hGamma(X) \triangle \Gamma^*_\beta(X)}
        \geq
        \frac{m\tau}{4}\sup_{w \in \class{W}}\mu^{\otimes N}\otimes\Prob_w^{\otimes n}\pare{w \neq \hat w}\\
        &\geq
        \frac{m\tau}{4}\pare{0.29 \bigvee 1 -\frac{\sum_{w \in \class{W} \setminus \{w'\}}\KL(\mu^{\otimes N}\otimes\Prob_w^{\otimes n}, \mu^{\otimes N}\otimes\Prob_{w'}^{\otimes n})}{\abs{\class{W} -1}\log\abs{\class{W}}}}\enspace.
    \end{align*}
    Since the marginal distribution of the vector $X \in \bbR^d$ is shared among the hypotheses, using the upper-bound on the $\KL$-divergence and the conditions on $\class{W}$ we get for some $C > 0$
    \begin{align*}
        \sup_{\Prob \in \class{P}} \Exp_{(\data_n, \data_N)}\Exp_{\Prob_X} \abs{\hGamma(X) \triangle \Gamma^*_\beta(X)} \geq \frac{m\tau}{4}\pare{1 - Cn\tau q^{-2\gamma}}\enspace.
    \end{align*}
    Finally, let $q = \floor{\bar{C}n^{1/(2\gamma + d)}}$, $\tau = \floor{C'q^{-d}}$ and $m = \floor{C''q^{d - \alpha \gamma}}$ for some $\bar{C}, C', C'' > 0$ small enough we get for some $C > 0$ and $c < 1$
    \begin{align*}
        \sup_{\Prob \in \class{P}} \Exp_{(\data_n, \data_N)}\Exp_{\Prob_X} \abs{\hGamma(X) \triangle \Gamma^*_\beta(X)}
        &\geq
        Cn^{-\alpha \gamma/(2\gamma + d)}\pare{1 - c}\enspace.
    \end{align*}
    One can easily verify that this choice of parameters $\tau, m, q$ is possible as long as $2 \ceil{\frac{\gamma}{2}} \alpha \leq d$ and clearly with our choice we have $\tau m = \bigO(q^{-\alpha \gamma})$.
    As already mentioned the lower bound for the excess risk and the discrepancy follows from Propositions~\ref{prop:distanceComparison} and~\ref{prop:order_risks}.

%%%%%%%%%%%%%%%%%%%%%%%%%%%%%%%%%%%%%%%%%%%%%%%%%%%%%%%%%%%%%%%%%%%%%%%%%%%%%%%
%%%%%%%%%%%%%%%%%%%%%%%%%%%%%%%%%%%%%%%%%%%%%%%%%%%%%%%%%%%%%%%%%%%%%%%%%%%%%%%
\section{Inconsistency of top-\texorpdfstring{$\beta$}{Lg} approach}
\label{app:fail_naive_top_k}
%%%%%%%%%%%%%%%%%%%%%%%%%%%%%%%%%%%%%%%%%%%%%%%%%%%%%%%%%%%%%%%%%%%%%%%%%%%%%%%
%%%%%%%%%%%%%%%%%%%%%%%%%%%%%%%%%%%%%%%%%%%%%%%%%%%%%%%%%%%%%%%%%%%%%%%%%%%%%%%
In this section we prove Proposition~\ref{prop:no_consistency_top_beta}.
The proof builds an explicit construction of a distribution $\Prob$ whose $\beta$-Oracle satisfies $\absin{\Gamma^*_\beta(x)} > \beta$ for all $x$ in some $A \subset \bbR^d$ with $\Prob_X(A) > 0$.
Clearly, if such a distribution exists then there is no estimator in $\hat\Upsilon_\beta$ that would consistently estimate this $\beta$-Oracle.
Let $\beta \in [0, \ldots, \floor{K / 2} - 1]$ be a fixed integer and $K \geq 3$.
For the proof of the theorem we shall construct one distribution $\Prob$ for which none of the estimators with a fixed information can perform well.
We start by specifying the marginal distribution of $X \in \bbR^d$.
We start the construction by specifying the density $\mu$ of the marginal distribution $\Prob_X$.
Define a disk in $\bbR^d$ for some positive $r \leq r'$ as $\class{D}(r, r') = \enscond{x \in \bbR^d}{r \leq \norm{x} \leq r'}$.
First of all fix some parameters $r_1 < r_2 < 2r_2 < r_3 < 2r_3$ which are independent from $n, N$.
% A particular choice of $(r_1, r_2, r_3)$ will be given later in the proof.
The density $\mu$ is supported on $\class{B}(0, r_1) \cup \class{D}(r_2, 2r_2) \cup \class{D}(r_3, 2r_3)$.

Moreover,
\begin{itemize}
    \item $\mu(x) = \frac{\tfrac{\beta}{\beta + 1} - \Leb(\class{B}(0, r_1))}{\Leb\pare{\class{D}(r_2, 2r_2)}}$ for all $x \in \class{D}(r_2, 2r_2)$,
    \item $\mu(x) = \frac{1}{(\beta + 1)\Leb\pare{\class{D}(r_3, 2r_3)}} $ for all $x \in \class{D}(r_3, 2r_3)$,
    \item $\mu(x) = 1$, for all $x \in \class{B}(0, r_1)$,
    \item $\mu(x) = 0$ otherwise,
\end{itemize}
where $r_1 > 0$ is chosen small enough to ensure that $\tfrac{\beta}{\beta + 1} - \Leb(\class{B}(0, r_1)) > 0$.
The regression function $p(\cdot) = (p_1(\cdot), \ldots, p_K(\cdot))^\top$ are defined as
\begin{align*}
    p_{1}(x) = \ldots = p_{\beta + 1}(x) &= \begin{cases}
            \frac{1}{2(\beta + 1)} + C_L\frac{1 - \cos\pare{\frac{2\pi}{r_1}\norm{x}}}{\beta + 1},&\quad x \in \class{B}(0, r_1)\\
            \frac{1}{2(\beta + 1)}+ \frac{g(x)}{\beta + 1}, &\quad x \in \class{D}(r_1, r_2)\\
            \frac{1}{\beta + 1} - C_L\frac{1 - \cos\pare{\frac{2\pi}{r_2}\norm{x}}}{\beta + 1},&\quad x \in \class{D}(r_2, 2r_2)\\
            \frac{1}{\beta + 1} - \frac{\xi(x)}{\beta + 1}, &\quad x\in \class{D}(2r_2, r_3)\\
            \frac{1}{4(\beta + 1)} - C_L\frac{1 - \cos\pare{\frac{2\pi}{r_3}\norm{x}}}{\beta + 1},&\quad x \in \bbR^d \setminus \class{B}(0, r_3)
                                            \end{cases}\enspace,\\
    p_{\beta + 2}(x) = \ldots = p_{K}(x) &= \begin{cases}
            \frac{1}{2(K - \beta - 1)} - C_L\frac{1 - \cos\pare{\frac{2\pi}{r_1}\norm{x}}}{\beta + 1},&\quad x \in \class{B}(0, r_1)\\
            \frac{1}{2(K - \beta - 1)} - \frac{g(x)}{K - \beta - 1}, &\quad x \in \class{D}(r_1, r_2)\\
            C_L\frac{1 - \cos\pare{\frac{2\pi}{r_2}\norm{x}}}{K - \beta - 1},&\quad x \in \class{D}(r_2, 2r_2)\\
            \frac{\xi(x)}{K - \beta - 1}, &\quad x\in \class{D}(2r_2, r_3)\\
            \frac{3}{4(K -\beta - 1)} + C_L\frac{1 - \cos\pare{\frac{2\pi}{r_3}\norm{x}}}{\beta + 1},&\quad x \in \bbR^d \setminus \class{B}(0, r_3)\\
                                            \end{cases}\enspace,
\end{align*}
where  the constant $C_L$ is chosen small enough to ensure that these
functions are $(\gamma, L)$-H\"older and have sufficiently small variation.
Consider an arbitrary infinitely many times differentiable function $v: \bbR \mapsto [0, 1]$ which satisfies $v(x) = 0$ for all $x \leq 0$ and $v(x) = 1$ for all $x \geq 1$.
Then, the functions $g(\cdot)$ and $\xi(\cdot)$ are defined as $g (x) = \frac{1}{2}v\pare{\frac{\norm{x} - r_1}{r_2 - r_1}}$, $\xi(x) = \frac{3}{4}v\pare{\frac{\norm{x} - 2r_2}{r_3 - 2r_2}}$.
The above construction defines a distribution $\Prob$ for which we have
\begin{align*}
    G^{-1}(\beta) &= \frac{1}{2(\beta + 1)}\\
    \Gamma^*_\beta(x) &= \begin{cases}
                            \{1, \ldots, \beta + 1\},&\quad x \in \class{B}(0, r_1) \bigcup \class{D}(r_2, 2r_2)\\
                            \emptyset, &\quad \text{otherwise}
                        \end{cases}\enspace.
\end{align*}
Indeed, let us evaluate the following quantity under the assumption that $\beta \leq \floor{K / 2} - 1$
\begin{align*}
    \sum_{k = 1}^K\int\1_{p_k(x) \geq G^{-1}(\beta)} \mu(x) dx
    &=
    (\beta + 1)\pare{\int_{\class{B}(0, r_1)}\mu(x)dx + \int_{\class{D}(r_2, 2r_2)}\mu(x)dx}\\
    &=
    (\beta + 1) \pare{\Leb\pare{\class{B}(0, r_1)} +\pare{\tfrac{\beta}{\beta + 1} - \Leb(\class{B}(0, r_1))}}\\
    &= \beta\enspace.
\end{align*}
Thus, using this distribution we can write for any classifier $\hGamma \in \hat\Upsilon_{\beta}$ with fixed cardinal
\begin{align*}
    \perf(\hGamma) - &\perf(\Gamma^*_\beta)
    =
    \int_{\bbR^d}\sum_{k = 1}^K\abs{p_k(x) - G^{-1}(\beta)}\1_{k \in \hGamma(x) \triangle \Gamma^*(x)}\mu(x)dx\\
    &\geq
    \int_{\class{D}(r_2, 2r_2)}\abs{\frac{1}{(\beta + 1)} - C_L\frac{1 - \cos\pare{\frac{2\pi}{r_2}\norm{x}}}{\beta + 1} - \frac{1}{2(\beta + 1)}}\mu(x)dx\\
    &=
    \int_{\class{D}(r_2, 2r_2)}\abs{\frac{1}{2(\beta + 1)} - C_L\frac{1 - \cos\pare{\frac{2\pi}{r_2}\norm{x}}}{\beta + 1} }\frac{\tfrac{\beta}{\beta + 1} - \Leb(\class{B}(0, r_1))}{\Leb\pare{\class{D}(r_2, 2r_2)}}dx\enspace,
\end{align*}
where the first inequality follows from the observation that for $x \in \class{D}(r_2, 2r_2)$ there is always at least one label $k$ such that $k \in \hGamma(x) \triangle \Gamma^*(x)$.
Thus, since the constant $C_L$ is chosen to satisfy $2C_L/(\beta + 1) \leq 1 / 4(\beta + 1)$ we have for any $\hGamma \in \hat\Upsilon_{\beta}$
\begin{align*}
    \perf(\hGamma) - \perf(\Gamma^*_\beta) \geq {\frac{\tfrac{\beta}{\beta + 1} - \Leb(\class{B}(0, r_1))}{4(\beta + 1)} }\enspace,
\end{align*}
If $r_1$ is such that $\Leb(\class{B}(0, r_1)) \leq \tfrac{\beta}{2(\beta + 1)}$ we get
\begin{align*}
    \perf(\hGamma) - \perf(\Gamma^*_\beta) \geq \frac{\beta}{8 (\beta + 1)^2},\quad \text{almost surely}\enspace.
\end{align*}
By construction, the regression vector is $(\gamma, L)$-H\"older and the density is lower- and upper-bounded by some positive constants.
Hence, it remains to check that the constructed distribution satisfies the $\alpha$-margin assumption. This can be achieved by an appropriate choice of $r_1$. Indeed, on the sets $\class{D}(r_2, 2r_2) \cup \class{D}(r_3, 2r_3)$ there is a ``corridor'' of constant size between the regression functions and the threshold $G^{-1}(\beta)$.
The threshold $G^{-1}(\beta)$ is only approached by the regression function on the set $\class{B}(0, r_1)$.
As all the parameters in our construction are independent from $n, N \in \bbN$ we can find a value $r_1$ being small enough so that the $\alpha$-margin assumption is verified for a fixed $\alpha > 0$.

\end{document}